\documentclass[a4paper,12pt]{article}

\usepackage[latin1]{inputenc} 
\usepackage[T1]{fontenc}      
\usepackage{geometry}         
\usepackage[francais, english]{babel}  
\usepackage{amsmath}
\usepackage{amsfonts}
\usepackage{amssymb}
\usepackage{fancyhdr}
\usepackage{url}
\usepackage{mathrsfs}
\usepackage{euscript}
\usepackage{enumerate}
\usepackage{cite}
\usepackage[all]{xy}
\usepackage{mathtools}

\newtheorem{df}{Definition}[section]

\newtheorem{fact}{Fact}
\newtheorem{prop}[df]{Proposition}
\newtheorem{thm}[df]{Theorem}

\newtheorem{lem}[df]{Lemma}

\newcommand{\prf}{\noindent \textit{Proof}}
\newtheorem{rmk}[df]{Remark}

\newcommand{\dbar}{\overline{\partial}}

\newcommand{\R}{\mathbb{R}}

\newcommand{\C}{\mathbb{C}}
\newcommand{\Z}{\mathbb{Z}}
\newcommand{\N}{\mathbb{N}}

\newcommand{\ric}{\operatorname{Ric}}
\newcommand{\tr}{\operatorname{tr}}
\renewcommand{\det}{\operatorname{det}} 
\newcommand{\euc}{\operatorname{euc}}
\newcommand{\mxd}{\operatorname{mxd}}

\newcommand{\vol}{\operatorname{vol}}

\newcommand{\vareps}{\varepsilon}
\newcommand{\eps}{\epsilon}

\newcommand{\vl}{\operatorname{Vol}}

\newcommand{\inj}{\operatorname{inj}}

\newcommand{\Sp}{\operatorname{Sp}}

\newcommand{\id}{\operatorname{id}}

\newcommand{\riem}{\operatorname{Rm}}
\newcommand{\e}{\mathbf{e}}
\newcommand{\f}{\mathbf{f}}
\newcommand{\cc}{\mathbf{c}}

\newcommand{\diff}{\mathbf{Diff}}
\newcommand{\met}{\mathbf{Met}}
\newcommand{\SU}{\operatorname{SU}}
\newcommand{\SO}{\operatorname{SO}}
\newcommand{\diag}{\operatorname{diag}}
\newcommand{\cqfd}{ \hfill $\square$ }
\newcommand{\sameauthor}{\leavevmode\vrule height 2pt depth -1.6pt width 40pt}

\title{\large FROM ALE TO ALF GRAVITATIONAL INSTANTONS. II}

\vspace{5pt}
\author{\normalsize\textsc{Hugues AUVRAY
          \footnote{This work was started while the author was visiting 
                    the University of Edinburgh during Spring 2012, 
                    supported by an FSMP grant, 
                    and completed during the author's stay at the MPIM Bonn as an EPDI post-doc visitor in 2013.}}}
\date{}

\begin{document}

\makeatletter
\renewcommand%
   {\section}%
   {%
   \@startsection{section}%
      {1}
      {0mm}%
      {\baselineskip}
      {0.5\baselineskip}%
      {\sc\large\centering}
   }%
\makeatother

\makeatletter
\renewcommand%
   {\subsubsection}%
   {%
   \@startsection{subsubsection}%
      {1}
      {0mm}%
      {1.25\baselineskip}
      {0.25\baselineskip}%
      {\sf\normalsize}
   }%
\makeatother

\makeatletter
\renewcommand%
   {\paragraph}%
   {%
   \@startsection{paragraph}%
      {4}%
      {0mm}%
      {0mm}%
      {-5pt}%
      {\it\normalsize}%
   }%
\makeatother

\maketitle

\renewcommand{\abstractname}{Abstract}
 \begin{abstract}
  This paper is the sequel of our previous article \textit{From ALE to ALF gravitational instantons}, 
  where we constructed ALF hyperkähler metrics on \textit{minimal resolutions} of $\C^2/\mathcal{D}_k$, 
  with $\mathcal{D}_k$ the binary dihedral group of order $4k$, $k\geq 2$. 
  In the present article we generalize the construction to \textit{smooth deformations} of this Kleinian singularity, 
  with help of the computation of the asymptotics of the ALE gravitational instantons. 
 \end{abstract}

 \begin{center}
   \rule{3cm}{0.25pt}
 \end{center}

~

\section*{Introduction}

A central question in Riemannian geometry in real dimension 4 is that of the comprehension of \textit{non-compact, complete, Ricci-flat manifolds}. 
These spaces indeed arise naturally in differential geometry as limiting spaces, 
after rescaling, of families of compact Einstein 4-manifolds; 
their knowledge might thus be of crucial interest for the study of Einstein compact 4-manifolds and especially of sequences of such spaces.   

The dimension 4 moreover allows one to specify the question to \textit{Ricci-flat Kähler}, 
or even \textit{hyperkähler}, non-compact, complete manifolds. 
If one adds furthermore a decay condition on the Riemannian curvature tensor, this leads to:
 \begin{df}[Gravitational instantons] 
  A gravitational instanton is a non-compact, complete, hyperkähler manifold $(X,g, I,J,K)$ of real 
  dimension 4, whose Riemannian curvature satisfies the following $L^2$ integral condition:  
  $\int_X|\riem^{g}|^2\rho\vol^g$ is finite, 
  where $\rho$ is a "ball volume growth ratio" function. 
 \end{df}
By a ball volume growth ratio function, we merely mean a function comparing the growth of balls 
in $X$ to the that of the flat model $(\R^4,\e)$, 
with $\e$ the standard euclidean metric on $\R^4$;  
more concretely, fixing $x\in X$ and setting $r=d_g(x,\cdot)$, one can take $\rho=\frac{r^4}{1+\vl(B_X(x,r))}$.

Let us add here that besides this differential-geometric definition, 
gravitational instantons also appear as fundamental objects in theoretical physics, 
in fields such as Quantum Gravity \cite{haw} or String Theory \cite{ch-hi, ch-ka}.

~

Recall that the hyperkähler condition implies the Ricci-flatness. 
Thus, the fundamental Bishop-Gromov theorem \cite{glp} implies that on gravitational instantons, the function $\rho$ mentioned above is
at least positively bounded below -- in other words, 
the growth of the ball volume is at most euclidean. 
If it is also bounded above, one deals with \textit{Asymptotically Locally Euclidean} instantons, 
or \textit{ALE instantons} for short. 
These hyperkähler manifolds are completely classified, after \cite{bkn} and \cite{kro2} 
(notice also the recent extension \cite{suv} to the Kähler Ricci-flat case); 
their comprehension actually goes deeper, 
since the classification corresponds to an exhaustive construction by Kronheimer \cite{kro1}, 
which is why we shall also often refer to these spaces as \textit{Kronheimer's instantons}. 
In a nutshell, the hyperkähler structures of these spaces are asymptotic to that of a quotient $\R^4/\Gamma$, 
with $\Gamma$ a finite subgroup of $\SU(2)=\Sp(1)$; 
when moreover $\Gamma$ is fixed, these spaces are all diffeomorphic to the \textit{minimal resolution of the Kleinian singularity} $\C^2/\Gamma$.  
Here $\C^2$ stands for $(\R^4,I_1)$ with $I_1$ the standard complex structure given by the coordinates $z_1=x_1+ix_2$, $z_2=x_3+ix_4$; 
we fix this notation, as well as that of $I_2$ and $I_3$ for the other two standard complex structures on $\R^4\cong\mathbb{H}$, 
given respectively by the coordinates $(x_1+ix_3,x_4+ix_2)$ and $(x_1+ix_4,x_2+ix_3)$. 
 
Now, still on gravitational instantons, a result by Minerbe \cite{min1} states the following quantification: 
\textit{if the ball volume growth is less than euclidean, i.e quartic, it is at most cubic}; 
in other words, if the comparison function $\rho$ is not bounded, it grows at most like $r$. 
One then speaks about \textit{Asymptotically Locally Flat}, or \textit{ALF}, 
gravitational instantons. 
The prototype of such metrics is (are) the Taub-NUT metric(s), living on $\R^4$. 
Roughly speaking, half of these spaces are classified, by Minerbe again \cite{min3}; 
their geometry at infinity is that of a circle fibration over $\R^3$, 
and they are explicitly described by the so-called \textit{Gibbons-Hawking ansatz}. 

~

\noindent
\textit{Results.} 
The only other possibility for the asymptotic geometry of the ALF gravitational instantons is that of a circle fibration over $\R^3/\pm$ \cite{min2}.   
This we illustrate by the following, which is one of the two main results of this paper:   
 \begin{thm} \label{thm_ALFintro}
  Let $\big(X,g, I_1^X,I_2^X,I_3^X)$ be an ALE gravitational instanton modelled on $\R^4/\mathcal{D}_k$, 
  in the sense that $X$ minus a compact subset is diffeomorphic to $\R^4/\mathcal{D}_k$ minus a ball, 
  and that the hyperkähler structure of $X$ is asymptotic to that of $\R^4/\mathcal{D}_k$ via the diffeomorphism in play; 
  here $\mathcal{D}_k$ is the binary dihedral group of order $4k$, $k\geq2$. 
  Then there exists on $X$ a family of hyperkähler structures $(g_m, J_{1,m}^X,J_{2,m}^X,J_{3,m}^X)_{m\in(0,\infty)}$, 
  such that:
   \begin{enumerate}
    \vspace{-3pt}
    \item the diffeomorphism above can be chosen so that $g_m$ is asymptotic to the Taub-NUT metric $\f_m$; 
    \item the Kähler classes associated to this hyperkähler structure are the same as those of the ALE hyperkähler structure, 
          and moreover $\vol^{g_m}=\vol^{g_X}$; 
    \item the curvature tensor $\riem^{g_m}$ has cubic decay. 
   \end{enumerate}
 \end{thm}

The metric $\f_m$ of this statement is the Taub-NUT metric \textit{of parameter} $m$, 
the squared inverse of which shall be interpreted modulo a multiplicative constant as the length at infinity of the fibres of a circle action $\f_m$. 
The asymptotics between the ALF metric $g_m$ and $\f_m$ are as follows: 
if $R$ denotes a distance function for $\f_m$, then $(g_m-\f_m)$ and $\nabla^{\f_m}(g_m-\f_m)$ are $O(R^{-2+\eps})$ for an arbitrarily small $\eps>0$. 
Of course $\f_m$ is also invariant under $\mathcal{D}_k$, so that it makes perfect sense on $\R^4/\mathcal{D}_k$. 

~

Before discussing in more details on how Theorem \ref{thm_ALFintro} is proved, 
we shall underline the following: 
its proof heavily relies on the computation of the asymptotics of the ALE instantons modelled on $\R^4/\mathcal{D}_k$. 
More precisely, Kronheimer's construction of these spaces allows one to write these asymptotics under the shape of a power series, 
the main term of which is the euclidean model, 
and this actually holds for any finite subgroup $\Gamma$ of $\SU(2)$ alluded above. 
Our second main result deals with the first non-vanishing terms of those expansions:
 \begin{thm} \label{thm_ALEintro}
  Let $\big(X,g, I_1^X,I_2^X,I_3^X)$ be an ALE gravitational instanton modelled on $\R^4/\Gamma$. 
  Then one can choose a diffeomorphism $\Phi$ between $X$ minus a compact subset and $\R^4$ minus a ball 
  such that:
   \begin{enumerate}
    \vspace{-3pt}
    \item $\Phi_*g_X-\e=h_X+O(r^{-6})$, $\Phi_*I_1^X-I_1=\iota_1^X+O(r^{-6})$ 
          and if $\omega_1^X=g_X(I_1^X\cdot, \cdot)$ and $\omega_1^{\e}=\e(I_1\cdot, \cdot)$, 
          then $\Phi_*\omega_1^X-\omega_1^{\e}=\varpi_1^X+O(r^{-6})$, 
          where $h_X$, $\iota_1^X$ and $\varpi_1^X$ admit  explicit formulas and are $O(r^{-4})$; 
          for instance $\varpi_1^X=-\sum_{j=1}^3c_j(X)dd^c_{I_j}(r^{-2})$ for some explicit constants $c_j(X)$. 
    \item when $\Gamma$ is not a cyclic subgroup of $\SU(2)$, the $O(r^{-6})$ of the previous point can be replaced by $O(r^{-8})$. 
   \end{enumerate}                                                                                             
 \end{thm}
 Here the $O$ are understood in an asymptotically euclidean context: 
 $\vareps$ is $O(r^{-a})$ if for any $\ell\geq0$, $\big|(\nabla^{\e})^{\ell}\vareps\big|_{\e}=O(r^{-a-\ell})$ near infinity. 

~

 \noindent
 \textit{Organization of the article.} 
 The paper is organized as follows. 
 It is divided into two parts, reflecting the dichotomy of Theorem \ref{thm_ALFintro} and \ref{thm_ALEintro}. 
 Part 1 is devoted to the proof of Theorem \ref{thm_ALFintro}. 
 We first draw in section \ref{section_prgrm} a detailed program of construction of our hyperkähler ALF metrics, 
 leading us to the expected result (Theorem \ref{thm_gRF}). 
 In section \ref{section_TN} are recalled essential facts on the Taub-NUT metric, seen as a Kähler metric on $\C^2$. 
 The construction itself occupies section \ref{section_gluing} and \ref{section_corrections}; 
 roughly speaking, it consists into a gluing of the Taub-NUT metric with the ALE metric of some ALE instanton, 
 which we subsequently correct into a Ricci-flat metric thanks to an appropriate Calabi-Yau theorem adapted to ALF geometry. 
 The concluding section \ref{section_lemmas}, 
 which is mainly computational, 
 deals with the proof of two technical lemmas useful to our construction. 

 In Part 2, which mostly does not depend on Part 1, after recalling some basic facts about Kronheimer's construction of ALE instantons, 
 we state Theorem \ref{thm_ALE}, which is a specified version of Theorem \ref{thm_ALEintro} 
 -- in particular we give the promised explicit formulas (section \ref{section_thmALE}). 
 We give further details on Kronheimer's construction and classification in section \ref{section_kroconst}, 
 where we also fix the diffeomorphism of Theorem \ref{thm_ALEintro}. 
 Then we compute the tensors $h_X$, $\iota_1^X$ and $\varpi_1^X$ in section \ref{section_hzeta}; 
 using similar techniques, we show in the following section \ref{section_ord3term} that the precision of the asymptotics is automatically improved 
 when $\Gamma$ is binary dihedral, tetrahedral, octahedral or icosahedral. 
 We develop in the last section \ref{section_cplxstr} a few informal digressive considerations on the approximation of complex structures 
 of certain ALE instantons by the standard $I_1$, 
 relied on links observed in the construction of Part 1.

~

 \noindent
 \textit{Comments.} 
  Let us start with a few words on previous constructions of ALF dihedral gravitational instantons. 
  Such objects are known to exist since the works \cite{ch-hi, ch-ka}, 
  in which twistor methods are employed. 
  The produced hyperkähler spaces are unfortunately not much explicit, 
  which partly motivates our construction, although we hope that actually, 
  \textit{both constructions produce the same families of ALF hyperkähler metrics}, 
  and more precisely that \textit{these constructions are exhaustive, in the sense the any ALF dihedral gravitational instanton fits into the produced examples} 
  (up to a tri-holomorphic isometry): 
  this folklore conjecture is the analogue of the classification of \cite{min3}. 
  Let us recall also our previous construction \cite{auv1}, 
  in which we restricted ourselves to \textit{resolutions of Kleinian dihedral singularities} for technical reasons. 
  We actually here follow the same guideline, bypassing those difficulties. 
  We would however like to emphasise their non-trivial character: 
  we indeed need Theorem \ref{thm_ALEintro} to tackle them, 
  and some rather important technical adjustments are still necessary to handle our construction (see e.g. Proposition \ref{prop_approx_omega1Y}). 
  Consequently, our approach in Part 1 is to focus on these adjustments, 
  recalling only the minimal material from \cite{auv1} to get to the envisaged end, i.e. Theorem \ref{thm_ALFintro}; 
  in particular, we only state and use an ALF Calabi-Yau type theorem which was a major result of that paper. 
  
  More closely to the statement of Theorem \ref{thm_ALFintro}, 
  notice that it is \textit{not of a perturbative nature}: 
  this corresponds to taking the parameter $m$ in the whole range $(0,+\infty)$. 
  The price to pay is somehow that so far, we do not control what happens when $m$ goes to 0.  
  We conjecture that the ALF hyperkähler structure converges in $C^{\infty}_{loc}$-topology to the original ALE one, as is the case on $\C^2$; 
  this question will be handled in a future article. 
  
  Now in Theorem \ref{thm_ALEintro}, the existence of the first order variation terms $h_X$, $\iota_1^X$ and $\varpi_1^X$ is of course not new,  
  as it is even known from Kronheimer's construction \cite{kro1} that the ALE hyperkähler structures 
  admit an expansion in terms of similar tensors. 
  What is new though is their explicit determination, 
  which we could only find in the simplest case that is the Eguchi-Hanson space (see e.g \cite[p.153]{joy}), 
  i.e. when $\Gamma=\mathcal{A}_2=\{\pm \id_{\C^2}\}$. 
  Notice at this point that as suggested by the statement of Theorem \ref{thm_ALEintro}, 
  the shapes of $h_X$, $\iota_1^X$ and $\varpi_1^X$ follow a general pattern which is only slightly affected by the group $\Gamma$; 
  up to a multiplicative constant, we can indeed compute them on the explicit Eguchi-Hanson example. 
  We think moreover that Theorem \ref{thm_ALEintro} is of further interest, 
  and that among others, the order of precision it brings could be useful in more general gluing constructions.

~

 \noindent
 \textit{Acknowledgements.} 
  I would like to thank my former PhD advisor O. Biquard for his wide advice, 
  and particularly for his insight into Theorem \ref{thm_ALEintro}. 
  I am also thankful to FSMP for its support during my visit in the University of Edinburgh, 
  and to the Max-Planck-Institut für Mathematik Bonn for its very stimulating environment. 

~

\section{Construction of ALF hyperkähler metrics}
  \subsection{Strategy of construction}   \label{section_prgrm}

\noindent
\textit{Outline of the strategy.} 
As described in \cite{auv1} and as we shall see in next section, 
one can describe the Taub-NUT metric on $\R^4$ as a $\mathcal{D}_k$-invariant hyperkähler metric with volume the standard $\Omega_{\e}$, 
Kähler for the standard complex structure $I_1$, 
and compute a somehow explicit potential, $\varphi$ say, for it. 

Now, given one of Kronheimer ALE gravitational instantons $\big(X, g_X, I_1^X, I_2^X, I_3^X\big)$ modelled on $\R^4/\mathcal{D}_k$, 
we have a diffeomorphism $\Phi_X$ between infinities of $X$ and $\R^4/\mathcal{D}_k$ such that ${\Phi_X}_*g_X$ is asymptotic to $\e$, 
and ${\Phi_X}_*I_1^X$ is asymptotic to $I_1$. 
It is this way quite natural to try to take $dI_1^Xd(\Phi_X^*\varphi)$ glued with $g_X$ as an ALF metric on $X$, 
before we correct it into a hyperkähler metric. 
This naive idea does work when $X$ is \textit{a minimal resolution of $(\C^2/\mathcal{D}_k, I_1)$ and $\Phi_X$ the associated map}: 
this is the purpose of \cite{auv1}. 
However this fails in the general case, where $X$ is a \textit{deformation of $(\C^2/\mathcal{D}_k, I_1)$}, 
without further precautions: 
the size of the Taub-NUT potential $\varphi$, roughly of order $r^4$ \textit{as well as its euclidean derivatives}, 
together with the error term ${\Phi_X}_*I_1^X - I_1$ on the complex structure, 
even make wrong the assertion that the rough candidate $dI_1^Xd(\Phi_X^*\varphi)$ is positive 
-- in the sense that $dI_1^Xd(\Phi_X^*\varphi)(\cdot,I_1^X\cdot)$ is a metric -- near the infinity of $X$. 
 
Fortunately, up to choosing a different complex structure on $X$ to work with, 
we can make the appropriate corrections on $\varphi$ so as to get a good enough ALF metric on $X$ to start with, 
and then run the same machinery as in \cite{auv1} up to minor but yet technical adjustments, so as to end up with Theorem \ref{thm_ALFintro}.  

~

\noindent
\textit{Detailed strategy, and involvements of the asymptotics of Kronheimer's instantons.} 
We shall now be more specific about the different steps involved in the program we are following throughout this part. 
\begin{enumerate}
  \vspace{-3pt}
   \item Let act $\SO(3)$ on the complex structure of $X$ as follows: 
         for $A=(a_{j\ell})\in \SO(3)$, define the triplet $(AI^X)_{\cdot}$ as 
          \begin{equation*}
           (AI^X)_{\cdot} = \big((AI^X)_{j}\big)_{j=1,2,3}= \big(a_{j1}I_1^X+a_{j2}I_2^X+a_{j3}I_3^X\big)_{j=1,2,3};
          \end{equation*}
         then $\big(X, g_X, (AI^X)_{1}, (AI^X)_{2},(AI^X)_{3}\big)$ is again hyperkähler, 
         and therefore an ALE gravitational instanton modelled on $\R^4/\mathcal{D}_k$. 
   \item With the model $\R^4/\mathcal{D}_k$ at infinity fixed, 
         Kronheimer's instantons are parame\-trised \cite{kro1} by a triplet $\zeta=(\zeta_1,\zeta_2,\zeta_3)\in\mathfrak{h}\otimes\R^3-D$, 
         where $\mathfrak{h}$ is a $(k+2)$-dimensional real vector space endowed with some scalar product $\langle\cdot,\cdot\rangle$, 
         and $D$ is a finite union of spaces $H\otimes\R^3$ with $H$ a hyperplane in $\mathfrak{h}$ 
         (as notation suggests, $\mathfrak{h}$ is a Lie algebra; 
         we will be more specific about its interpretation in part \ref{part_ALE}).  
         This parametrisation is compatible with the $\SO(3)$-action of Point 1. in the sense that if $\zeta$ is the parameter associated to 
         $\big(X, g_X, I_1^X, I_2^X, I_3^X\big)$, and if $\big(Y, g_Y, I_1^Y, I_2^Y, I_3^Y\big)$ is the instanton associated to $A\zeta$, 
         defined by:
          \begin{equation}   \label{eqn_Azeta}
           A\zeta = \big((A\zeta)_{j}\big)_{j=1,2,3}= \big(a_{j1}\zeta_1+a_{j2}\zeta_2+a_{j3}\zeta_3\big)_{j=1,2,3}, 
          \end{equation}
         then there exists an isometry which is moreover tri-holomorphic between $\big(X, g_X, (AI^X)_{1}, (AI^X)_{2},(AI^X)_{3}\big)$ 
         and $\big(Y, g_Y, I_1^Y, I_2^Y, I_3^Y\big)$:  
         this is Lemma \ref{lem_zetaAzeta}, stated and proved in Part 2.
         Defined this way, $A\zeta$ is of course still in $\mathfrak{h}\otimes\R^3-D$; 
         otherwise $A\zeta\in H\otimes\R^3$ for one of the hyperplanes $H$ mentioned above, 
         and thus $\zeta=A^t(A\zeta)\in H\otimes\R^3$, which would be absurd.
     
   \item 
         In general, one can take the diffeomorphism $\Phi_X$ between infinities of $X$ and $\R^4/\mathcal{D}_k$ so that 
         ${\Phi_X}_*I_1^X - I_1=O(r^{-4})$ with according decay on derivatives, which is not good enough for our construction, 
         see \cite[p.17-18]{auv1}.  
         We can nonetheless improve the precision thanks to the following two lemmas:
          \begin{lem}  \label{lem_beth}
           If $\xi\in\mathfrak{h}\otimes\R^3-D$ is such that $|\xi_2|^2-|\xi_3|^3=\langle\xi_2,\xi_3\rangle= 0$, 
           and $\big(Y, g_Y, I_1^Y, I_2^Y, I_3^Y\big)$ is the associated ALE instanton, then
           one can choose $\Phi_Y$ such that there exists a diffeomorphism $\beth=\beth_{\xi}$ of $\R^4$ commuting with the action of $\mathcal{D}_k$, 
           and such that: 
            \begin{equation*}
             \Big|(\nabla_{\e})^{\ell}\big({\Phi_Y}_*I_1^Y - \beth^*I_1\big)\Big|_{\e}=O(r^{-8-\ell}) \quad \text{for all }\ell\geq0. 
            \end{equation*}
           Moreover, the shape of $\beth$ is given by: $\beth(z_1,z_2)=\big(1+\frac{a}{\kappa+r^4}\big)(z_1,z_2)$, 
           where $\kappa,a\in\R$, and $(z_1,z_2)$ are the standard complex coordinates on $(\C^2,I_1)$, 
           and $\big|(\nabla_{\e})^{\ell}\big(\Omega_{\e} - \beth^*\Omega_{\e}\big)\big|_{\e}=O(r^{-8-\ell})$ for all $\ell\geq0$. 
          \end{lem}
          \begin{lem} \label{lem_AZ}
           For any $\zeta\in \mathfrak{h}\otimes\R^3$, 
           there exists $A\in SO(3)$ such that $|(A\zeta)_2|^2-|(A\zeta)_3|^2=\langle(A\zeta)_2,(A\zeta)_3\rangle=0$.
          \end{lem}
         Lemma \ref{lem_beth}, \textit{which relies on our analysis of the asymptotics of Kronheimer's instantons}, 
         is proved in section \ref{section_gluing}, assuming a general statement for this asymptotics that is seen in Part 2; 
         Lemma \ref{lem_AZ}, which is elementary, is proved at the end of this section. 
 
   \item We choose $A$ as in Lemma \ref{lem_AZ}, consider the instanton $Y$ associated to $\xi=A\zeta\in\mathfrak{h}\otimes\R^3-D$, 
         and apply the original program, i.e. the gluing of the Kähler forms and the corrections of the prototypical ALF metric into a hyperkähler metric, 
         to the potential $\varphi^{\flat}:=\beth^*\varphi$ instead of $\varphi$.  
         Thanks to the better coincidence of the complex structures, the rough candidate $dI_1^Xd(\Phi_X^*\varphi^{\flat})$ is now positive at infinity, 
         and actually also rather close to $\f^{\flat}:=\beth^*\f$, with $\f$ the Taub-NUT metric on $\R^4$, Kähler for $I_1$. 
         
         We should moreover specify here that the gluing requires a precise description of the Kähler form $\omega_1^Y:= g_Y\big(I_1^Y\cdot,\cdot\big)$, 
         \textit{which is again part of the analysis of the asymptotics of Kronheimer's ALE instantons}.
          
         We get this way after corrections a Ricci-flat, actually hyperkähler, 
         manifold $\big(Y,g'_Y, I^Y_1, J^Y_2, J^Y_3\big)$, 
         with ${\Phi_Y}_*g'_Y$ asymptotic to $\f^{\flat}$, 
         and $\big[g'_Y(I^Y_1\cdot,\cdot)\big]=[g_Y(I^Y_1\cdot,\cdot)]$; 
         the construction gives moreover $\big[g'_Y(J^Y_j\cdot,\cdot)\big]=[g_Y(I^Y_j\cdot,\cdot)]$, $j=2,3$.  

   \item We let $A^t=A^{-1}$ act back on the previous data to come back to $X$, 
         and end up with a hyperkähler manifold $\big(X,g'_X, J^X_1, J^X_2, J^X_3\big)$, with $\big[g'_X(J^X_j\cdot,\cdot)\big]=[g_X(I^X_j\cdot,\cdot)]$, $j=1,2,3$, 
         and ${\Phi_X}_*g'_X$ asymptotic to $\f^{\flat}$, 
         provided that $\Phi_X$ is the composite of $\Phi_Y$ and the tri-holomorphic isometry of Point 2 .  
  \end{enumerate}

We shall also add that we can play on the metric $\f$ in this construction. 
Indeed, $\f$ is invariant under some fixed circle action on $\R^4$, 
and the length for $\f$ of the fibres of this action tends to some constant $L>0$ at infinity.  
We can make this length vary in the whole $(0,\infty)$ and keep the same volume form for $\f$; 
given $m\in (0,\infty)$ that we call the "\textit{mass parameter}", 
we then denote by $\f_m$ the Taub-NUT metric giving length $L(m)=\pi\sqrt{\frac{2}{m}}$ to the fibres at infinity, and of volume form $\Omega_{\e}$ 
(the choice of the parameter $m$ instead of $L$ will become clear below).  

We can then sum up our construction by the following statement, which is the main result of this part:
 \begin{thm}  \label{thm_gRF}
  Consider an ALE gravitational instanton $\big(X, g_X, I_1^X, I_2^X, I_3^X\big)$ modelled on $\R^4/\mathcal{D}_k$. 
  Then there exists a one-parameter family $\big(g'_{X,m}, J_{1,m}^X, J_{2,m}^X, J_{3,m}^X\big)$ 
  of smooth hyperkähler metrics on $X$ such that:
   \begin{itemize}
    \vspace{-3pt}
    \item for any $m$, one has the equality $\big[g'_{X,m}\big( J_{1,m}^X\cdot,\cdot\big)\big]=\big[g_X\big(I_1^X\cdot,\cdot\big)\big]$ of Kähler classes 
          for $j=1,2,3$; 
    \item for any $m$, $g'_{X,m}$ and $g_X$ have the same volume form; 
    \item when $m$ is fixed, $g'_{X,m}$ is ALF in the sense that one has the asymptotics 
           \begin{equation*}
            \Big|(\nabla^{\f_m^{\flat}})^{\ell}\big({\Phi_X}_{*}g'_{X,m} - \f_m^{\flat}\big) \Big|_{\f_m^{\flat}} = O\big(R^{-1-\delta}\big), \qquad \ell=0,1,
           \end{equation*}
          for any $\delta\in (0,1)$ and that $\riem^{g'_{X,m}}$ has cubic decay at infinity.  
   \end{itemize} 
     Here $R$ is a distance function for $\f_m^{\flat}$, 
     and $\Phi_X$ is an ALE diffeomorphism between infinities of $X$ and $\R^4/\mathcal{D}_k$, 
     in the sense that $|{\Phi_X}_*g_X-\e|_{\e}$, $|{\Phi_X}_*I_j^X-I_j|_{\e}=O(r^{-4})$, 
     with according decay on derivatives.  
 \end{thm}

  In this statement, $\f_m^{\flat}=\beth^*\f_m$, where $\beth=\beth_{A\zeta}$ is given by Lemma \ref{lem_beth}, 
  $\zeta\in \mathfrak{h}\otimes\R^3-D$ is the parameter associated to $\big(X, g_X, I_1^X, I_2^X, I_3^X\big)$, 
  and $A$ is chosen as in Lemma \ref{lem_AZ}.    
  There might a slight ambiguity here, since different $A\in \SO(3)$ could do -- 
  namely, given $\zeta$ as in Lemma \ref{lem_AZ}, there may be many $A$ satisfying its conclusions; 
  we will see however in Remarks \ref{rmk_AZ} and \ref{rmk_beth} that $\beth$ as we construct it is not affected by this choice. 

~

Points 1. and 5. of our program above do not need further developments. 
We postpone the tri-holomorphic isometry of Point 2. to Part 2, paragraph \ref{subsec_deg2hom}, 
since it will be easier to tackle with a little extra notions on Kronheimer's classification of ALE gravitational instantons. 
As for Point 3., as mentioned already, the proof of \ref{lem_beth} is given in section \ref{section_gluing} assuming results from part 2; 
apart from the proof of Lemma \ref{lem_AZ} which we shall settle now, 
our main task in the current part is thus the gluing and the subsequent corrections stated in Point 4., 
to which we devote sections \ref{section_gluing} and \ref{section_corrections} below, 
after recalling a few useful facts on the Taub-NUT metric seen as a Kähler metric on $(\C^2, I_1)$ 
in next section.  

~

\noindent
\textit{Proof of Lemma \ref{lem_AZ}}. 
 Fix $\zeta\in\mathfrak{h}\otimes\R^3$, 
 and define the matrix $Z(\zeta)=\big(\langle \zeta_j,\zeta_{\ell}\rangle\big)_{1\leq j,\ell\leq3}$ of its scalar products. 
 It is therefore elementary matrix calculus to check that the $\SO(3)$-action defined by \eqref{eqn_Azeta} 
 and referred to in the statement of the lemma translates into: $Z(A\xi)= AZ(\xi)A^t$. 
 
 We thus want to find $A\in \SO(3)$ such that $AZ(\zeta)A^t$ has shape: 
   \begin{equation}  \label{eqn_matrixshape}
    \begin{pmatrix}
     \mu & *       & *       \\
     *   & \lambda & 0       \\
     *   &    0    & \lambda
    \end{pmatrix}.
   \end{equation}

  Since $Z=Z(\zeta)$ is symmetric, 
  there exists $O\in \SO(3)$ such that $OZO^t=\diag(\lambda_1,\lambda_2,\lambda_3)$, 
  and we now look for $Q\in\SO(3)$ such that $Q\diag(\lambda_1,\lambda_2,\lambda_3)Q^t$ has shape \eqref{eqn_matrixshape}; 
  setting then $A=QO$ leads us to the conclusion. 
  If two of the $\lambda_j$ are the same then we are done, 
  up to letting act one of the permutation matrices 
  $\Big(\begin{smallmatrix} 1 & 0 &0 \\ 0 & 0 & 1 \\ 0 & -1 & 0 \end{smallmatrix}\Big)$, 
  $\Big(\begin{smallmatrix} 0 & 1 &0 \\ -1 & 0 & 0 \\ 0 & 0 & 1 \end{smallmatrix}\Big)$ and 
  $\Big(\begin{smallmatrix} 0 & 0 &1 \\  0 & 1 & 0 \\ -1 & 0 & 0 \end{smallmatrix}\Big)$. 
  Up to this action again, we can therefore assume $\lambda_1>\lambda_2>\lambda_3$. 
  
  Setting 
   \begin{equation*}
    Q=\begin{pmatrix}
       \big(\frac{\lambda_1-\lambda_2}{\lambda_1-\lambda_3}\big)^{1/2} & 0 
                         & \big(\frac{\lambda_2-\lambda_3}{\lambda_1-\lambda_3}\big)^{1/2}  \\
        0                & 1                                            & 0                  \\
       -\big(\frac{\lambda_2-\lambda_3}{\lambda_1-\lambda_3}\big)^{1/2} & 0 
                         & \big(\frac{\lambda_1-\lambda_2}{\lambda_1-\lambda_3}\big)^{1/2} 
      \end{pmatrix},
   \end{equation*}
 a direct computation gives 
 $Q\diag(\lambda_1,\lambda_2,\lambda_3)Q^t
     =\Big(\begin{smallmatrix} \lambda_1+\lambda_3-\lambda_2 &    0      & -\Lambda   \\ 
                                                           0 & \lambda_2 & 0          \\ 
                                                    -\Lambda &    0      & \lambda_2
           \end{smallmatrix}\Big) $, 
 where $\Lambda=(\lambda_1-\lambda_2)^{1/2}(\lambda_2-\lambda_3)^{1/2}$. 
 \cqfd

 \begin{rmk}  \label{rmk_AZ}
  Let us give a brief idea about how such a matrix $Q$ can be found. 
  This is actually what one can get by writing down the three relevant coefficients of $Q\diag(\lambda_1,\lambda_2,\lambda_3)Q^t$ 
  for $Q=(q_{j\ell})_{1\leq j,\ell\leq3}\in \SO(3)$, which leads to the underdetermined system 
   \begin{equation*}
    \left\{\begin{aligned}
            (\lambda_1-\lambda_2)q_{22}q_{32}+ (\lambda_1-\lambda_3)q_{23}q_{33}&      = 0, \\
            (\lambda_1-\lambda_2)q_{22}^2   + (\lambda_1-\lambda_3)q_{23}^2 \,\quad   & =  (\lambda_1-\lambda_2)q_{32}^2 +    (\lambda_1-\lambda_3)q_{33}^2,
            \end{aligned}
    \right.
  \end{equation*}
  to which one adds the arbitrary extra two conditions $q_{32}=q_{23}=0$. 
  Still keeping the same notations, 
  one can show that the only possibilities for writing $Q\diag(\lambda_1,\lambda_2,\lambda_3)Q^t$ under shape \eqref{eqn_matrixshape} 
  are the $\Big(\begin{smallmatrix} \lambda_1+\lambda_3-\lambda_2 &   \Lambda \cos\phi      & \Lambda \sin\phi  \\ 
                                                 \Lambda \cos\phi & \lambda_2               & 0                    \\ 
                                                 \Lambda \sin\phi &    0                    & \lambda_2
           \end{smallmatrix}\Big)$, 
  $\phi\in\R$, and again $\lambda_1\geq \lambda_2\geq \lambda_3$. 
 \end{rmk}

 \subsection{The Taub-NUT metric as a Kähler metric on $(\C^2,I_1)$} \label{section_TN}
  Before we proceed to the gluing of the Taub-NUT metric with the ALE metric of one of Kronheimer's instantons, 
  we recall a few facts about this very Taub-NUT metric on $\C^2$, 
  that will be used in the analytic upcoming sections \ref{section_gluing} and \ref{section_corrections}. 
  
  \subsubsection{Gibbons-Hawking versus LeBrun ansätze}
   \noindent
   \textit{Gibbons-Hawking ansatz.} 
   As recalled in the Introduction, 
   the Taub-NUT metric on $\R^4$ is often described via the Gibbons-Hawking ansatz as follows: 
   given $m\in (0,\infty)$, 
    \begin{equation}  \label{eqn_defTN}
     \f_m = V(dy_1^2+dy_2^2+dy_3^2)+ V^{-1}\eta^2,
    \end{equation}
   where $(y_1,y_2,y_3)$ is a circle fibration of $\R^4\backslash\{0\}$ over $\R^3\backslash\{0\}$, 
   $V$ is the function $\frac{1+4mR}{2R}$ (harmonic in the $y_j$ coordinates) with $R^2=y_1^2+y_2^2+y_3^2$, 
   and where $\eta$ is a connection 1-form for this fibration such that $d\eta=*_{\R^3}dV$. 
   Thus defined, the metric $\f_m$ confers length $\pi\sqrt{2/m}$ to the fibres at infinity, and is hermitian for the almost-complex structures 
    \begin{equation*} 
     J_a:\left\{ \begin{aligned}
                  Vdy_{a} &\longmapsto \eta, \\
                  dy_{b} &\longmapsto  dy_{c},
                 \end{aligned} \right.
    \end{equation*}
  with $(a,b,c)\in\big\{(1,2,3),(2,3,1),(3,1,2)\big\}$.  
  These are in fact complex structures, verifying the quaternionic relations $J_aJ_bJ_c=-1$, for which $\f_m$ is Kähler, thanks to the harmonicity of $V$: 
  $\f_m$ is thus hyperkähler. 
  One checks moreover that this way, the metric $\f_m$ and the complex structures extend as such through $0\in \R^4$. 
  The reader is referred to \cite{gh, leb} for details, 
  as we shall now switch point of view to a description better adapted to our construction.

 ~
  
 \noindent
 \textit{LeBrun's potential.} 
  As depicted in \cite{auv1}, after \cite{leb}, 
  one can give a more concrete support of this description,  
  through which, among others, the complex structure $J_1$ mentioned above is the standard $I_1$ on $\C^2$, 
  and $\vol^{\f_m}= \Omega_{\e}$, the standard euclidean volume form. 
  One starts with the following implicit formulas:
   \begin{equation}   \label{eqn_lebformulas}
    \begin{aligned}
     |z_1|&= e^{m(u^2-v^2)}u, \\
     |z_2|&= e^{m(v^2-u^2)}v,
    \end{aligned}
   \end{equation}
  defining functions $u,v:\C^2\to \R$, invariant under the circle action $e^{i\theta}\cdot(z_1,z_2)= (e^{i\theta}z_1,e^{-i\theta}z_2)$, 
  making $\mathbb{S}^1$ as a subgroup of $\SU(2)$; 
  notice the role of $m$ in these formulas, which enlightens our choice of taking it as the parameter of the upcoming construction. 
  One then sets $y_1=\frac{1}{2}(u^2-v^2)$, $y_2+iy_3= -iz_1z_2$, $R=\frac{1}{2}(u^2+v^2)=\big(y_1^2+y_2^2+y_3^2\big)^{1/2}$; 
  these are $\mathbb{S}^1$-invariant functions, 
  making $(y_1,y_2,y_3)$ as a principal-$\mathbb{S}^1$ fibration $\C^2\to\R^3$ away from the origins. 
  One finally defines:
    \begin{equation}  \label{eqn_TNpot}
    \varphi_m:= \frac{1}{4}\big(u^2+v^2+m(u^4+v^4)\big)= \frac{1}{2}\big(R + m(R^2+y_1^2)\big).
   \end{equation}
  One can then checks that $dd^c_{I_1}\varphi_m$ is positive is the sense of $I_1$-hermitian 2-forms, 
  and that $(dd^c_{I_1}\varphi_m)^2=2\Omega_{\e}$. 
  If one sets moreover $V= \frac{1+4mR}{2R}$, and $\eta=I_1Vdy_1$,  
  noticing by passing that $\eta$ is then a connection 1-form for the fibration, $d\eta=*_{\R^3}dV$, 
  one has: $\f_m := V(dy_1^2+dy_2^2+dy_3^2)+ V^{-1}\eta^2= (dd^c_{I_1}\varphi_m)(\cdot,I_1\cdot)$. 
  This metric is well-defined at $0\in \C^2$, as $(dd^c_{I_1}\varphi_m)(\cdot,I_1\cdot)=\e$ at that point. 
  
  The metric $\f_m$ is therefore Kähler for $I_1$ with volume form $\vol^{\f_m}=\Omega_{\e}$ on the whole $\C^2$; 
  by the standard properties of Kähler metrics, it is thus Ricci-flat.  
  One recovers a complete hyperkähler data after checking that the defining equations
   \begin{equation}  \label{eqn_defJj}
    \f_m(J_j\cdot,\cdot) = \omega_j^{\e}, \qquad \text{where} \qquad\omega^{\e}_j=\e(I_j\cdot,\cdot), 
     \qquad j=2,3,
   \end{equation}
  with $I_2$, $I_3$ the other two standard complex structures on $\R^4\cong \mathbb{H}$, 
  give rise to complex structures, verifying respectively $Vdy_j\mapsto\eta$, 
  $dy_k\mapsto dy_i$ for $(i,j,k)\in \big\{(2,3,1), (3,1,2)\big\}$,  
  as well as the quaternionic relations together with $I_1$.

 ~
  
  Let us now give a look at the length of the $\mathbb{S}^1$-fibres at infinity. 
  Consider the vector field 
  $\xi:=i\big(z_1\frac{\partial}{\partial z_1}-\overline{z_1}\frac{\partial}{\partial \overline{z_1}}
              -z_2\frac{\partial}{\partial z_2}+\overline{z_2}\frac{\partial}{\partial \overline{z_2}}\big)$ 
  giving the infinitesimal action of $\mathbb{S}^1$. 
  One has $dy_1(-I_1\xi)= V^{-1}$,  
  thus $\eta(\xi)=1$, and $dy_j(\xi)=0$, $j=1,2,3$; 
  since $R$ is $\mathbb{S}^1$-invariant, the length of the fibres is just $2\pi V^{-1/2}$, which tends to $\pi\sqrt{2/m}$. 
   \begin{rmk}
    One can check that even if we can let $m$ vary, this description actually leads to essentially one metric; 
    indeed, if $\kappa_s$ is the dilation of factor $s>0$ of $\R^4$, 
    one gets with help of formulas \eqref{eqn_lebformulas} and \eqref{eqn_TNpot}
    the following homogeneity property: $\kappa_s^*\f_m= s^2\f_{ms^2}$, 
    which is of course coherent with the length of the fibres at infinity and the fact that $\vol^{\f_m}=\vol^{\f_{ms^2}}=\Omega_{\e}$. 
   \end{rmk}

~

  \textit{\textbf{From now on, we see the mass parameter $m$ as fixed, and we drop the indexes $m$ when there is no risk of confusion. }}  

~
 
  \noindent
  \textit{The Taub-NUT metric and the action of the binary dihedral group on $\C^2$.} 
  For $k\geq 2$, which we fix until the end of this part, 
  the action of the \textit{binary dihedral group} $\mathcal{D}_k$ \textit{of order} $4k$ 
  seen as a subgroup of $\SU(2)$ is generated by the matrices $\zeta_k:= \Big(\begin{smallmatrix} e^{i\pi/k} & 0 \\ 0 & e^{-i\pi/k} \end{smallmatrix}\Big)$ 
  and $\tau:= \big(\begin{smallmatrix} 0 & -1 \\ 1 & 0 \end{smallmatrix}\big)$. 
  One has: $\zeta_k^*y_j=y_j$, $j=1,2,3$, and thus $\zeta_k^*R=R$, and $\zeta_k^*\eta=\eta$, 
  whereas: $\tau^*y_j=-y_j$, $j=1,2,3$, thus $\tau^*R=R$, and $\tau^*\eta=\eta$.  
  The Taub-NUT metric $\f$ is therefore $\mathcal{D}_k$-invariant, and descends smoothly to $(\R^4\backslash\{0\})/\mathcal{D}_k$: 
  this is the metric we are going to glue at infinity of $\mathcal{D}_k$-ALE instantons in the next section. 
  Before though, we need a few more analytical tools for the Taub-NUT metric as we describe it.   
  

 \subsubsection{Orthonormal frames, covariant derivatives and curvature}  \label{subsection_TN}
  In addition to the above relations between the vector field $\xi$, and the 1-forms $\eta$ and $dy_j$, $j=1,2,3$, 
  one has that the data 
   \begin{equation}   \label{eqn_ei}
     (e_0,e_1,e_2,e_3):=\big(V^{1/2}\xi, -I_1V^{1/2}\xi, V^{-1/2}\zeta, V^{-1/2}I_1\zeta\big),  
   \end{equation}
  is the dual frame of the orthonormal frame of 1-forms 
   \begin{equation}  \label{eqn_ei*}
    (e_0^*,e_1^*,e_2^*,e_3^*):=\big(V^{-1/2}\eta, V^{1/2}dy_1, V^{1/2}dy_2, V^{1/2}dy_3\big)
   \end{equation}
   on $\C^2\backslash\{0\}$, 
  provided that the vector field $\zeta$ is defined by:
   \begin{equation}  \label{eqn_defzeta}
    \zeta:= \frac{1}{2iR}\Big(e^{4my_1}\big(z_2\frac{\partial}{\partial\overline{z_1}}-\partial\overline{z_2}\frac{\partial}{\partial z_1}\big)
                              + e^{-4my_1}\big(z_1\frac{\partial}{\partial\overline{z_2}}-\partial\overline{z_1}\frac{\partial}{\partial z_2}\Big), 
   \end{equation}
 see \cite[Lemma 1.11]{auv1}; 
 we keep the notations $(e_j)_{j=0,\dots,3}$ and $(e_j^*)_{j=0,\dots,3}$ throughout this part. 
 An explicit computation made in \cite[\S 1.2.3]{auv1} then gives us the estimate 
  \begin{equation*}
   \big|(\nabla^{\f})^{\ell}e_j\big|_{\f}= O\big(R^{-1-\ell}\big)\quad\text{near infinity for all }\ell\geq1\text{ and }j=0,\dots,3.
  \end{equation*}
 Notice that consequently, 
 for all $\ell\geq0$, $\big|(\nabla^{\f})^{\ell}\riem^{\f}\big|_{\f}=O\big(R^{-3-\ell})$; 
 this estimate, done using the Gibbons-Hawking ansatz, can also be found in \cite[\S 1.0.3]{min1}. 
  
~

 We close this section by two further useful estimates, 
 which may give an idea of the geometric gap between $\e$ and $\f$: 
 first, at the level of distance functions, 
 rearranging formulas \eqref{eqn_lebformulas} gives: $R\leq 2r^2$, which is sharp is general; 
 second, there exists $C=C(m)>0$ such that outside the unit ball of $\C^2$, 
 $C^{-1}r^{-2}\e\leq \f\leq Cr^2\e$, and again this is sharp in general.

\subsection{Gluing the Taub-NUT metric to an ALE metric} \label{section_gluing}
   
  As is usual when gluing Kähler metrics, 
  we shall work on potentials to glue the ALF model-metric to an ALE one.  
  The previous section gives us the potential $\varphi$ for the ALF metric (equation \eqref{eqn_TNpot}); 
  the following paragraph provides us with a sharp enough potential for the ALE metric. 

  \subsubsection{Approximation of the ALE Kähler form as a complex hessian}
   
   \noindent
   \textit{Asymptotics of the Kähler form and the complex structure.} 
   In view of Step 3. and 4. of the program developed in section \ref{section_prgrm}, 
   since we are performing our gluing on some specific ALE instantons, we fix \textit{for the rest of this part} 
    \begin{equation}   \label{eqn_xi}
      \xi\in \mathfrak{h}-D, \qquad\text{such that:}\quad |\xi_2|^2-|\xi_3|^2=\langle\xi_2,\xi_3\rangle=0 
    \end{equation}
   and consider the associated ALE instanton $\big(Y, g_Y, I_1^Y, I_2^Y, I_3^Y\big)$.  
   Lemma \ref{lem_beth} gives an ALE diffeomorphism $\Phi_Y: Y\backslash K\rightarrow (\R^4\backslash B)/\mathcal{D}_k$,   
   where $K$ is some compact subset of $Y$ and $B$ a ball in $\R^4$ centred at the origin;   
   recall that by "ALE diffeomorphism" we mean that for all $\ell\geq 0$, 
    \begin{equation*}
     \big|(\nabla^{\e})^{\ell}({\Phi_Y}_*g_Y-\e)\big|_{\e} = O(r^{-4-\ell}),  
    \end{equation*}
   and likewise on the complex structures. 
   Before using the more specific properties of $\Phi_Y$ at the level of complex structures, 
   let us mention the following:
   %
   we want to proceed to a gluing of Kähler metrics,   
   and the convenient way of doing so is to glue \textit{the Kähler forms, via their potentials}. 
   We already have a candidate for the potential of an ALF metric at infinity at hand: 
   as evoked, this would be ${\Phi_Y}_*\varphi^{\flat}$ (see Point 4. in section \ref{section_prgrm}). 
   Conversely, we need to kill the ALE metric near infinity, and for this we want a sharp enough potential, 
   in a sense that we make clear below, see Proposition \ref{prop_gluing}. 
   We thus need for this a sharp knowledge of \textit{the Kähler form} $\omega_1^Y:=g_Y(I_1^Y\cdot, \cdot)$, 
   and since we are about to compute $I_1^Y$-complex hessians as well, 
   we also need a precise description of this complex structure. 
   These are given by the following, from which Lemma \ref{lem_beth} actually follows as we shall see at the end of this section, with the same $\Phi_Y$:
    \begin{lem}  \label{lem_omega1Y}
     One can choose the ALE diffeomorphism $\Phi_Y$ such that 
      \begin{equation}
       {\Phi_Y}_* \omega_1^Y = \omega_1^{\e} 
                  - c\big(|\xi_1|^2\theta_1 + \langle\xi_1,\xi_2\rangle\theta_2 + \langle\xi_1,\xi_3\rangle\theta_3\big) + O(r^{-8})
      \end{equation}
     where $c>0$ is some universal constant, $\theta_j= \frac{1}{4}dd^c_{I_j}\big(r^{-4}\big)$, $j=1,2,3$, on the one hand, 
     and if $\iota_1^Y$ denotes ${\Phi_Y}_*I_1^Y-I_1$, then it is given by:
      \begin{equation}
       \e(\iota_1^Y\cdot,\cdot)= - c\big(|\xi_2|^2+|\xi_3|^2\big)\frac{rdr\cdot \alpha_1}{r^6} + O(r^{-8})
      \end{equation}
     where $c$ is the same constant as above and $\alpha_1=I_1rdr$, on the other hand. 
     
     We can moreover assume that ${\Phi_Y}_*\Omega_Y=\Omega_{\e}$, where $\Omega_{Y}=\vol^{g_Y}$. 
    \end{lem}
 In this statement the error terms $O(r^{-8})$ are understood in the "euclidean way",  
 namely for any $\ell\geq0$, the $\ell$th $\nabla^{\e}$-derivatives of these tensors are $O(r^{-8-\ell})$. 
 This lemma requires further notions on Kronheimer's, 
 and is more precisely a direct application of Theorem \ref{thm_ALE} of Part 2 to $Y$ with $\xi$ verifying \eqref{eqn_xi}.  
 Notice however the order of the error term, which is $-8$ whereas $-6$ would be expected, 
 if one think for instance to the Eguchi-Hanson metric (\cite[Ex. 7.2.2]{joy}); 
 this is crucial in proving Lemma \ref{lem_beth}, and this is specific to (groups containing) dihedral binary groups. 
 Besides, the assertion on the volume forms is only needed in next paragraph. 
 
 ~

 \noindent
 \textit{Approximating $\omega_1^Y$ as an $I_1^Y$-complex hessian.} 
 We shall see for now how Lemma \ref{lem_omega1Y} allows us to approximate the Kähler form $\omega_1^Y$ as an $I_1^Y$-complex hessian, 
 with respect to the Taub-NUT metric pushed-forward to $Y$:
  \begin{prop}  \label{prop_approx_omega1Y}
   Take $\Phi_Y$ as in Lemmas \ref{lem_beth} and \ref{lem_omega1Y}, 
   and denote by $\tilde{\f}$ a smooth extension of ${\Phi_Y}^*\f$ on $Y$. 
   Then there exists a function $\Psi$ on $Y$ such that near infinity, 
    \begin{equation}  \label{eqn_approx_omega1Y}
     \big|(\nabla^{\tilde{\f}})^{\ell}\big(\omega_1^Y-dd^c_{I_1^Y}\Psi\big)\big|_{\tilde{\f}}=O(R^{-2}), \qquad \ell=0,1,2. 
    \end{equation}
   More precisely, $\Psi$ can be decomposed as a sum ${\Phi_Y}^*\Psi_{\euc}+{\Phi_Y}^*\Psi_{\mxd}$, 
   where on the one hand, $\Psi_{\euc}=O(r^2)$, $\big|\Psi_{\euc}\big|_{\e}=O(r)$, and 
    \begin{equation}  \label{eqn_approx_omega1Y1}
     \big|(\nabla^{\e})^{\ell}\big(\omega_1^{\e}-c|\xi_1|^2\theta_1-dd^c_{{\Phi_Y}_*I_1^Y} \Psi_{\euc}\big)\big|_{\e}= O(r^{-8-\ell})
       \qquad\text{for all } \ell\geq0, 
    \end{equation}
   and on the other hand, $\Psi_{\mxd}, \big|d\Psi_{\mxd}\big|_{\f}=O(R^{-1})$, and 
    \begin{equation}  \label{eqn_approx_omega1Y2}
     \big|(\nabla^{\f})^{\ell}
          \big(-c(\langle\xi_1,\xi_2\rangle\theta_2 + \langle\xi_1,\xi_3\rangle\theta_3)-dd^c_{{\Phi_Y}_*I_1^Y} \Psi_{\mxd}\big) \big)\big|_{\f}=O(R^{-2}), 
     \quad \ell=0,1,2.
    \end{equation}

  \end{prop}
 \prf. Notice that once the statement on $\Psi_{\euc}$ (the "euclidean component" of $\Psi$) and $\Psi_{\mxd}$ (the "mixed component") 
  are known, estimates \eqref{eqn_approx_omega1Y} follow at once by transposition of estimates \eqref{eqn_approx_omega1Y1} and \eqref{eqn_approx_omega1Y2}
  and of the expansion of $\omega_1^Y$ stated in Lemma \ref{lem_omega1Y}  to $Y$, 
  keeping the following fact in mind: 
  
   \begin{fact}
    If $\alpha$ is a tensor of type $(2,0)$, $(1,1)$ or $(0,2)$ such that $\big|(\nabla^{\e})^{\ell}\alpha\big|_{\e}=O(r^{-2a-\ell})$, $a\geq1$, $\ell=0,1,2$, on $\R^4$, 
  then $\big|(\nabla^{\f})^{\ell}\alpha\big|_{\f}=O(R^{1-a})$, $\ell=0,1,2$.
   \end{fact}

  This fact takes into account estimates such as $R=O(r^2)$ and $C^{-1}r^{-2}\e\leq \f\leq Cr^2\e$ at level $\ell=0$, 
  and roughly says that by passing from euclidean to ALF geometry we do not win regularity by differentiation anymore;  
  it is proved in \cite[p.22]{auv1}.  
  
  We hence come to the statements on $\Psi_{\euc}$ and $\Psi_{\mxd}$. 
  We consider before starting a large constant $K$ such that the image of $\Phi_Y$ is contained in both $\{r\geq K\}\subset\R^4/\mathcal{D}_k$ 
  and $\{R\geq K\}\subset\R^4/\mathcal{D}_k$, and define a cut-off function $\chi:\R\to [0,1]$ such that:
   \begin{equation} 
    \left\{\begin{aligned}
            \chi(t) &=0 \quad \text{on} \quad t\leq K-1, \\
            \chi(t) &=1 \quad \text{on} \quad t\geq K,  
           \end{aligned}
    \right.
   \end{equation}
  which will be useful when defining functions to be pulled-back to $Y$ via $\Phi_Y$. 

~

  \noindent 
  \textit{The euclidean component $\Psi_{\euc}$.} 
   In an asymptotically euclidean setting, a natural first candidate for the potential of a Kähler form is $\frac{1}{4}r^2$. 
   Now remember we are working with $I_1^Y$ -- or more exactly with ${\Phi_Y}_*I_1^Y$, 
   but we forget about the push-forward here for simplicity of notation; 
   following Lemma \ref{lem_omega1Y}, a straightforward computation gives, near infinity in $\R^4$:
    \begin{align*}
     dd^c_{I_1^Y}\Big(\frac{1}{4}r^2\Big)&= \frac{1}{2}d\big[(I_1+\iota_1^Y)rdr\big] 
                                          = \frac{1}{2}d\big[\alpha_1 + c(|\xi_2|^2+|\xi_3|^2)r^{-4}\alpha_1+O(r^{-7})\big] \\
                                         &= \omega_1^{\e} - c(|\xi_2|^2+|\xi_3|^2)\theta_1 +O(r^{-8}), 
    \end{align*}
   where the $O$ are understood in the euclidean way. 
   On the other hand observe that $I_1d(r^{-2})=-2 r^{-4}\alpha_1$, and thus
    \begin{equation*}
     dd^c_{I_1^Y}(r^{-2})= d\big[(I_1+\iota_1^Y)d(r^{-2})\big]= d\big[-2r^{-4}\alpha_1+O(r^{-7})\big]=4\theta_1 +O(r^{-8}). 
    \end{equation*}
   Now define 
    \begin{equation*}
     \Psi_{\euc} = \frac{1}{4}\chi(r)\big(r^2 + c(|\xi_2|^2+|\xi_3|^2-|\xi_1|^2)r^{-2}\big); 
    \end{equation*}
   on $\R^4/\mathcal{D}_k$ (it is $\mathcal{D}_k$-invariant); 
   it has support in the image of $\Phi_Y$, has the growth stated in the lemma as well as its differential,  
   and by the previous two estimates we get that $\omega_1^{\e}-c|\xi_1|^2\theta_1-dd^c_{{\Phi_Y}_*I_1^Y} \Psi_{\euc}=O(r^{-8})$ 
   for $\e$ with according decay on the derivatives, as wanted. 

 ~
  
  \noindent 
  \textit{The mixed component $\Psi_{\mxd}$.} 
   The main reason why we could construct $\Psi_{\euc}$ such as to reach estimates \eqref{eqn_approx_omega1Y1} 
   is essentially that $\theta_1$ can be realised as an $I_1$-complex hessian, at least away from 0. 
   Now realising $\theta_2$ and $\theta_3$ as $I_1$-complex hessians as well does not seem possible: see \cite[p.202]{joy} on that matter. 
   Nonetheless, $\theta_2$ and $\theta_3$ may not be so problematic when looked at via $\f$. 
   We can indeed approximate them precisely enough with respect to this metric by the $I_1$ or $I_1^Y$-complex hessians 
   of some well-chosen $\mathcal{D}_k$-invariant functions, 
   provided that we partially leave the euclidean world and use also functions coming from Taub-NUT geometry, e.g. $y_1$ and $R$ 
   (hence the previous dichotomy "euclidean/mixed"):
    \begin{lem}  \label{lem_est_psidiff}
     Consider the complex valued function
      \begin{equation*}
        \psi_{c} := -2\frac{(y_2+iy_3)\sinh (4my_1)}{r^2R}
      \end{equation*}
     on $\R^4\backslash\{0\}$. 
     Then near infinity:
      \begin{enumerate}
       \item $\big|(\nabla^{\f})^{\ell}\psi_{c}\big|_{\f}= O(R^{-1})$ for $\ell=0,\dots,4$;  
       \item $\big|(\nabla^{\f})^{\ell}\big(dd^c_{I_1}\psi_c-(\theta_2+i\theta_3)\big)\big|_{\f} = O(R^{-2})$ 
          for $\ell=0,1,2$, and these estimates hold for $I_1^Y$ as well. 
      \end{enumerate}
    \end{lem}
    The proof of this crucial lemma is essentially computational, 
    which is why we postpone it to section \ref{section_lemmas}. 
    For now set $\psi_2= \mathfrak{Re}(\psi_c)$ and $\psi_3= \mathfrak{Im}(\psi_c)$, and define 
     \begin{equation*}
      \Psi_{\mxd} := -c\chi(R)(\langle\xi_1,\xi_2\rangle \psi_2 + \langle\xi_1,\xi_3\rangle\psi_3). 
     \end{equation*}
    In view of Lemma \ref{lem_est_psidiff}, 
    such a function, defined on the image of $\Phi_Y$, verifies the growth assertions of Proposition \ref{prop_approx_omega1Y}, 
    as well as the estimates \eqref{eqn_approx_omega1Y2}: Proposition \ref{prop_approx_omega1Y} is proved. 
  \cqfd

~

  We are now in position to perform the gluing advertised in Point 4. of the program of section \ref{section_prgrm}. 
  This is done in next paragraph to which the reader may jump directly, 
  since we conclude the current paragraph by the proof of Lemma \ref{lem_beth}, assuming Lemma \ref{lem_omega1Y} 
  (and more precisely the assertion on $I_1^Y$ in that statement). 

~ 
  
  \prf \textit{of Lemma \ref{lem_beth} following Lemma \ref{lem_omega1Y}}. 
   We fix $\Phi_Y$ as in Lemma \ref{lem_omega1Y}; 
   we work on $\R^4$, and to simplify notations we forget about the push-forwards by $\Phi_Y$. 

   We are thus looking for a diffeomorphism $\beth$ of $\R^4$ such that 
   $\big|I_1^Y - \beth^*I_1\big|_{\e}=O(r^{-8})$, 
   with according decay on euclidean derivatives -- 
   until the end of this proof we forget about ALF geometry and stick to the euclidean setting;  
   we will thus content ourselves with using $O$ in this euclidean meaning. 
   An explicit formula is given for $\beth$ in the statement of Lemma \ref{lem_beth}, 
   which is: $\beth(z_1,z_2)\mapsto \big(1+\frac{a}{\kappa+r^4}\big)(z_1,z_2)$ with $(z_1,z_2)$ the standard complex coordinates on $(\C^2, I_1)$; 
   since the value of $\kappa$ does not affect asymptotic considerations -- changing $\kappa$ only contributes as a $O(r^{-8})$ --, 
   we could thus, up to determining the value of the constant $a$, 
   simply check that such a $\beth$ meets our requirement, in light of the asymptotics for $I_1^Y$ stated in Lemma \ref{lem_omega1Y}. 
   
   We prefer nonetheless the following more constructive approach. 
   If we are to look at some $\beth$ as in the statement, 
  we should certainly take it with shape $(z_1,z_2)\mapsto (z_1+\vareps_1,z_2+\vareps_2)$, 
  with $\vareps_j=O(r^{-4})$, $j=1,2$. 
  The condition $I_1^Y-\beth^*I_1=O(r^{-8})$ can be rewritten as a condition on $\vareps_1$: 
  $I_1^Y\beth^*dz_1=\beth^*(I_1dz_1)+O(r^{-8})= \beth^*(idz_1)+O(r^{-8})$, 
  i.e. $I_1^Y(dz_1+d\vareps_1)= i(dz_1+d\vareps_1)+O(r^{-8})$. 
  Recall the writing $I_1^Y=I_1+\iota_1^{Y}$;  
  the previous condition hence gives us: 
  $I_1dz_1+I_1d\vareps_1+ \iota_1^{Y}dz_1=i(dz_1+d\vareps_1)+O(r^{-8}) $. 
  Since $I_1dz_1=idz_1$ and $I_1d\vareps_1=I_1(\partial\vareps_1+\dbar\vareps_1)=i(\partial\vareps_1-\dbar\vareps_1)$ 
  (with $\partial$ and $\dbar$ those attached to $I_1$), 
  the final condition is: $2i\dbar\vareps_1=\iota_1^{Y}dz_1+O(r^{-8})=-dz_1(\iota_1^{Y}\cdot)+O(r^{-8})$. 

 Set $\alpha_j=I_jrdr$, $j=2,3$; 
 from Lemma \ref{lem_omega1Y}, $rdr(\iota_1^{Y}\cdot)=\tfrac{c(|\xi_2|^2+|\xi_3|^2)}{r^4}\alpha_1+O(r^{-7})$, 
 $\alpha_1(\iota_1^{Y}\cdot)=\tfrac{c(|\xi_2|^2+|\xi_3|^2)}{r^4}rdr+O(r^{-7})$, 
 $\alpha_2(\iota_1^{Y}\cdot),\alpha_3(\iota_1^{Y}\cdot)=O(r^{-7})$. 
 Hence from the equality 
  \begin{equation*}
   dz_1=\frac{1}{r^2}\big[z_1(rdr+i\alpha_1)-\overline{z_2}(\alpha_2+i\alpha_3)\big],
  \end{equation*} 
 we have: $dz_1(\iota_1^{Y}\cdot)=\tfrac{c(|\xi_2|^2+|\xi_3|^2)}{r^6}z_1(\alpha_1+irdr)+O(r^{-8})=
             \tfrac{ic(|\xi_2|^2+|\xi_3|^2)}{r^6}\big(z_1^2d\overline{z_1}+z_1z_2d\overline{z_2}\big)+O(r^{-8})$. 
 At last we must thus solve 
  \begin{equation*}
   2\dbar\vareps_1 = -  \frac{c(|\xi_2|^2+|\xi_3|^2)}{r^6}\big(z_1^2d\overline{z_1}+z_1z_2d\overline{z_2}\big)+O(r^{-8}). 
  \end{equation*}
 One easily checks that $\vareps_1= \frac{c(|\xi_2|^2+|\xi_3|^2)z_1}{4r^4}$ is an exact solution. 
 A similar analysis leads us to $\vareps_2= \frac{c(|\xi_2|^2+|\xi_3|^2)z_2}{4r^4}$, 
 and one checks easily that this way, one has indeed $I_1^Y-\beth^*I_1=O(r^{-8})$. 
 The last point to be dealt with is the singularity of the $\vareps_j$ at 0; 
 one can nonetheless take instead $\vareps_j= \frac{cz_1}{1+4r^4}$.  
 One can even take $\vareps_j= \frac{c(|\xi_2|^2+|\xi_3|^2)z_1}{4(\kappa+r^4)}$ 
 with $\kappa\geq1$ large enough so that $\beth=\id_{\C^2}+(\vareps_1,\vareps_2)$ is a diffeomorphism of $\C^2$; 
 we leave it to the reader as an exercise to check that 
 $\kappa=20c(|\xi_2|^2+|\xi_3|^2)$ seems is sufficient).   

 The estimate $\beth^*\Omega_{\e}-\Omega_{\e}=O(r^{-8})$ amounts to seeing that 
 $\mathfrak{Re}\big(\tfrac{\partial\vareps_1}{\partial z_1}+\tfrac{\partial\vareps_2}{\partial z_2}\big)=O(r^{-8})$: 
 extend $id(z_1+\vareps_1)\wedge d(\overline{z_1}+\overline{\vareps_1})\wedge 
          id(z_2+\vareps_2)\wedge d(\overline{z_2}+\overline{\vareps_2})$, 
 and look at the linear terms in $\vareps_1$, $\vareps_2$. 
 Since after multiplication by $a:=\frac{c(|\xi_2|^2+|\xi_3|^2)}{4}$ the error would again be $O(r^{-8})$, 
 we can do this computation with $\tfrac{z_1}{r^4}$ and $\tfrac{z_2}{r^4}$ playing the respective roles of $\vareps_1$ 
 and $\vareps_2$. 
 Now $\tfrac{\partial}{\partial z_j}\big(\tfrac{z_j}{r^4}\big)= \tfrac{1}{r^4}-\tfrac{2|z_j|^2}{r^6}$, $j=1,2$. 
 Since these are real, we only need to compute the sum 
 $\tfrac{\partial}{\partial z_1}\big(\tfrac{z_1}{r^4}\big)+\tfrac{\partial}{\partial z_2}\big(\tfrac{z_2}{r^4}\big)$, 
 which is $\tfrac{2}{r^4}-\tfrac{2|z_1|^2}{r^6}-\tfrac{2|z_2|^2}{r^6}=0$. 
 
 The $\mathcal{D}_k$-invariance of $\beth$ thus constituted is clear. 
 \cqfd

 \begin{rmk} \label{rmk_beth}
  According to the preceding proof, 
  $\beth$ as we construct it depends only 
  on $c(|\xi_2|^2+|\xi_3|^2)$. 
  If now $\xi$ is chosen as an $A\zeta$, $A\in \SO(3)$, $\zeta\in \mathfrak{h}-D$, 
  so as to satisfy condition \eqref{eqn_xi} as is evoked in Point 3. in the program of section \ref{section_prgrm}, 
  by Remark \ref{rmk_AZ}, $|\xi_2|^2=|\xi_3|^2$ does not depend on $A$, 
  and has to be the middle eigenvalue of the matrix $(\langle \zeta_j,\zeta_{\ell}\rangle)$.
  Consequently, $\beth=\beth_{A\zeta}$ does not depend on $A\in \SO(3)$. 
 \end{rmk}

  \subsubsection{The gluing}  \label{subsec_gluing}

   We keep the notations of the previous paragraph: 
   $\big(Y, g_Y, (I^Y_j)_{j=1,2,3}\big)$ is a $\mathcal{D}_k$-ALE instanton with parameter $\xi$ verifying \eqref{eqn_xi}, 
   $\Phi_Y$ an asymptotic isometry between infinities of $Y$ and $\R^4/\mathcal{D}_k$ fixed by Lemma \ref{lem_omega1Y}, 
   and $\beth$ is given by Lemma \ref{lem_beth} which we may also see as as diffeomorphism of $(\R^4\backslash\{0\})/\mathcal{D}_k$. 
   
   As alluded above, the form we want to glue $\omega_1^{Y}=g_{Y}(I_1^{Y}\cdot,\cdot)$ with at infinity is $dI_1^{Y}d\varphi^{\flat}$, 
   where $\varphi=\varphi_m$ is LeBrun's $I_1$-potential for $\f$ given by \eqref{eqn_TNpot}, 
   and where $\varphi^{\flat}=\beth^*\varphi$. 
   We set likewise $\f^{\flat}= \beth^*\f$, both on $\R^4$ and its quotient. 
   Recall that $\Psi=\Psi_{\euc}+\Psi_{\mxd}$ is defined in Proposition \ref{prop_approx_omega1Y} as an approximate $I_1^Y$-complex hessian of $\omega_1^{Y}$.  
   The following proposition, which is the analogue of \cite[Prop. 2.3]{auv1}, explains how to glue $dI_1^{Y}d\varphi^{\flat}$ to $\omega_1^{Y}$, so as to obtain an ALF metric on $Y$ at the end; 
   notice that we lose one order of precision in the asymptotics though:
    \begin{prop}  \label{prop_gluing}
      Take $K\geq0$ so that the identification $\Phi_{Y}$ between infinities of $\R^4/\mathcal{D}_k$ and $Y$ is defined on $\varphi\geq K$.
      Consider $r_0 \gg 1$, $\beta\in(0,1]$ and set 
      \begin{equation*}
       \Phi_m^{\flat} =  \kappa\circ\big(\varphi^{\flat}+\Psi_{\mxd}-K\big)- \chi\big((r-r_0)^{\beta}\big)\widetilde{\Psi}_{\euc},
      \end{equation*}
      where $\kappa:\R\to \R$ is a convex function vanishing on $(-\infty,0]$ and equal to $\id_{\R}$ on $[1,\infty)$, 
      $\chi$ is the cut-off function $\frac{d\kappa}{dt}$, and 
       $\widetilde{\Psi}_{\euc}:=\chi(r-r_0)\Psi_{\euc}$. 
     Then if the parameters $K$ and $r_0$ (resp. $\beta$) are chosen large enough (resp. small enough), 
     the symmetric 2-tensor $g_m$ associated via $I_1^{Y}$ to the $I_1^{Y}$-(1,1)-form
   \begin{equation*}
    \omega_m:= \omega_1^{Y}+dd^c_{I_1^{Y}}\Phi_m^{\flat}
   \end{equation*}
  is well-defined on the whole $Y$, is a Kähler metric for $I_1^{Y}$, is ALF in the sense that 
   $ \big|(\nabla^{\f^{\flat}})^{\ell}(g_m-\f^{\flat})\big|_{\f^{\flat}}=O\big(R^{-2}\big)$ 
   for $\ell=0,1,2$,
     and its volume form $\Omega_m$ verifies
   \begin{equation}   \label{eqn_estvolform}
    \big|(\nabla^{\f^{\flat}})^{\ell}(\Omega_m-\Omega_{Y})\big|_{\f^{\flat}}=O\big(R^{-2}\big)
   \end{equation}
for $\ell=0,1,2$, where $\Omega_{}$ is the volume form of the ALE metric $g_{Y}$.  
 \end{prop}
 \prf. 
  The positivity of $g_m$ and its asymptotics follow from the arguments analogous to those developed in the proof of Proposition 2.3 in \cite{auv1}, 
  which is why we are brief. 
  We nonetheless underline a few necessary adjustment; 
  for example, we will keep the following comparison between $\f$ and its correction $\f^{\flat}= \beth^*\f$ in mind:
   \begin{lem}   \label{lem_f/fflat}
    For $\ell=0,1,2$, we have: 
    $\big|(\nabla^{\f})^{\ell}(\f^{\flat}-\f)\big|_{\f}=O(R^{-1})$ on $\R^4$. 
    Moreover $\beth^*R=R+O(R^{-1})$. 
   \end{lem}

  We first consider the closed $I_1^Y$-hermitian form $dd^c_{I_1^{Y}}\kappa\circ\big(\varphi^{\flat}+\Psi_{\mxd}-K\big)$ on $Y$. 
  Even though $K$ is not fixed yet, this is equal to $dd^c_{I_1^{Y}}\big(\varphi^{\flat}+\Psi_{\mxd}\big)$ on $\{\varphi^{\flat}+\Psi_{\mxd}\geq K+1\}$ 
  seen on $Y$ via $\Phi_Y$ -- 
  this is possible for $K$ large enough since $\varphi^{\flat}+\Psi_{\mxd}$ is proper on $\R^4$ as $\varphi^{\flat}\geq \beth^*{R}\sim R$ (by Lemma \ref{lem_f/fflat})
  and $\Psi_{\mxd}=O(R^{-1})$. 
  Moreover $\kappa$ is convex, 
  and thus $dd^c_{I_1^{Y}}\kappa\circ\big(\varphi^{\flat}+\Psi_{\mxd}-K\big)$ is non-negative wherever $dd^c_{I_1^{Y}}\big(\varphi^{\flat}+\Psi_{\mxd}\big)$ is, 
  which is indeed the case near infinity. 
  Since indeed $\big|dd^c_{I_1^{Y}}\Psi_{\mxd}\big|_{\f}=O(R^{-1})$, 
  this claim will be checked if we prove the estimate:
   \begin{equation} \label{eqn_estmte1}
    \Big|dd^c_{I_1^{Y}}\varphi^{\flat} 
           - \frac{1}{2}\big[\f^{\flat}\big(I_1^Y\cdot,\cdot\big)- \f^{\flat}\big(\cdot,I_1^Y\cdot\big)\big]\Big|_{\f^{\flat}}=O(R^{-2})
   \end{equation}
  as $\varpi_{\f^{\flat}}:=\frac{1}{2}\big[\f^{\flat}\big(I_1^Y\cdot,\cdot\big)
                                                       - \f^{\flat}\big(\cdot,I_1^Y\cdot\big)\big]$ 
  is nothing but the $I_1^Y$-hermitian form associated to the $I_1^Y$-hermitian metric 
  $\frac{1}{2}\big[\f^{\flat}+ \f^{\flat}\big(I_1^Y\cdot,I_1^Y\cdot\big)\big]$ -- notice $\varpi_{\f^{\flat}}$ is not closed in general. 
  Pushing-forward by $\beth$, this amounts to 
   \begin{equation*}
    \Big|dd^c_{\beth_*I_1^{Y}}\varphi 
           - \frac{1}{2}\big[\f\big(\beth_*I_1^Y\cdot,\cdot\big)- \f\big(\cdot,\beth_*I_1^Y\cdot\big)\big]\Big|_{\f^{\flat}}=O(R^{-2}). 
   \end{equation*}
  Now $dd^c_{\beth_*I_1^{Y}}\varphi=d\beth_*I_1^{Y}d\varphi=\omega_{\f}+ d\jmath d\varphi $, 
  where $\omega_{\f}= \f(I_1\cdot,\cdot)$ and $\jmath=\beth_*I_1^{Y}-I_1$. 
  Let us estimate $ d\jmath d\varphi$;  
  by Lemma \ref{lem_beth} and by the analogue of the fact raised in the proof of Proposition \ref{prop_approx_omega1Y} for $(1,1)$-tensors, 
  for all $\ell\geq0 $, $\big|(\nabla^{\f})^{\ell}\jmath\big|_{\f}= O(R^{-3})$, 
  whereas $\big|(\nabla^{\f})^{\ell}\varphi\big|_{\f}= O(R^{2-\ell})$; 
  therefore $\big| d\jmath d\varphi\big|_{\f}=O(R^{-2})$. 
  On the other hand, still from $\beth_*I_1^{Y}=I_1+\jmath$,  
  $\f\big(\beth_*I_1^Y\cdot,\cdot\big)- \f\big(\cdot,\beth_*I_1^Y\cdot\big)= 2\omega_{\f}+\f(\jmath\cdot, \cdot)-\f(\cdot,\jmath \cdot)$. 
  The error term $\f(\jmath\cdot, \cdot)-\f(\cdot,\jmath \cdot)$ is controlled by $|\jmath|_{\f}$, which is $O(R^{-3})$. 
  We have thus proved estimate \eqref{eqn_estmte1}. 
  Thanks to the general formal formula 
   \begin{equation}  \label{eqn_nablagh}
   \nabla^{g+h}T = \nabla^g T + (g+h)^{-1}* \nabla^{g}h*T, 
  \end{equation}
 (see e.g. \cite[p.21]{gv}) for any metrics $g$ and $g+h$ ($h$ is thus seen as a perturbation) and any tensor $T$, 
  with Lemma \ref{lem_f/fflat} take $g=\f$, $g+h=\f^{\flat}$ and $T$ the tensor in play, 
  we prove with the same techniques an estimate similar to \eqref{eqn_estmte1} up to order 2, that is: 
   \begin{equation*}
    \big|(\nabla^{\f^{\flat}})^{\ell}\big(dd^c_{\beth_*I_1^{Y}}\varphi - \varpi_{\f^{\flat}}\big)\big|_{\f^{\flat}}
             =O(R^{-2}), 
   \end{equation*}
 for $\ell=1,2$. 
 If therefore $K$ is chosen large enough, 
 and taking moreover the contribution of $\Psi_{\mxd}$ into account, 
 $\omega_1^{Y}+dd^c_{I_1^{Y}}\big(\varphi^{\flat}+\Psi_{\mxd}-K\big)$ is well-defined and is an $I_1^Y$-Kähler form, 
 and is equal to $\big(\omega_1^{Y}-dd^c_{I_1^{Y}}\Psi_{\mxd}\big)+\varpi_{\f^{\flat}}$ 
 up to a $O(R^{-2})$ error at orders 0,1 and 2 for $\f^{\flat}$; 
 \textit{we fix such a $K$ once for all}. 
 
~

 We know deal with the summand $-\chi\big((r-r_0)^{\beth}\big)\widetilde{\psi}_{\euc}$ of $\Phi_m^{\flat}$, 
 which is meant to kill the ALE part the Kähler form we reached, 
 or equivalently of the $I_1^Y$-hermitian form $\big(\omega_1^{Y}-dd^c_{I_1^{Y}}\Psi_{\mxd}\big)+\varpi_{\f^{\flat}}$. 
 As before, there are two issues here: 
 the positivity of the resulting $I_1^Y$-(1,1) form on $Y$, 
 and its asymptotics. 
 
 About the latter, since we are only looking at what happens near infinity, 
 notice they are independent of $r_0$ and $\beta$. 
 Indeed, for any value of these parameters, and provided that $r_0$ is chosen much larger than $K$, 
 we have on $r\geq r_0$, by definition of $\Phi_m^{\flat}$, 
  \begin{equation*}
   \omega_1^{Y}+dd^c_{I_1^{Y}}\Phi_m^{\flat}=\big(\omega_1^{Y}-dd^c_{I_1^{Y}}\Psi\big)
                                              + dd^c_{\beth_*I_1^{Y}}\varphi,
  \end{equation*}
 with $\Psi$ that of Proposition \ref{prop_approx_omega1Y}; 
 the parenthesis in the right-hand side is thus $O(R^{-2})$ for $\f$ by this proposition, 
 and again this holds for $\f^{\flat}$ by Lemma \ref{lem_f/fflat}. 
 We have already dealt with the asymptotics of $dd^c_{\beth_*I_1^{Y}}\varphi$ in the previous step, 
 and know they verify they announced estimates, 
 i.e. the metric associated to this Kähler form via $I_1^Y$ differs from $\f^{\flat}$ up to order 2 by a $O(R^{-2})$ error. 
 
 We are therefore left with the positivity assertion, 
 which has to be proved carefully since we essentially have to subtract a metric to another one, 
 hence our use of the to parameters $r_0$ and $\beta$. 
 This is however perfectly similar to what is done in \cite[p.21]{auv1}, 
 and consists here into the following:
  \begin{itemize}
   \vspace{-3pt}
   \item take $r_0$ so that on $r\geq r_0$, 
         $dd^c_{I_1^{Y}}\big(\varphi^{\flat}+\Psi_{\mxd}-K\big)\geq \frac{1}{2} \varpi_{\f^{\flat}}$, 
   \item consider the remaining part $\omega_1^Y-dd^c_{I_1^{Y}}\big[\chi\big((r-r_0)^{\beta}\big)\widetilde{\psi}_{\euc}\big]$, 
         which can be rewritten as
         \begin{equation*}
          \chi\big((r-r_0)^{\beta}\big)\big(\omega_1^Y-dd^c_{I_1^{Y}}\widetilde{\Psi}_{\euc}\big)
          +\big(1-\chi\big((r-r_0)^{\beth}\big)\big) + R_{\beta}, 
         \end{equation*}
         where $R_{\beta}$ vanishes outside of $\{r_0\leq r\leq r_0+1\}$, 
         and $|R_{\beta}|_{\e}\leq C\beta$ for some constant $C=C(r_0)$; 
         
   \item $\chi\big((r-r_0)^{\beta}\big)\big(\omega_1^Y-dd^c_{I_1^{Y}}\widetilde{\Psi}_{\euc}\big)$, 
         which is $O(r^{-4})$ for $\e$, that is $O(r^{-2})$ i.e. $O(R^{-1})$ for $\f$ or $\f^{\flat}$, 
         and vanishes outside $\{r\geq r_0\}$; 
         one can thus fix $r_0$ large enough so that it this 2-form is $\geq -\frac{1}{6}\varpi_{\f^{\flat}}$ 
         everywhere on $Y$
         
   \item fix finally $\beta$ so that$|R_{\beta}|_{\e}$ small enough to say that 
         $|R_{\beta}|_{\f^{\flat}}\leq \frac{1}{6}\varpi_{\f^{\flat}}$ where it may not vanish, 
         i.e. on $\{r_0\leq r\leq r_0+1\}$; 
         this way 
         $\omega_1^Y-dd^c_{I_1^{Y}}\big[\chi\big((r-r_0)^{\beta}\big)\widetilde{\psi}_{\euc}\big]
          \geq -\frac{1}{6}\varpi_{\f^{\flat}} -\frac{1}{6}\varpi_{\f^{\flat}}
          = -\frac{1}{3}\varpi_{\f^{\flat}}$, 
         and therefore $\omega_1^Y+dd_{I_1^Y}\Phi_m^{\flat}
                        \geq \frac{1}{2}\varpi_{\f^{\flat}}-\frac{1}{3}\varpi_{\f^{\flat}}
                        =\frac{1}{2}\varpi_{\f^{\flat}}-\frac{1}{3}\varpi_{\f^{\flat}}$ 
         on $\{r\geq r_0\}$, 
         whereas it is equal to $\omega_1^Y+dd^c_{I_1^{Y}}\big(\varphi^{\flat}+\Psi_{\mxd}-K\big)\geq0$ 
         on $Y\backslash\{r\geq r_0\}$, hence the desired positivity assertion.
  \end{itemize}
 
 The last part of the statements concerns volume forms, 
 and is a direct consequence of the estimates on the metrics, after seeing that 
 (on $\R^4$, say; recall that ${\Phi_Y}_*\Omega_Y=\Omega_{\e}$)
 $\vol^{\f^{\flat}}-\Omega_Y=\beth^*\vol^{\f}-\Omega_{\e}=\beth^*\Omega_{\e}-\Omega_{\e}$, 
 which can be written as $\vareps\Omega_{_e}$ with $\big|(\nabla^{\e})^{\ell}\vareps\big|_{\e}=O(r^{-8})$ ($\ell\geq0$)
 by Lemma \ref{lem_beth}. 
 This converts into $\big|(\nabla^{\f^{\flat}})^{\ell}\vareps\big|_{\f^{\flat}}=O(R^{-4})$, $\ell\geq0$, 
 which is better than wanted. 
 \cqfd

 \subsection{Corrections on the glued metric}  \label{section_corrections}
  \subsubsection{A Calabi-Yau type theorem}    \label{subsec_CYALF}
  We want to correct our $I_1^Y$-Kähler metric $g_m$ from Proposition \ref{prop_gluing} into a Ricci-flat Kähler metric. 
  For this it is sufficient to correct it into an $I_1^Y$-Kähler metric with volume form $\Omega_Y$, 
  since this is the volume of the $I_1^Y$-Kähler metric $g_Y$ 
  -- and once the complex structure is fixed, it is well-known that the Ricci tensor of a Kähler metric
  depends only on its volume form. 
  As suggested by program ending to Theorem \ref{thm_gRF}, 
  at the level of $I_1^Y$-Kähler forms, we want to stay in the same class; 
  in other words, we are looking for the $I_1^Y$-complex hessian of some function to be our the desired correction.   
  
  The tool we are willing to use to determine this function is the following ALF Calabi-Yau type theorem:  
   \begin{thm}  \label{thm_CYALF}
    Let $\alpha,\beta\in (0,1)$ and let $(\EuScript{Y},g_{\EuScript{Y}},J_{\EuScript{Y}}, \omega_{\EuScript{Y}})$ 
    an ALF Kähler 4-manifold of dihedral type of order $(3,\alpha,\beta)$. 
    Let $f$ a smooth function in $C^{3,\alpha}_{\beta+2}(\EuScript{Y},g_{\EuScript{Y}})$.  
    Then there exists a smooth function $\varphi\in  C^{5,\alpha}_{\beta}(\EuScript{Y},g_{\EuScript{Y}})$ 
    such that $\omega_{\EuScript{Y}}+dd^c_{J_{\EuScript{Y}}}\varphi$ is Kähler, 
    and $(\omega_{\EuScript{Y}}+dd^c_{J_{\EuScript{Y}}}\varphi)^2=e^f\omega_{\EuScript{Y}}^2$. 
   \end{thm}
  The weighted Hölder spaces of this statement follow the classical definition, 
  and are the analogues of those defined in next paragraph for $g_m$ on $Y$. 
  Let us now make the following remark: 
  since we want to construct a metric with volume form $\Omega_Y$, 
  this is tempting to take $f=\log\big(\frac{\Omega_Y}{\vol_{g_m}}\big)$ to apply Theorem \ref{thm_CYALF}. 
  But so far we only control such an $f$ up to two derivatives (see Proposition \ref{prop_gluing}, estimates \eqref{eqn_estvolform}). 
  
  The other issue is that "ALF Kähler manifold of dihedral type" can be taken so as to mean that outside a compact subset, 
  $\EuScript{Y}$ is diffeomorphic to the complement of a ball in $\R^4/\mathcal{D}_k$, 
  and that one can choose the diffeomorphism $\Phi_{\EuScript{Y}}$
  between infinities of $\EuScript{Y}$ and $\R^4/\mathcal{D}_k$ such that 
  for all $\ell=0,\dots,3$, 
  $\big|(\nabla^{g_{\EuScript{Y}}})\big({\Phi_{\EuScript{Y}}}_*g_{\EuScript{Y}}-\f\big)\big|_{g_{\EuScript{Y}}}
        =O(\rho^{-\beta-\ell})$, 
  in addition with a similar statement of the $\alpha$-Hölder derivative of $({\Phi_{\EuScript{Y}}}_*g_{\EuScript{Y}}-\f)$, 
  and an analogous statement on the complex structures ${\Phi_{\EuScript{Y}}}_*J^{\EuScript{Y}}$ and $I_1$. 
  Here again, a reading of Proposition \ref{prop_gluing} indicates us that the asymptotics at our disposal 
  do not allow us to take $\Phi_{\EuScript{Y}}=\Phi_Y$. 
  
  We remedy to those problems as follows. 
  First we correct $g_m$ into an $I_1^Y$-Kähler metric with volume form $\Omega_Y$ 
  -- which is nothing but a \textit{Ricci-flat} $I_1^Y$-Kähler metric -- \textit{outside a compact subset of $Y$}, 
  which gives us an $f$ with compact support; 
  then we put this corrected metric into so-called \textit{Bianchi gauge} with respect to ${\Phi_Y}^*\f^{\flat}$, 
  which corresponds to correct $\Phi_Y$ itself so as to fit into the definition of an ALF Kähler manifold of dihedral type. 
  
  We conclude this paragraph with a word on Theorem \ref{thm_CYALF}. 
  Strictly speaking, it is proved in \cite{auv1} -- see also the seminal works \cite{yau, tian-yau1 , tian-yau2}, or the more recent \cite{hein} -- 
  with Hölder regularity of infinite order instead of $(3,\alpha)$. 
  A careful reading of the proof though shows the statement we propose here is also valid.

 \subsubsection{Ricci-flatness outside a compact subset}   \label{subsec_RFnearinfinity}
  To correct $g_m$ into an $I_1^Y$-Kähler metric with volume form $\Omega_Y$  outside a compact subset of $Y$, 
  we proceed as in \cite[\S 2.3]{auv1}, 
  with the loss of one order in the asymptotics, both in the precision and the order of differentiation. 
  Namely we start with defining on $Y$ \textit{weighted Hölder spaces}
   \begin{equation}  \label{eqn_dfwghtHspce}
    C^{\ell,\alpha}_{\delta}(Y,g_m)
     :=\big\{f\in C^{k,\alpha}_{loc}\big| \|f\|_{C^{\ell,\alpha}_{\delta}}<\infty\big\}, 
   \end{equation}
  for $k\in\N$, $\alpha\in(0,1]$, $\delta\in\R$, and where
   \begin{equation*}
    \begin{aligned}
     \|f\|_{C^{\ell,\alpha}_{\delta}}&:= 
          \big\|R^{\delta}f \big\|_{C^0} 
          \cdots +  \big\|R^{\delta+\ell}(\nabla^{g_m})^{\ell} f \big\|_{C^0}\\
          &+\sup_{\stackrel{(x,y)\in Y,}{d_{g_m}(x,y)<\inj_{g_m}} }
          \bigg|\max\big(R(x)^{\ell+\alpha+\delta},R(y)^{\ell+\alpha+\delta}\big)
               \frac{(\nabla^{g_m})^{\ell} f(x)-(\nabla^{g_m})^{\ell} f(y)}{d_{g_m}(x,y)^{\alpha}}\bigg|_{g_m}
    \end{aligned}
   \end{equation*}
  with $R$ a positive smooth extension of ${\Phi_Y}^*R$ on $Y$, 
  and $C^0$-norms of the tensors computed with $g_m$. 
  
  We then state the following, indicating the type of functions which can help correcting $\omega_m$ in the sense raised above:
   \begin{prop}  \label{prop_RFnearinfinity}
    Fix $(\alpha_1,\delta_1)\in(0,1)^2$ such that $\alpha_1+\delta_1<1$, $\delta_1>\frac{3}{4}$. 
    There exists a smooth function $\psi\in C^{2,\alpha_1}_{\delta_1-1}\cap C^{3,\alpha_1}_{\delta_1-2}$ 
    such that $\omega_{\psi}:=\omega_m+ dd^c_{I_1^Y}\psi$ is Kähler for $I_1^Y$, 
    and such that $\frac{1}{2}\omega_{\psi}^2=\Omega_Y$ outside a compact set. 
   \end{prop}
   \prf. 
    It is completely analogous to the proof of Proposition 2.7 of \cite{auv1}, 
    namely if $\chi$ is a cut-off function as in Proposition \ref{prop_gluing} and setting $\chi_{R_1}=\chi(R-R_1)$, 
    one solves the problem $\big(\omega_m+dd^c_{I_1^Y}\psi\big)^2=(1-\chi_{R_1})\omega_m^2+2\chi_{R_1}\Omega_Y$ 
    for $R_1$ large enough with help of the implicit functions theorem. 
    This is manageable, since: 
     \begin{itemize}
      \vspace{-3pt}
      \item $\omega_m^2-\big((1-\chi_{R_1})\omega_m^2+2\chi_{R_1}\Omega_Y\big)= \chi_{R_1}(\omega_m^2-2\Omega_Y)$, 
            and $\Big\|\chi_{R_1}\frac{\omega_m^2-2\Omega_Y}{\Omega_Y}\Big\|_{C^{k,\alpha_1}_{\delta_1}}$ 
            tends to 0 as $R_1$ goes to $\infty$ thanks to estimates \eqref{eqn_estvolform} for $k=0,1$; 
      \item the linearisation of the operators
            $C^{2+\vareps,\alpha_1}_{\delta_1-1-\vareps}\to C^{\vareps,\alpha_1}_{\delta_1+1-\vareps}$, 
            $\vareps=0,1$, $\psi\mapsto (\omega_m+dd^c_{I_1^Y}\psi)^2/\omega_m^2$, at $\psi=0$, 
            are the scalar Laplacians 
            $\Delta_{g_m}:C^{2+\vareps,\alpha_1}_{\delta_1-1-\vareps}\to C^{\vareps,\alpha_1}_{\delta_1+1-\vareps}$.  
            These are surjective, with kernel reduced to constant functions, 
            according to the appendix of \cite{bm}, 
            and using that $(Y,g_m)$ is asymptotically a circle fibration over $\R^3/\pm$. 
     \end{itemize}
    Once $R_1$ is chosen large enough so that one can apply the implicit function theorem simultaneously and that
    $\psi$ is fixed in $C^{2,\alpha_1}_{\delta_1-1}\cap C^{3,\alpha_1}_{\delta_1-2}$ so that 
    $\big(\omega_m+dd^c_{I_1^Y}\psi\big)^2=(1-\chi_{R_1})\omega_m^2+2\chi_{R_1}\Omega_Y$, 
    the last point to be checked is the positivity of $\omega_{\psi}:=\omega_m+dd^c_{I_1^Y}\psi$. 
    Since $dd^c_{I_1^Y}\psi=O(R^{-\delta_1})$, 
    $\omega_{\psi}$ is asymptotic to $\omega_m$, hence positive near infinity. 
    Since its determinant $\frac{(1-\chi_{R_1})\omega_m^2+2\chi_{R_1}\Omega_Y}{\omega_m^2}$ relatively to $\omega_m$  never vanishes, 
    it is positive on the whole $Y$. 
    The smoothness of $\psi$ is local. 
    \cqfd

   \subsubsection{Bianchi gauge for $\omega_{\psi}$}   \label{subsec_Bgge}
    \noindent
    \textit{Motivation.}
     We are know willing to deduce regularity statements on $g_{\psi}$, using its Ricci-flatness near infinity.  
     Nonetheless this cannot be done immediately.  
     The reason is that the Ricci-flatness condition is invariant under diffeomorphisms, 
     and consequently \textit{the linearisation of the Ricci tensor seen as an operator on metrics is not elliptic}, 
     which is problematic when looking for regularity.  
          
     One can however bypass this difficulty by fixing \textit{a gauge}, 
     which infinitesimally corresponds to looking at metrics with good diffeomorphisms. 
     We introduce the diffeomorphisms we shall work with in next paragraph; 
     then the gauge is  fixed, and regularity is deduced from this process (Propositions \ref{prop_gge} and \ref{prop_gge2}). 
     Notice that the Ricci-flatness of $\omega_{\psi}$ is an indispensable prerequisite in this procedure, 
     since the gauge alone is not enough in general to obtain the regularity statement we are seeking here. 
     
     ~
   
    \noindent
    \textit{ALF diffeomorphisms of $\C^2$.} 
     The class of diffeomorphism we work with to perform our gauge enters into the following definition; 
     we define the dual frames $(e_0^{\flat}, \dots, e_3^{\flat})$ 
     and $\big((e_0^{*})^{\flat}, \dots, (e_3^{*})^{\flat}\big)$ as the pull-backs by $\beth$ of the frames 
     $(e_i)$ and $(e_i^*)$ defined in section \ref{subsection_TN} by formulas \eqref{eqn_ei}, \eqref{eqn_ei*}. 
      \begin{df}  \label{df_TNdiffeos}
       Let $(\ell,\alpha)\in\N^*\times(0,1)$, and let $\nu>-1$. 
       We denote by $\diff^{\ell,\alpha}_{\nu}$ the class of diffeomorphisms $\phi$ of $\C^2$ such that:
        \begin{itemize}
         \vspace{-3pt}
         \item $\phi$ has regularity $(\ell,\alpha)$;
         \item there exists a constant $C$ such that for any $x\in \C^2$, 
               $d_{\f^{\flat}}\big(x,\phi(x)\big)\leq C\big(1+R(x)\big)^{-\nu}$;
         \item let $R_0\geq 1$ such that for any $x\in\{R\geq R_0\}$, $d_{\f^{\flat}}\big(0,\phi(x)\big)\geq1$. 
          Denote by $\gamma_x:[0,1]\to\C^2$ a minimizing geodesic for $\f^{\flat}$ joining $\phi(x)$ to $x$ 
          and by $p_{\gamma_x}$ the parallel transport along $\gamma_x$. 
          Consider the maps $\phi_{ij}: \{R\geq R_0\} \to \R$ given by 
       \begin{equation*}
        \phi_{ij}(x)=(e_{i}^{\flat})^*\big((T_x\phi\circ p_{\gamma_x}(1)-\id_{\C^2})(e_{j}^{\flat})\big),  
        \qquad i,j=0,\dots,3, 
       \end{equation*}
      and extend them smoothly in $\{R\leq R_0\}$. 
      We then ask: $\phi_{ij}\in C^{\ell-1,\alpha}_{\nu+1}(\C^2, \f^{\flat})$.  
    \end{itemize}
   We endow $\diff^{\ell,\alpha}_{\nu}$ with the natural topology.
   
   We moreover denote by $\big(\diff^{\ell,\alpha}_{\nu}\big)^{\mathcal{D}_k}$ 
   the set of diffeomorphisms of $\diff^{\ell,\alpha}_{\nu}$ commuting with the action of $\mathcal{D}_k$. 
      \end{df}
      
   The Hölder spaces are those defined for $\f^{\flat}$ on $\C^2$, in the same way as those of defining equation \eqref{eqn_dfwghtHspce}. 
   Notice that we authorise the distance between a point and its image to go to $\infty$; 
   nonetheless the rate of blow-up we allow makes clear the existence of the $R_0$ of the third item, since $\nu>-1$.

     ~
   
  \noindent
  \textit{Diffeomorphisms as Riemannian exponential maps. }
   We now "parametrise" our diffeomorphisms via vector fields:
    \begin{lem}   \label{lem_expdiffeo}
     There exists a neighbourhood $\mathscr{V}^{\ell,\alpha}_{\nu}$ of 0 in $C^{\ell,\alpha}_{\nu}(\C^2, \f^{\flat})$ such that
     for any $Z$ in that neighbourhood, the map $\phi_Z : x\longmapsto \exp^{\f^{\flat}}_x\big(Z(x)\big)$ 
     is in $\diff^{\ell,\alpha}_{\nu}$. 
    %
   \end{lem}
  The weighted spaces of vector fields are defined analogously to that of the previous paragraph, or equivalently:
  $Z\in C^{\ell,\alpha}_{\nu}(\C^2, \f^{\flat})$ if and only if $Z\in C^{\ell,\alpha}_{loc}$ and 
  $\chi(R)(e_i^{\flat})^*(Z)\in C^{\ell,\alpha}_{\nu}(\C^2, \f^{\flat})$, $i=1,\dots,3$ (with $\chi$ a cut-off function as in \ref{prop_gluing}). 
  A similar statement with $\mathcal{D}_k$-invariant vector fields, 
  and diffeomorphisms commuting with the action of $\mathcal{D}_k$ of course holds. 
  We simply call $\big(\mathscr{V}^{\ell,\alpha}_{\nu}\big)^{\mathcal{D}_k}$ the neighbourhood of 0 in $\big(C^{\ell,\alpha}_{\nu}(\C^2, \f^{\flat})\big)^{\mathcal{D}_k}$, 
  the $\mathcal{D}_k$-invariant vector fields of $C^{\ell,\alpha}_{\nu}(\C^2, \f^{\flat})$. 
  Notice finally that for a genuine parametrisation, 
  we would also need the surjectivity and the injectivity of $Z\mapsto \phi_Z$ from $\mathscr{V}^{\ell,\alpha}_{\nu}$ onto its image. 
  We do not need however this degree of precision, since as seen in Proposition \ref{prop_gge} below, 
  it is enough for us to realise sufficiently many diffeomorphisms of $\diff^{\ell,\alpha}_{\nu}$ under the shape $\phi_Z$. 

 ~

  \prf. The regularity assertions are rather standard. 
   We shall nonetheless pay a particular attention to the fact that we authorise vector fields blowing up at infinity,  
   when verifying the injectivity of $\phi_Z$ for a given $Z$ 
   close to 0 in $C^{\ell,\alpha}_{\nu}$; the key is the decay of the derivatives of $Z$ at infinity. 
   Suppose $(\ell,\alpha)=(1,0)$ to fix ideas. 
   For the injectivity of $\phi_Z$ with fixed $Z\in C^{1,0}_{\nu}$ and $\|Z\|_{C^{1,0}_{\nu}}\leq1$ say, 
   we claim that there exists a constant $C$ independent of $Z$ such that for any triple $(x,y,z)$ such that $\phi_Z(x)=\phi_Z(y)=:z$, 
    \begin{equation*}
     d_{\f^{\flat}}(x,y) \leq C\big(1+R(z)\big)^{-4-3\nu}\|Z\|_{C^{1,0}_{\nu}} d_{\f^{\flat}}(x,y), 
    \end{equation*}
   from which the injectivity of $\phi_Z$ follows at once provided $\|Z\|_{C^{1,0}_{\nu}}$ is small enough. 
   We reach this claim thanks to the estimate $\big|\riem^{\f^{\flat}}\big|=O(R^{-3})$, as follows. 
   For $x,y$ as in the claim, call respectively $\gamma_x$ and $\gamma_y$ the geodesics $t\mapsto \exp_x^{\f^{\flat}}\big(tZ(x)\big)$ 
   and $t\mapsto \exp_y^{\f^{\flat}}\big(tZ(y)\big)$, and denote by $p_{\gamma_x}$, $p_{\gamma_y}$ the attached parallel transports. 
   Using \cite[Prop. 6.6]{bk}, control first $d_{\f^{\flat}}(x,y)$ by 
   $\big|p_{\gamma_x}(1)\big(Z(x)\big)-p_{\gamma_y}(1)\big(Z(y)\big)\big|_{\f^{\flat}}\big(1+R(z)\big)^{-3-2\nu}$.
   Then control $\big|p_{\gamma_x}(1)\big(Z(x)\big)-p_{\gamma_y}(1)\big(Z(y)\big)\big|_{\f^{\flat}}$  by
   $d_{\f^{\flat}}(x,y)\big(1+R(z)\big)^{-1-\nu}\|Z\|_{C^{1,0}_{\nu}}$;  
   for this interpolate between $\gamma_{x}$ and $\gamma_y$ 
   by $\gamma_s(t):=\exp^{\f^{\flat}}_{\alpha(s)}\big[tZ\big(\alpha(s)\big)\big]$, 
   where $\alpha$ is a minimizing geodesic for $\f^{\flat}$ joining $x$ and $y$. 
   This is where one uses the estimates on the derivatives of $Z$. 
    \cqfd
  
 ~

   \noindent 
   \textit{The gauge.} 
   Denote by $B^h=\delta^{h}+\frac{1}{2}\tr^h$ the Bianchi operator associated to any smooth metric $h$ on $\R^4$. 
   The gauge process now states as:
   \begin{prop} \label{prop_gge}
    Let $(\alpha_2,\delta_2)\in(0,1)^2$ such that $\alpha_2<\alpha_1$, $\tfrac{3}{4}<\delta_2<\delta_1$, with $(\alpha_1,\delta_1)$ 
    fixed in Proposition \ref{prop_RFnearinfinity}. 
    There exists a smooth diffeomorphism $\phi\in  \big(\diff^{1,\alpha_2}_{\delta_2-1}\big)^{\mathcal{D}_k}$ 
    hence descending to $\R^4/\mathcal{D}_k$ such that
     \begin{equation}  \label{eq_gge2}
      B^{\phi^*\f^{\flat}}\big((\Phi_Y)_*g_{\psi}\big) = 0
     \end{equation}
    near infinity on $\C^2$,  where $g_{\psi}$ stands for the $I_1^Y$-Kähler metric associated to the Kähler form $\omega_{\psi}$ of Proposition \ref{prop_RFnearinfinity}. 
    As a consequence, $\f^{\flat}-(\phi\circ\Phi_Y)_*g_{\psi}\in C^{1, \alpha_2}_{\delta_2}\big(X, \f^{\flat}\big)$.  
  \end{prop}
   \prf. 
    Fix $(\alpha_2,\delta_2)$ as in the statement, and consider the map
    \begin{equation*}
     \begin{aligned} 
      \Xi\,:\, (\mathscr{V}^{2,\alpha_2}_{\delta_2-1})^{\mathcal{D}_k}\times \met^{1,\alpha_2}_{\delta_2}(\f^{\flat})^{\mathcal{D}_k}\, 
                                                                                                               &\longrightarrow \, C^{0,\alpha_2}_{\delta_2+1}(T^*\C^2,\f^{\flat})^{\mathcal{D}_k}\   \\
                              (Z,g)      \quad  \,\,\, \,  \quad    \qquad      \quad        &\longmapsto     \,           B^{\phi_Z^*\f^{\flat}}(g),
   \end{aligned}
  \end{equation*}
  where $\met^{1,\alpha_2}_{\delta_2}(\f^{\flat})^{\mathcal{D}_k}$ denotes 
  the set of $\mathcal{D}_k$-invariant metrics $g$ on $\R^4$ such that $g-\f^{\flat}\in C^{\infty}_{\delta_2}\big(X, \f^{\flat}\big)$   
  -- we shall see further why we anticipate the action of $\mathcal{D}_k$ at this point of our discussion. 

  We would like to solve the equation 
   \begin{equation}  \label{eqn_gge}
    B^{\phi_Z^*\f^{\flat}}(g_{\psi})=0,  \qquad \text{i.e} \qquad \Xi(Z,g_{\psi})=0, 
   \end{equation}
  and for this use the implicit function theorem near $(0,\f^{\flat})$, 
  since the differential of $\Xi$ with respect to its first argument is $(\nabla^{\f^{\flat}})^*\nabla^{\f^{\flat}}$, 
  which as we shall see enjoys isomorphism properties. 
  Forgetting that $g_m$ may not be defined via $\Phi_Y$ on the whole $\C^2$, 
  if we are to do so nonetheless, we need to make $g_m$ arbitrarily close to $\f^{\flat}$ in $C^{\infty}_{\delta_2}\big(X, \f^{\flat}\big)$. 
  Since we only want equation \eqref{eqn_gge} to be solved near infinity, 
  instead of $g_m$ we consider, 
  for $\chi$ a cut-off function as in Proposition \ref{prop_gluing}, 
  the metric
   \begin{equation*}
    g_{R_2}:= \chi(R-R_2)g_{\psi} + \big(1-\chi(R-R_2)\big)\f^{\flat},  
   \end{equation*}
  which makes sense via $\Phi_Y$ on the whole $\C^2$ provided $R_2$ is large enough; 
  since $g_{\psi}$ is $C^{1,\alpha_1}_{\delta_1}$ close to $\f^{\flat}$ at infinity, 
  we have that $\|g_{R_2}-\f^{\flat}\|_{C^{1,\alpha_2}_{\delta_2}}(\f^{\flat})$ goes to 0 when $R_2$ goes to $\infty$. 
  We are thus left with checking the isomorphism assertion on 
   \begin{equation*}
    (\nabla^{\f^{\flat}})^*\nabla^{\f^{\flat}}=\frac{\partial\Xi}{\partial Z}\Big|_{(0,\f^{\flat})} : 
     (C^{2,\alpha_2}_{\delta_2-1})^{\mathcal{D}_k}\longrightarrow (C^{0,\alpha_2}_{\delta_2+1})^{\mathcal{D}_k}
     \cong \big(C^{0,\alpha_2}_{\delta_2+1}(T^*\C^2)\big)^{\mathcal{D}_k}
   \end{equation*}
  where the isomorphism on the right is just the duality for $\f^{\flat}$. 
  We shall moreover replace $\f^{\flat}$ by $\f$ (and the weighted spaces subsequently) since these are diffeomorphic to each other. 
  Now surjectivity follows from that of $(\nabla^{\f})^*\nabla^{\f}$ between $C^{2,\alpha_2}_{\delta_2-1}(T\C^2,\f)$ and $C^{0,\alpha_2}_{\delta_2+1}(T\C^2,\f)$, 
  which amounts by the theory of self-adjoint operators on weighted spaces to the injectivity of this operator on $C^{0,\alpha_2}_{3-(\delta_2+1)}=C^{0,\alpha_2}_{2-\delta_2}$, 
  and is seen in the proof of Proposition 2.11 in \cite{auv1} on the bigger space $C^{0,\alpha_2}_{1-\delta_2}$; 
  this relies on an easy ellipticity argument, 
  and an integration by parts using than in ALF geometry, 
  a sphere of large radius $t$ has its volume in $t^2$. 

  The slightly newer point here is the injectivity of 
   $(\nabla^{\f})^*\nabla^{\f}: (C^{2,\alpha_2}_{\delta_2-1})^{\mathcal{D}_k}\to (C^{0,\alpha_2}_{\delta_2+1})^{\mathcal{D}_k}$; 
  this is also where the invariance under the action of $\mathcal{D}_k$ is needed. 
  Let thus $v\in (C^{2,\alpha_2}_{\delta_2-1})^{\mathcal{D}_k}$ such that $(\nabla^{\f})^*\nabla^{\f}v=0$. 
  Forget momentarily about the $\mathcal{D}_k$-invariance of $v$, 
  and write it $\sum_{i=0}^3 v^i e_i$, on $\{R\geq 1\}$, say, with $(e_i)$ the frame defined by \eqref{eqn_ei}; 
  each $v^i$ is thus in $C^{2,\alpha_2}_{\delta_2-1}$ near infinity. 
  Now since $(\nabla^{\f})^{\ell}e_i=O(R^{-1-\ell})$, $\ell=1,2,3$ and because $(\nabla^{\f})^*\nabla^{\f}(v^ie_i)$ is equal to $(\Delta_{\f}v^i)e_i$ plus a linear combination of 
  the $\nabla^{\f}_{e_j} v^i \nabla^{\f}e_k$ and the $v^i(\nabla^{\f})_{e_j,e_k}e_i $, we get that each $v^i$ is the solution of a Dirichlet problem
   \begin{equation*}
    \left\{\begin{aligned}
           \Delta_{\f}v^i  &= w^i \in C^{0,\alpha_2}_{\delta_2+2}\quad \text{on }\{R\geq 1\}, \\
           v^i|_{\{R= 1\}}&\in C^{2,\alpha_2}.
          \end{aligned}
    \right.
   \end{equation*}
  Recall that $\Delta_{\f}: C^{2,\alpha_2}_{\delta_2-1}\to C^{0,\alpha_2}_{\delta_2+1}$ is surjective, with kernel reduced to the constants 
  (see the proof of Proposition \ref{prop_RFnearinfinity}), 
  and that $\Delta_{\f}: C^{2,\alpha_2}_{\delta_2}\to C^{0,\alpha_2}_{\delta_2+2}$ is an isomorphism (see e.g. \cite[App.]{bm}); 
  those properties transfer to Dirichlet problems to tell us that each $v^i$ writes $c_i+u_i$, with $c_i$ a constant and $u_i\in C^{2,\alpha_2}_{\delta_2}$. 
  Thus $v$ is asymptotic to $\sum_ic_ie_i$; but $v$ is $\mathcal{D}_k$-invariant, whereas $\tau^*e_i=-e_i$ for any $i$. 
  This forces the $c_i$ to be 0, and as a result $v\in C^{2,\alpha_2}_{\delta_2}$. 
  Finally, we know that $(\nabla^{\f})^*\nabla^{\f}$ on this latter space, so $v=0$; 
  in other words, the action of $\mathcal{D}_k$ allows us to say that 0 is no more a critical weight for $(\nabla^{\f})^*\nabla^{\f}$ on vector fields. 

  The smoothness of $Z$, and therefore that of $\phi_Z$, is purely local. 
  \cqfd


~

 \noindent
 \textit{Regularity of $g_{\psi}$.} 
  We conclude this paragraph by the following statement, which finally allows us to apply Theorem \ref{thm_CYALF}: 
  \begin{prop}  \label{prop_gge2}
   With the same notation as in Proposition \ref{prop_gge}, 
   $\f^{\flat}-(\phi\circ\Phi_Y)_*g_{\psi}\in C^{3,\alpha_2}_{4\delta_2-3}\big(X, \f^{\flat}\big)$ near infinity, 
   and in particular, $\big|\riem^{g_{\psi}}\big|_{g_{\psi}}=O(R^{-2-(4\delta_2-3)})$.     
  \end{prop}
   \prf. 
    The assertion on the curvature of $g_{\psi}$ directly follows from the estimate stated on $\vareps:=g_{\psi}-\phi^*\f^{\flat}$ (or $\phi_*\vareps$), 
    and the fact that $\big|\riem^{\f^{\flat}}\big|_{\f^{\flat}}=O(R^{-3})$. 
    For the regularity statement on $\vareps$, 
    proceed as for \cite[Cor. 2.12]{auv1}, 
    with one order less in the approximation, to deduce, from the \textit{two} equations $\ric^{g_{\psi}}=0$, 
    $B^{\f^{\flat}}(\phi_*g_{\psi})=0$ which are verified near infinity: 
     \begin{equation}  \label{eqn_crl_gge}
      \frac{1}{2}\mathscr{L}_{\phi^*\f^{\flat}}\vareps + \vareps\star \partial^2\vareps 
          = Q(\vareps,\partial\vareps)
     \end{equation}
    where $\mathscr{L}_{\phi^*\f^{\flat}}$ is the Lichnerowicz laplacian of $\phi^*\f^{\flat}$, 
    the symbols $\star$ denote algebraic operations, and $Q$ is at least quadratic in its arguments, 
    and can be factorised by $\vareps\star\partial\vareps$. 
    Since $\vareps \in C^{1,\alpha_2}_{\delta_2}$ 
    -- computed with respect to $\phi^*\f^{\flat}$ or $g_{\psi}$, 
    which does not matter because of the size of the error term $\vareps$ itself  --, 
    the right-hand-side of \eqref{eqn_crl_gge} is in $C^{0,\alpha_2}_{\delta_2}\big(X, \f^{\flat}\big)$. 
    Again since $\vareps \in C^{1,\alpha_2}_{2\delta_2+1}$, 
    the linear operator $\eta\to \frac{1}{2}\mathscr{L}_{\phi^*\f^{\flat}}\eta + \vareps\star \partial^2\eta$ 
    is elliptic and one can draw for this operator weighted estimates similar to those for $\mathscr{L}_{\phi^*\f^{\flat}}$. 
    From this we deduce that $\vareps \in C^{2,\alpha_2}_{2\delta_2-1}$;  
    we conclude by repeating this argument, 
    giving us $\vareps \in C^{3,\alpha_2}_{2(2\delta_2-1)-1}=C^{2,\alpha_2}_{3\delta_2-3}$. 
   \cqfd 
    
  \subsubsection{Conclusion: proof of Theorem \ref{thm_gRF}}
   Provided that we take $\alpha=\alpha_2$ and $\beta= 4\delta_2-3$ (which is positive since $\delta_2>\frac{3}{4}$), 
   to fulfil completely the requirements of Theorem \ref{thm_CYALF}, 
   we are only left with checking that $(\phi\circ\Phi_Y)_*I_1^Y$ is also $C^{3,\alpha}_{\beta}$ 
   close to the complex structure $I_1^{\flat}:=\beth^*I_1$. 
   
   The estimate $(\phi\circ\Phi_Y)_*I_1^Y-I_1^{\flat}\in C^{0}_{\beta}$ follows easily form the decomposition 
   $(\Phi_Y)_*I_1^Y-\phi^*I_1^{\flat}=\big((\Phi_Y)_*I_1^Y-I_1\big)+(I_1-I_1^{\flat})+(I_1^{\flat}-\phi^*I_1^{\flat})$, 
   and from the estimates $\big|(\Phi_Y)_*I_1^Y-I_1\big|_{\e}=O(r^{-4})$ and 
   $\big|I_1-I_1^{\flat}\big|_{\e}=O(r^{-4})$ converted into $\big|(\Phi_Y)_*I_1^Y-I_1\big|_{\f}, \big|I_1-I_1^{\flat}\big|_{\f^{\flat}}=O(R^{-1})$, 
   and $\big|I_1^{\flat}-\phi^*I_1^{\flat}\big|_{\f^{\flat}}=O(R^{-\delta_2})$ following from $\phi\in \diff^{1,\alpha_2}_{\delta_2-1}$. 
   
   For higher order estimates, 
   remember that $g_{\psi}$ is Kähler for $I_1^Y$, and $\f^{\flat}$ for $I_1^{\flat}$. 
   It is thus enough for instance to evaluate the successive $(\nabla^{\f^{\flat}})^{\ell}\big((\phi\circ\Phi_Y)_*I_1^Y\big)$. 
   In view of formula \ref{eqn_nablagh}, 
   we thus write formally for $\ell=1$
    \begin{equation*}
     \nabla^{\f^{\flat}}\big((\phi\circ\Phi_Y)_*I_1^Y\big) 
             = \underbrace{\nabla^{\f^{\flat}}\big((\phi\circ\Phi_Y)_*I_1^Y\big)}_{=0\text{ since }g_{\psi}\text{ is }I_1^Y\text{-Kähler}}
               + (\f^{\flat})^{-1}*\nabla^{g_{\psi}}(\f^{\flat}-g_{\psi})*\big((\phi\circ\Phi_Y)_*I_1^Y\big),
    \end{equation*}
   which gives easily $\nabla^{\f^{\flat}}\big((\phi\circ\Phi_Y)_*I_1^Y\big)\in C^{0}_{\beta+1}$ in view of $(\f^{\flat}-g_{\psi})\in C^{0}_{\beta}$. 
   For $\ell\geq 2$, simply use inductively formula \eqref{eqn_nablagh}. 

   As sketched in the introduction of this section, 
   we now apply Theorem \ref{thm_CYALF}, 
   with $\big(\EuScript{Y},g_{\EuScript{Y}},J_{\EuScript{Y}}, \omega_{\EuScript{Y}}\big)=(Y,g_{\psi}, I_1^Y, \omega_{\psi})$, 
   and $f=\log\big(\frac{\Omega_Y}{\vol^{g_{\psi}}}\big)$, which is smooth and has compact support. 
   This gives us an $I_1^Y$-metric $g_{RF,m}$ on $Y$, with volume form $\Omega_Y$ and which is thus Ricci-flat, 
   and with Kähler form $\omega_{\psi}+dd^c_{I_1^Y}\varphi$ for some smooth $\varphi\in C^{5,\alpha}_{\beta}(Y,g_{\psi})$ with $\beta$ close to 1. 
   
   We need two more complex structures for Theorem \ref{thm_gRF}. 
   Recall we have two more symplectic forms coming with the ALE hyperkähler structure $\big(Y,g_Y, I_1^Y, I_2^Y, I_3^Y\big)$, 
   namely $\omega_2^Y:=g_Y\big(I_2^Y\cdot,\cdot\big)$ and $\omega_3^Y:=g_Y\big(I_3^Y\cdot,\cdot\big)$. 
   We simply define $J^Y_2$ and $J^Y_3$ as the almost-complex structures verifying 
   $g_{RF,m}\big(J_2^Y\cdot,\cdot\big)= \omega_2^Y$ and $g_{RF,m}\big(J_3^Y\cdot,\cdot\big)= \omega_3^Y$, 
   just as in \eqref{eqn_defJj}; 
   one then checks these are complex structures thanks using that the symplectic holomorphic 2-form $\omega_2^Y+i\omega_3^Y$ 
   is $g_{RF,m}$-parallel (as $(\omega_2^Y+i\omega_3^Y)\wedge(\omega_2^Y-i\omega_3^Y)= 4\Omega_Y=4\vol^{g_{RF,m}}$), 
   and that together with $I_1^Y$ they verify the quaternionic relations.  

   The cubic  decay of $\riem^{g_{RF,m}}$ comes as follows: 
   first, an over-quadratic decay is easily deduced from $(g_{\psi}-g_{RF,m})\in C^2_{\beta+2}\big(Y,g_{\psi}\big)$ and $\riem^{g_{\psi}}=O(R^{-2-\beta})$ (Proposition \ref{prop_gge2}). 
   Then a result of Minerbe \cite[Cor. A.2]{min2} asserts that we automatically end up with a cubic rate decay of the curvature.   
    \cqfd

  \subsection{Verification of the technical Lemmas \ref{lem_est_psidiff} and \ref{lem_f/fflat}}  \label{section_lemmas}
   We conclude this part by the left-over proofs of Lemmas \ref{lem_est_psidiff} and \ref{lem_f/fflat}, 
   both useful in the gluing performed in section \ref{section_gluing}. 
   Recall that on the one hand, Lemma \ref{lem_est_psidiff} is about verifying the asymptotics at different orders of a function $\psi_c$, 
   the hessian of which is meant to approximate the 2-form $\theta_2+i\theta_3$ in the Taub-NUT framework, 
   although such an approximation is likely to be vain in the euclidean setting;  
   and that on the other hand, Lemma \ref{lem_f/fflat} consists in saying that even though $\f^{\flat}=\beth^*\f$,  
   with $\beth$ a diffeomorphism of $\R^4$ better adapted to the euclidean scope, 
   the transition between $\f$ and $\f^{\flat}$ is relatively harmless. 
 
   \subsubsection{Proof of Lemma \ref{lem_est_psidiff}}   \label{subsec_lem_est_psidiff}
    \noindent
    \textit{Asymptotics of $\psi_c$ and its successive derivatives. }
     We first look at the first point of the statement of Lemma \ref{lem_est_psidiff}. 
     Since $\psi_c$ is $\mathbb{S}^1$-invariant when looked at on $\C^2$ 
     (recall that the $\mathbb{S}^1$-action on $\C^2$ is given by $\alpha\cdot(z_1,z_2)=(e^{i\alpha}z_1, e^{-i\alpha}z_2)$), 
     or in other words is a function of $y_1$, $y_2$, $y_3$ (recall in particular that $2r^2=R\cosh(4my_1)+y_1\sinh(4my_1)$, 
     following formulas \eqref{eqn_lebformulas} and the definitions of $y_1$ and $R$ given in paragraph \ref{subsection_TN}), 
     we have: $d\psi_c= \tfrac{\partial\psi_c}{\partial y_1}dy_1+ \tfrac{\partial\psi_c}{\partial y_2}dy_2+ \tfrac{\partial\psi_c}{\partial y_3}dy_3$, 
     and one can see as well this partial derivatives as functions of the $y_j$ only. 
     If we thus prove here that for any $p$, $q$, $s\geq0$ such that $p+q+s\leq4$, 
      \begin{equation}  \label{eqn_est_psidiff}
        \frac{\partial^{p+q+s}\psi_{c}}{\partial y_1^p\partial y_2^q \partial y_3^s} = O(R^{-1-q-s}),  
       \end{equation} 
     we will get the desired estimates, since we moreover know that 
     $\big|(\nabla^{\f})^{\ell}dy_j\big|_{\f}=O(R^{-1-\ell})$ for all $\ell\geq 1$ and $j=1,2,3$. 
     
     The estimate \eqref{eqn_est_psidiff} at order $0$ is immediate, since $\sinh (4my_1)= O(R^{-1}r^2)$ -- 
     this follows from the identity $2r^2=R\cosh(4my_1)+y_1\sinh(4my_1)$. 
     What is thus clearly to be seen is that each time we differentiate with respect to $y_2$ or $y_3$, we win an $R^{-1}$, 
     and each we differentiate with respect to $y_1$, we lose nothing. 
     Let us how it goes at order 1, that is when $p+q+s=1$. 
     If $p=1$ and $q=s=0$, then (near infinity, where $\chi(R)\equiv 1$): 
      \begin{equation*}
         \frac{\partial\psi_c}{\partial y_1}
               =-4 (y_2+iy_3)\Big(\frac{4m\cosh (4my_1)}{2Rr^2}- \frac{y_1\sinh (4my_1)}{2r^2R^3} 
                                - \frac{\sinh (4my_1)}{4r^4R}\frac{\partial(2r^2)}{\partial y_1} \Big)
      \end{equation*}
     and $\tfrac{\partial(2r^2)}{\partial y_1}= 2V\big(y_1\cosh(4my_1)+R\sinh(4my_1)\big) $ 
     (recall that $V=\tfrac{1+4mR}{2R}$), so that, 
     after simplifying:
       \begin{equation*}
        \frac{\partial\psi_c}{\partial y_1}
           =-4 (y_2+iy_3)\Big( \frac{1}{r^4} -\frac{1}{R^3} + \frac{1}{4r^4R}\Big), 
      \end{equation*}
     and this is $O(R^{-1})$, since $r^{-2}=O(R^{-1})$ (as $R=O(r^2)$). 
 
     If $q=1$ and $p=s=0$, then 
      \begin{equation*}
       \frac{\partial\psi_c}{\partial y_2}
               = -2\frac{\sinh (4my_1)}{r^2R} 
                 -2 (y_2+iy_3)\sinh (4my_1)\Big( \frac{y_2}{r^2R^3} + \frac{y_2\cosh (4my_1)}{2r^4R^2}\Big), 
      \end{equation*}
     since $\tfrac{\partial(2r^2)}{\partial y_2}= \tfrac{y_2}{R}\cosh (4my_1)$. 
     As $\sinh(4my_1)$ and $\cosh(4my_1)$ are $O(r^{2}R^{-1})$, we end up with 
     $\frac{\partial\psi_c}{\partial y_2}= O\big(r^2/(R^2r^2)\big)+O\big(R\cdot r^2/R\cdot (r^{-2}R^{-2}+r^2/R\cdot r^{-4}R^{-1})\big)=O(R^{-2}) $. 
     The case $s=1$ and $p=s=0$, i.e. the estimate of $\tfrac{\partial\psi_c}{\partial y_3}$, 
     is done by substituting $y_3$ to $y_2$. 
 
     In a nutshell, we win one order each time we differentiate $y_2$, $y_3$, $R$ and $r^2$ with respect to $y_2$ or $y_3$, 
     which moreover kills functions of $y_1$ such as $\sinh(4my_1)$; 
     we win one order as well when differentiating $y_2$, $y_3$ and $R$ with respect to $y_1$, 
     but this does not hold any more for $r^2$ or functions like  $\sinh(4my_1)$. 
     More formally, using explicit formulas for the $\tfrac{\partial(2r^2)}{\partial y_j}$, $j=1,2,3$, 
     we can easily prove by induction that for any $p,q,s$ there exists a polynomial $Q_{p,q,s}$ 
     of total degree $\leq (1+p+q+s)$ in their first two variables, 
     and $2+3p+2(q+s)$ in total, such that: 
      \begin{equation*}
       \frac{\partial^{p+q+s}\psi_{c}}{\partial y_1^p\partial y_2^q \partial y_3^s} 
                   = \frac{Q_{p,q,s}\big(Re^{\pm 4my_1},y_1e^{\pm4my_1}, R,y_1,y_2,y_3\big)}{(2r^2)^{1+p+q+s}R^{2(1+p+q+s)}}  
      \end{equation*}
     (for instance, 
     $Q_{1,0,0}\big(Re^{\pm 4my_1},y_1e^{\pm4my_1}, R,y_1,y_2,y_3\big)=4(y_2+iy_3)\big[\big(R\cosh(4my_1)+y_1\sinh(4my_1)\big)^2-R^2-4R^3\big]$). 
     If now $P(\xi_1,\xi_2,\eta_1,\dots,\eta_4)=\xi_1^{a_1}\xi_2^{a_2}\eta_1^{b_1}\cdots\eta_4^{b_4}$ 
     is one of the monomials appearing in $Q_{p,q,s}$ and $a:=a_1+a_2$, 
     $b:=b_1+\cdots+b_4$ so that $a\leq 2(1+p+q+s)$ and $a+b\leq 2+3p+2(q+s)$, 
     since $Re^{\pm 4my_1}, y_1e^{\pm 4my_1}=O(r^2)$, we get that: 
      \begin{equation*}
       \frac{P\big(Re^{\pm 4my_1},y_1e^{\pm4my_1}, R,y_1,y_2,y_3\big)}{(r^2)^{1+p+q+s}R^{2+2(p+q+s)}}
                                                     =O\bigg(\frac{(r^2)^{a}R^{b}}{(r^2)^{1+p+q+s}R^{2(1+p+q+s)}}\bigg), 
      \end{equation*}
     and this is $O\big(r^{2a-2(1+p+q+s)}R^{b-2(1+p+q+s)}\big)$; 
     since $a\leq 1+p+q+s$ and $r^{-2}=O(R^{-1})$, 
     this is finally $O\big(R^{a+b-3(1+p+q+s)}\big)$, which in turn is $O(R^{-(1+q+s)})$ since 
     $a+b\leq 2+3p+2(q+s)$. 
     Therefore $\frac{\partial^{p+q+s}\psi_{c}}{\partial y_1^p\partial y_2^q \partial y_3^s}=O(R^{-(1+q+s)})$, 
     and this settles the proof of point 1. of the statement. 
   
 ~
 
  \noindent
  \textit{Asymptotics of $\theta_2+i\theta_3$, and comparison with $dd^c_{I_1}\psi_c$ and $dd^c_{I_1^Y}\psi_c$. } 
     We thus now come to point 2. of this statement. 
     We do it for $\ell=0$; it will become clear from this that the subsequent estimates could be dealt with in an analogous way. 
     Our strategy for proving the desired estimate is the following: 
     first we restrict ourselves to $dd^c_{I_1}\psi_c$; 
     next we decompose $dd^c_{I_1}\psi_c-(\theta_2+i\theta_3)$ into its $dy_1\wedge\eta$-component and 
     its $dy_1\wedge\eta$-free component; 
     we then observe that the $dy_1\wedge\eta$-free components of both $dd^c_{I_1}\psi_c$ and $(\theta_2+i\theta_3)$ 
     have the size we want, 
     whereas we need to look at the $dy_1\wedge\eta$-component of \textit{the very difference} 
     $\big[dd^c_{I_1}\psi_c-(\theta_2+i\theta_3)\big]$ 
     to reach the desired estimate. 
     We conclude by collecting together these estimates, 
     and settling the case of the error term $d(I_1^{Y}-I_1)d\psi_c$.
 
     Since $\psi_c$ is $\mathbb{S}^1$-invariant, 
      \begin{align*}
       dd^c_{I_1}\psi_c 
         =& V^{-1}\bigg(\frac{\partial^2\psi_c}{\partial y_1^2}
                       -V^{-1}\frac{\partial V }{\partial y_1}\frac{\partial \psi_c }{\partial y_1}\bigg)dy_1\wedge\eta
           +\bigg(\frac{\partial^2\psi_c}{\partial y_2^2}+\frac{\partial^2\psi_c}{\partial y_3^2}
                  +V^{-1}\frac{\partial V }{\partial y_1}\frac{\partial \psi_c }{\partial y_1}\bigg)dy_2\wedge dy_3  \\
          &+V^{-1}\bigg(\frac{\partial^2\psi_c}{\partial y_1\partial y_2}
                 -V^{-1}\frac{\partial V }{\partial y_2}\frac{\partial \psi_c }{\partial y_1}\bigg)
                   (dy_2\wedge \eta-Vdy_3\wedge dy_1)                                                                  \\
          &+V^{-1}\bigg(\frac{\partial^2\psi_c}{\partial y_1\partial y_3}
                  -V^{-1}\frac{\partial V }{\partial y_3}\frac{\partial \psi_c }{\partial y_1}\bigg)
                   (dy_3\wedge\eta-Vdy_1\wedge dy_2)
      \end{align*}
     and since $(\xi, -I_1V\xi, \zeta, I_1\zeta)$ is the dual frame of $(\eta,dy_1,dy_2,dy_3)$ 
     and $(\theta_2+i\theta_3)$ is $(1,1)$ for $I_1$, 
     \begin{equation}  \label{eqn_theta2+itheta3}
      \begin{aligned} 
        \theta_2+i\theta_3 = &V(\theta_2+i\theta_3)(\xi,I_1\xi) dy_1\wedge\eta 
                              +(\theta_2+i\theta_3)(\zeta,I_1\zeta)dy_2\wedge dy_3 \\
                             &+(\theta_2+i\theta_3)(\xi,I_1\zeta)(Vdy_1\wedge dy_2-dy_3\wedge \eta) \\
                             &+(\theta_2+i\theta_3)(\xi,\zeta)(Vdy_3\wedge dy_1-dy_2\wedge \eta). 
      \end{aligned}
     \end{equation}
    We already know that (on $R\geq K$), 
    $\tfrac{\partial\psi_c}{\partial y_1} = -4(y_2+iy_3)\big( \tfrac{m}{r^4} -\tfrac{1}{R^3} + \tfrac{1}{4r^4R}\big)$, 
    thus (recall that $\tfrac{\partial(2r^2)}{\partial y_1}=V(|z_1|^2-|z_2|^2) $):
      \begin{equation*}
       \frac{\partial^2\psi_c}{\partial y_1^2} 
            = -4(y_2+iy_3)\Big(-\frac{2mV(|z_1|^2-|z_2|^2)}{r^6} +\frac{3y_1}{R^5} - \frac{y_1}{4r^4R^3}
                         -\frac{V(|z_1|^2-|z_2|^2)}{4r^6R}\Big), 
      \end{equation*}
     the main term of which is $\frac{8mV(y_2+iy_3)(|z_1|^2-|z_2|^2)}{r^6}$, in the sense that 
     it is $O(R^{-1})$, whereas the other summands are $O(R^{-2})$. 
     Moreover, from the estimates of Point 1. and the fact that $\tfrac{\partial V}{\partial y_j}=O(R^{-2})$, 
     $j=1,2,3$, we get that:
      \begin{equation*}
       dd^c_{I_1}\psi_c = \frac{8mV(y_2+iy_3)(|z_1|^2-|z_2|^2)}{r^6}dy_1\wedge\eta +O(R^{-2}). 
      \end{equation*}
     when estimated with respect to $\f$. 
 
     Now recall that $\alpha_j=I_jrdr$, and observe that:
      \begin{align*}
       \theta_2+i \theta_3 &= \frac{rdr\wedge\alpha_2-\alpha_3\wedge\alpha_1+irdr\wedge\alpha_3-i\alpha_1\wedge\alpha_2}{r^6} 
                            = \frac{(rdr-i\alpha_1)\wedge(\alpha_2+i\alpha_3)}{r^6}                                          \\
                           &= \frac{(z_1d\overline{z_1}+z_2d\overline{z_2})\wedge (-z_2dz_1+z_1dz_2)}{r^6}                   \\
                           &= \frac{z_1z_2(dz_1\wedge d\overline{z_1}-dz_2\wedge d\overline{z_2})
                                   +z_2^2dz_1\wedge d\overline{z_2}-z_1^2dz_1\wedge d\overline{z_2}}{r^6}
                            = \frac{\vartheta\wedge \phi}{r^6},                   \\
      \end{align*}
     if we set $\vartheta=z_1d\overline{z_1}+z_2d\overline{z_2}$ and $\phi=-z_2dz_1+z_1dz_2$.  
     Direct computations give: 
      \begin{align*}
       &\vartheta(\xi) = -(|z_1|^2-|z_2|^2),         &  &\vartheta(\zeta) = \frac{2z_1z_2}{iR}\cosh(4my_1),  \\
       &\phi(\xi)= - 2iz_1z_2,                       &  &\phi(\zeta) = -\frac{y_1}{2iR}.  
      \end{align*}
     In particular, $\vartheta(\xi)= O(r^2)$, $\vartheta(\zeta)= O(r^2R^{-1})$, 
     $\phi(\xi)=O(R)$ and $\phi(\zeta) = O(1)$. 
     Moreover, since $\vartheta$ (resp. $\phi$) 
     is (0,1) (resp. (1,0)) for $I_1$, 
     $(\theta_2+i\theta_3)(\xi,I_1\xi)=-\tfrac{2i}{r^6}\vartheta(\xi)\phi(\xi)=\tfrac{8mz_1z_2(|z_1|^2-|z_2|^2)}{r^6}$. 
     Therefore, from \ref{eqn_theta2+itheta3} and since $\vartheta$ (resp. $\phi$) 
     has type (0,1) (resp. (1,0)) for $I_1$, using $r^{-2}=O(R^{-1})$ when necessary:  
      \begin{equation*}
       \theta_2+i\theta_3 = \frac{8mVz_1z_2(|z_1|^2-|z_2|^2)}{r^6}dy_1\wedge\eta +O(R^{-2}). 
      \end{equation*}
     with respect to $\f$. 
     Since $y_2+iy_3=-iz_1z_2$, we thus have (up to a multiplicative constant) 
     $\big|dd^c_{I_1}\psi_c-(\theta_2+i\theta_3)\big|_{\f}=O(R^{-2})$.  
  
     We set $\iota_1^Y:=I_1^{Y}-I_1$, and conclude with an estimate on $\big|d(I_1^{Y}-I_1)d\psi_c\big|_{\f}=|d\iota_1^{Y}d\psi_c|_{\f}$, 
     which is controlled by $|\iota_1^{Y}|_{\f}|\nabla^{\f}d\psi_c|_{\f}+|\nabla^{\f}\iota_1^{Y}|_{\f}|d\psi_c|_{\f}$. 
     But $|\iota_1^{Y}|_{\f}$ and $|\nabla^{\f}\iota_1^{Y}|_{\f}$ are $O(r^{-2})$ hence $O(R^{-1})$ 
     (see e.g. the proof of Proposition \ref{prop_gluing}), 
     and $|d\psi_c|_{\f}$ and $|\nabla^{\f}d\psi_c|_{\f}$ are $O(R^{-1})$ as well from Point 1., 
     and as a result $\big|d(I_1^{Y}-I_1)d\psi_c\big|_{\f}=O(R^{-2})$. 
 
     This settles the case $\ell=0$ of the statement. 
     Cases $\ell=1$ and $2$ are done in the same way, 
     noticing in particular that when letting $\nabla^{\f}$ act on the $(\nabla^{\f})^{j}\psi_c$ 
     or the $(\nabla^{\f})^{j}\iota_1^{Y}$, 
     we keep the same order of precision. 
    \cqfd

 \begin{rmk}
  The function $\psi_c$ is not so small with respect to $\e$, at least at positive orders; 
  for instance, the best we seem able to do on its differential is $|d\psi_c|_{\e}=O(rR^{-1})$. 
 \end{rmk}

   \subsubsection{Comparison between $\f$ and $\f^{\flat}$: proof of Lemma \ref{lem_f/fflat}}       
     Before comparing the metrics, and for this the 1-forms $dy_j^{\flat}:=\beth^{*}dy_j$, $j=1,2,3$, 
     and $\eta^{\flat}:=\beth^{*}\eta$ to their natural ("unflat") analogues, 
     we shall compare the $y_j^{\flat}:=\beth^{*}dy_j$ to the $y_j$, $j=1,2, 3$: 
      \begin{lem}   \label{lem_y_jflat}
       We have $y_j^{\flat}-y_j=O(R^{-1})$, $j=1,2,3$. 
       Consequently if $R^{\flat}:=\beth^*R$, then $R^{\flat}-R=O(R^{-1})$. 
      \end{lem}
     \prf\textit{ of Lemma \ref{lem_y_jflat} -- estimates on $(y_2^{\flat}-y_2)$ and $(y_3^{\flat}-y_3)$}. 
     Since $y_2=\tfrac{1}{2i}\big(z_1z_2-\overline{z_1z_2}\big)$, 
     and $\beth(z_1,z_2)=\alpha(z_1,z_2)$ with $\alpha=1+O(r^{-4})$, 
     it is clear that $y_2^{\flat}= \alpha^2y_2=y_2+O(y_2r^{-4})$, 
     that is $y_2^{\flat}-y_2=O(Rr^{-4})$, and this is $O(R^{-1})$ 
     (recall that $R=O(r^2)$). 
  
     Similarly, $y_3=-\tfrac{1}{2}\big(z_1z_2+\overline{z_1z_2}\big)$, 
     thus $y_3^{\flat}-y_3=y_3(\alpha^2-1)$, which is $O(R^{-1})$. 
  
   ~

      \noindent
      \textit{Estimate on $(y_1^{\flat}-y_1)$.}  
      The case of $y_1^{\flat}$ is slightly more subtle, and for this we shall use the very definition of $y_1$. 
      We fix $(z_1,z_2)\in \C^2$. 
      Since $\beth^*z_1 = \alpha z_1$, $\beth^*z_2 = \alpha z_2$, 
      if one sets $u^{\flat}=\beth^*u$ and $v^{\flat}=\beth^*v$, 
      LeBrun's formulas \eqref{eqn_lebformulas} become: 
       \begin{equation}  \label{eqn_lebrunflat}
        \begin{aligned} 
          \alpha^2 |z_1|^2 &= e^{2m [(u^{\flat})^2-(v^{\flat})^2] } (u^{\flat})^2, \\
          \alpha^2 |z_2|^2 &= e^{2m [(v^{\flat})^2-(u^{\flat})^2] } (v^{\flat})^2,
        \end{aligned}
       \end{equation}
     which we rewrite as:
       \begin{equation}  \label{eqn_lebrunflat2}
        \begin{aligned} 
           |z_1|^2 &= e^{2m\alpha^2[(u^{\flat}/\alpha)^2-(v^{\flat}/\alpha)^2] } (u^{\flat}/\alpha)^2, \\
           |z_2|^2 &= e^{2m\alpha^2[(v^{\flat}/\alpha)^2-(u^{\flat}/\alpha)^2] } (v^{\flat}/\alpha)^2.
        \end{aligned}
       \end{equation}
     These are precisely the equations verified by $u_{m\alpha^2}$ and $v_{m\alpha^2}$ 
     instead of $\tfrac{u^{\flat}}{\alpha}$ and $\tfrac{v^{\flat}}{\alpha}$; 
     by uniqueness of the solutions when $|z_1|$ and $|z_2|$ are fixed, 
     $\tfrac{u^{\flat}}{\alpha}=u_{m\alpha^2}$ and $\tfrac{v^{\flat}}{\alpha}=v_{m\alpha^2}$, 
     that is: 
     $u^{\flat}=\alpha u_{m\alpha^2}$ and $v^{\flat}=\alpha v_{m\alpha^2}$, and consequently
     $y_1^{\flat}=\frac{1}{2}[(u^{\flat})^2+(v^{\flat})^2]=\frac{\alpha^2}{2}(u_{m\alpha^2}^2+v_{m\alpha^2}^2)
                 = \alpha^2 y_{1,m\alpha^2}$. 
             
     Now again with $(z_1,z_2)$ fixed, 
     differentiating LeBrun's equations \textit{with respect to the mass parameter}, $\mu$ say, 
     since we also see $m$ as fixed, and rearranging them gives:
      \begin{equation}  \label{eqn_dy1dm}
       \frac{\partial y_{1,\mu}}{\partial \mu}= -\frac{4R_{\mu}y_{1,\mu}}{1+4\mu R_{\mu}};  
      \end{equation}
     in particular $y_{1,\mu}$ is a non-increasing (resp. non-decreasing) function of $\mu$ on $\{|z_1|\leq |z_2|\}$ 
     (resp. on $\{|z_2|\leq |z_1|\}$). 
 
     Since $\alpha\geq1$, we have for instance on $\{|z_1|\leq |z_2|\}$ the estimate:
      \begin{equation*}
       0 \leq y_{1,m}-y_{1,m\alpha^2} =\int_{m}^{m\alpha^2} \frac{4R_{\mu}y_{1,\mu}}{1+4\mu R_{\mu}} d\mu 
         \leq y_{1,m}\int_{m}^{m\alpha^2} \frac{d\mu}{\mu} 
         =2y_{1,m} \log \alpha,  
      \end{equation*}
     and similarly $0\leq y_{1,m\alpha^2}-y_{1,m}\leq -2y_{1,m} \log \alpha$ on $\{|z_2|\leq |z_1|\}$. 
     Since in both cases $\log \alpha=O(r^{-4})=O(R^{-2})$, we have:
      \begin{equation*}
       y_{1,m\alpha^2}-y_{1,m} = O(y_{1,m}R^{-2}) =O(R^{-1}). 
      \end{equation*}
     Therefore $y_1^{\flat}-y_1=\alpha^2 (y_{1,m\alpha^2}-y_1)+(\alpha^2-1)y_1=O(R^{-1})$ as claimed, 
     since $\alpha-1=O(r^{-4})=O(R^{-2})$ and in particular $\alpha\sim1$ near infinity. 
 
     The estimate $R^{\flat}-R=O(R^{-1})$ comes as follows: 
     $(R^{\flat}-R)(R^{\flat}+R)=(R^{\flat})^2-R^2
                                = (y_1^{\flat})^2-y_1^2+(y_2^{\flat})^2-y_2^2+(y_3^{\flat})^2-y_3^2= O(1)$ 
     from the previous estimates, and thus 
     $R^{\flat}-R=O\big(\frac{1}{R^{\flat}+R}\big)$, which in particular is $O(R^{-1})$. \hfill $\blacksquare$
 
 ~
 
     \noindent
     \textit{Estimates on the $dy_j^{\flat}-dy_j$, $j=1,2,3$, and $\eta^{\flat}-\eta$.} 
     We come back to the proof of Lemma \ref{lem_f/fflat} itself, and start with analysing 
     the transition involved by $\beth$ at the level of 1-forms. 
     We adopt by places the following elementary strategy to evaluate the gap between our fundamental 1-forms 
     and their pull-backs by $\beth$: 
     for $\gamma$ one of the $dy_j$ or $\eta$, 
     we write
      \begin{equation*}
       \gamma^{\flat}= \gamma^{\flat}(\xi)\eta + V^{-1} \gamma^{\flat}(-I_1\xi)dy_1 
                      + \gamma^{\flat}(\zeta)dy_2 +\gamma(I_1\zeta)^{\flat}dy_3,   
      \end{equation*}
     and then evaluate the differences $\gamma^{\flat}(\xi)-\gamma(\xi)$ and the subsequent ones. 
     We start with the easy cases of $dy_2$ and $dy_3$; 
     for more concision, we use the complex expression $\gamma=dy_2+idy_3$. 
 
     Keep the notation $\beth(z_1,z_2)=\alpha (z_1,z_2)$; 
     then $\beth^*(dy_2+idy_3)= d\big(\alpha^2(y_2+iy_3)\big)=\alpha^2(dy_2+idy_3)+(y_2+iy_3)d(\alpha^2)$. 
     Since $\alpha\sim1$, we focus on $d(\alpha^2)$, or rather on $d\alpha$. 
     As $\alpha$ is invariant under the usual action of $\mathbb{S}^1$, we already know that $d\alpha(\xi)=0$. 
     Moreover, by definition, setting $a=\frac{\|\Gamma\|(|\zeta_2|^2+|\zeta_3|^2)}{4}$,
      \begin{equation}  \label{eqn_dalpha}
       d\alpha= -2a\frac{r^2d(r^2)}{(\kappa+r^4)^2},
      \end{equation}
     which we keep under this shape since 
     $d(r^2)=\overline{z_1}dz_1+z_1d\overline{z_1}+\overline{z_2}dz_2+z_2d\overline{z_2}$ is easy to evaluate against $I_1\xi$, 
     $\zeta$ and $I_1\zeta$. 
     As a matter of fact, all computations done:
      \begin{equation} \label{eqn_dalpha123}
       \begin{aligned}
       &d\alpha(-I_1\xi)= -4a\frac{|z_1|^4-|z_2|^4}{(\kappa+r^4)^2}, & 
          & d\alpha(\zeta)= \frac{-8ar^2}{(\kappa+r^4)^2}\frac{y_2\cosh(4my_1)}{R},  \\
       &d\alpha(I_1\zeta)= \frac{-8ar^2}{(\kappa+r^4)^2}\frac{y_3\cosh(4my_1)}{R}. &  &
       \end{aligned}
      \end{equation}
     In particular, $d\alpha(-I_1\xi)=O(r^{-4})=O(R^{-2})$, 
     and $d\alpha(\zeta)=O(R^{-1}r^{-4})$ and $d\alpha(I_1\zeta)=O(R^{-1}r^{-4})$, which are $O(R^{-3})$. 
     Since $\alpha\sim1$ and $y_2+iy_3=O(R)$, 
     we end up with $(dy_2^{\flat}+idy_3^{\flat})(-I_1\zeta)=O(R^{-2}) $, 
     $(dy_2^{\flat}+idy_3^{\flat})(\zeta)=1+O(R^{-1})$ 
     and $(dy_2^{\flat}+idy_3^{I_1\flat})(\zeta)=i+O(R^{-1})$. 
     In other words, 
      \begin{equation*}
       \big|(dy_2^{\flat}+idy_3^{\flat})-(dy_2+idy_3)\big|_{\f}= O(R^{-1}). 
      \end{equation*}

     Analogously to what is done above, the estimate on $dy_1^{\flat}-dy_1$ requires little extra care. 
     First, likewise $y_1$, $y_1^{\flat}$ is invariant under the action of $\mathbb{S}^1$, 
     since $\beth$ commutes to this action; 
     therefore $dy_1^{\flat}(\xi)=0$. 
     Next, pulling-back the known expression for $dy_1$ (proof of Proposition 1.9 in \cite{auv1}) by $\beth$ gives:
      \begin{equation*}
       dy_1^{\flat} = \frac{1}{4mR^{\flat}}\big(e^{-4my_1^{\flat}}d(\alpha^2|z_1|^2)-e^{4my_1^{\flat}}d(\alpha^2|z_2|^2)\big)
      \end{equation*}
     When evaluating $dy_1^{\flat}$, we decompose the term $e^{-4my_1^{\flat}}d(\alpha^2|z_1|^2)-e^{4my_1^{\flat}}d(\alpha^2|z_2|^2)$ 
     into $\sigma:=\alpha^2\big(e^{-4my_1^{\flat}}d(|z_1|^2)-e^{4my_1^{\flat}}d(|z_2|^2)$ 
     and $\rho:=\big(e^{-4my_1^{\flat}}|z_1|^2-e^{4my_1^{\flat}}|z_2|^2\big)d(\alpha^2)
          =\alpha^{-2}\big((u^{\flat})^2-(v^{\flat})^2\big)d(\alpha^2)=4\alpha^{-1}y_1^{\flat}d\alpha$. 
 
     Now $\sigma(-I_1\xi)=2\alpha^2\big(|z_1|^2e^{-4my_1^{\flat}}+|z_2|^2e^{4my_1^{\flat}}\big)=4R^{\flat}$, 
     and by \eqref{eqn_dalpha123}, 
     $\rho(-I_1\xi)=4\alpha^{-1}y_1^{\flat}d\alpha(-I_1\xi)=-16a\alpha^{-1}y_1^{\flat}\frac{|z_1|^4-|z_2|^4}{(\kappa+r^4)^2}$; 
     this way 
      \begin{equation}   \label{eqn_dy1flat1}
       dy_1^{\flat}(-I_1\xi)= (V^{\flat})^{-1}- 8a\alpha^{-1}\frac{y_1^{\flat}}{1+4mR^{\flat}}\frac{|z_1|^4-|z_2|^4}{(\kappa+r^4)^2},  
      \end{equation}
     where $V^{\flat}=\beth^*V=\frac{1+4mR^{\flat}}{2R^{\flat}}$. 
     Since the last summand is $O(r^{-4})$ and thus $O(R^{-2})$, 
     and 
     $(V^{\flat})^{-1}-V^{-1}=\frac{2R^{\flat}}{1+4mR^{\flat}}-\frac{2R}{1+4mR}=2\frac{R^{\flat}-R}{(1+4mR^{\flat})(1+4mR)}=O(R^{-3})$, 
     we have $dy_1^{\flat}(-I_1\xi)=V+O(R^{-2})$. 
 
     Moreover 
     $\sigma(\zeta)=\frac{\alpha^2}{2iR}
                     \big(e^{4m(y_1-y_1^{\flat})}(z_1z_2-\overline{z_1z_2})-e^{-4m(y_1-y_1^{\flat})}(z_1z_2-\overline{z_1z_2})\big)
                   = \alpha^2\frac{y_2}{R}\sinh[4m(y_1-y_1^{\flat})]$,  
     and $\rho(\zeta)=4\alpha^{-1}y_1^{\flat}d\alpha(\zeta)
                     = -32a\alpha^{-1}\frac{r^2}{(\kappa+r^4)^2}\frac{y_1^{\flat}y_2\cosh(4my_1)}{R}$ by \eqref{eqn_dalpha123}. 
     Thus 
      \begin{equation}   \label{eqn_dy1flat2}
       dy_1^{\flat}(\zeta)= \alpha^2\frac{y_2}{2R(1+4mR^{\flat})}\sinh[4m(y_1-y_1^{\flat})]
                            -16a\alpha^{-1}\frac{r^2}{(\kappa+r^4)^2}\frac{y_1^{\flat}y_2\cosh(4my_1)}{R(1+4mR^{\flat})}; 
      \end{equation}
     since $y_1-y_1^{\flat}=O(R^{-1})$, the first summand is $O(R^{-2})$, 
     whereas since $\cosh(4my_1)=O(r^2R^{-1})$, the second summand is $O(R^{-1}r^{-4})$, that is $O(R^{-3})$, 
     and as a result $dy_1^{\flat}(\zeta)=O(R^{-2})$. 
     Similarly $dy_1^{\flat}(I_1\zeta)=O(R^{-2})$ (just replace $y_2$ by $y_3$ in the last equality above). 
 
      ~
 
     \noindent
     \textit{Estimate on $\eta^{\flat}$.} 
     We conclude our estimate of $|\f^{\flat}-\f|_{\f}$ by the estimate on $\eta^{\flat}$. 
     We start with a formula for $\eta^{\flat}$; 
     since on $\{z_1\neq0\}$, 
     $\frac{d(\beth^*\overline{z_1})}{\beth^*\overline{z_1}}-\frac{d(\beth^*z_1)}{\beth^*z_1}
        = \frac{d(\alpha^2\overline{z_1})}{\alpha^2\overline{z_1}}  -\frac{d(\alpha^2z_1)}{\alpha^2z_1}
        = \frac{d\overline{z_1}}{\overline{z_1}}+\frac{d(\alpha^2)}{\alpha^2}
           -\frac{dz_1}{z_1}- \frac{d(\alpha^2)}{\alpha^2}
        = \frac{d\overline{z_1}}{\overline{z_1}}-\frac{dz_1}{z_1}$, 
     and similarly $\beth^*\big(\frac{d\overline{z_2}}{\overline{z_2}}-\frac{dz_2}{z_2}\big)
                    =\frac{d\overline{z_2}}{\overline{z_2}}-\frac{dz_2}{z_2}$ on $\{z_2\neq 0\}$, 
     we have on $\{z_1z_2\neq0\}$, 
     according to the identity 
     $\eta^{\flat} = \frac{i}{4R} \big[u^2\big(\frac{d\overline{z_1}}{\overline{z_1}}-\frac{dz_1}{z_1}\big)
                                                -v^2\big(\frac{d\overline{z_2}}{\overline{z_2}}-\frac{dz_2}{z_2}\big) \big]$ 
     (\cite[Lemma 1.6]{auv1}): 
      \begin{equation}  \label{eqn_etaflat}
       \eta^{\flat} = \frac{i}{4R^{\flat}} \Big[(u^{\flat})^2\Big(\frac{d\overline{z_1}}{\overline{z_1}}-\frac{dz_1}{z_1}\Big)
                                                -(v^{\flat})^2\Big(\frac{d\overline{z_2}}{\overline{z_2}}-\frac{dz_2}{z_2}\Big) \Big]
      \end{equation}
     From this we compute $\eta^{\flat}(\xi)=1$ and $\eta^{\flat}(-I_1\xi)=0$. 
     We also compute $\eta^{\flat}(\zeta)$ as follows:
      \begin{align*}
       \eta^{\flat}(\zeta) 
          =& \frac{i}{4R^{\flat}}\frac{1}{2iR}
             \Big[(u^{\flat})^2e^{4my_1}\Big(\frac{z_2}{\overline{z_1}}-\frac{\overline{z_2}}{z_1}\Big)
                  -(v^{\flat})^2e^{-4my_1}\Big(\frac{z_1}{\overline{z_2}}-\frac{\overline{z_1}}{z_2}\Big) \Big] \\
          =& \frac{i}{4R^{\flat}}\frac{\alpha^2}{2iR}
             \Big[(u^{\flat})^2e^{4my_1}\frac{z_1z_2-\overline{z_1}\overline{z_2}}{\alpha^2|z_1|^2}
                  -(v^{\flat})^2e^{-4my_1}\frac{z_1z_2-\overline{z_1}\overline{z_2}}{\alpha^2|z_2|^2}\Big]              \\
          =& \frac{i\alpha^2 y_2}{2R^{\flat}R} \sinh\big[4m(y_1-y_1^{\flat})\big], 
      \end{align*}
     since from the pulled-back LeBrun's equations \eqref{eqn_lebrunflat}, 
     $\frac{(u^{\flat})^2}{\alpha^2|z_1|^2}=e^{-4my_1^{\flat}}$ and 
     $\frac{(v^{\flat})^2}{\alpha^2|z_2|^2}=e^{4my_1^{\flat}}$. 
     Similarly $\eta^{\flat}(I_1\zeta)=\frac{i\alpha^2 y_3}{2R^{\flat}R} \sinh[4m(y_1-y_1^{\flat})]$, 
     and since $(y_1-y_1^{\flat})=O(R^{-1})$, 
     both $\eta^{\flat}(\zeta)$ and $\eta^{\flat}(I_1\zeta)$ are $O(R^{-2})$. 
     Collecting together those estimates, we get that 
      \begin{equation*}
       |\eta^{\flat}-\eta|_{\f} = O(R^{-2}), 
      \end{equation*}
     which is better than needed. 

     Recall that $\f=V(dy_1^2+dy_2^2+dy_3^2)+V^{-1}\eta$; 
     since $V^{-1}-(V^{\flat})^{-1}$, and similarly $V-V^{-1}$ are $O(R^{-3})$, 
     in view on the estimates we have just proved on the $dy_j-dy_j^{\flat}$ and $\eta^{\flat}-\eta$, 
     we have:
    \begin{equation*}
       |\f^{\flat}-\f|_{\f} = O(R^{-1}).
    \end{equation*}

  ~
 
 \noindent
 \textit{Estimate on $\nabla^{\f}(\f-\f^{\flat})$.} 
 We now aim to prove that $\big|\nabla^{\f}(\f-\f^{\flat})\big|_{\f}=O(R^{-1})$, 
 which is the same as proving that $\big|\nabla^{\f}\f^{\flat}\big|_{\flat}=O(R^{-1})$. 
 In view of the previous estimates on $V-V^{\flat}$,  $V^{-1}-(V^{\flat})^{-1}$, 
 on the $dy_j-dy_j^{\flat}$ and on $\eta-\eta^{\flat}$, and since the $\nabla^{\f}dy_j$ and $\nabla^{\f}\eta$ are $O(R^{-2})$ for $\f$, 
 it will be sufficient for our purpose to see that 
 the $\nabla^{\f}(dy_j-dy_j^{\flat})$ and $\nabla^{\f}(\eta-\eta^{\flat})$ are $O(R^{-1})$ for $\f$. 

 We start with $\nabla^{\f}(dy_2-dy_2^{\flat})$ and $\nabla^{\f}(dy_3-dy_3^{\flat})$. 
 We have $d(y_2+iy_3)-d(y_2^{\flat}+iy_3^{\flat})=(\alpha^2-1)d(y_2+iy_3)+2(y_2+iy_3)\alpha d\alpha$, 
 we know that $\alpha-1=O(r^{-4})=O(R^{-2})$, 
 and we actually proved that $|d\alpha|_{\f}=O(r^{-4})=O(R^{-2})$. 
 Similarly, we will be done if we prove that $|\nabla^{\f}d\alpha|_{\f}$ is still $O(r^{-4})$. 

 Since $\alpha$ is $\S^1$-invariant, 
 $d\alpha=\tfrac{\partial \alpha}{\partial y_1}dy_1+\tfrac{\partial \alpha}{\partial y_2}dy_2+\tfrac{\partial \alpha}{\partial y_3}dy_3$; 
 the $\tfrac{\partial \alpha}{\partial y_j}$ are $\S^1$-invariant as well, and thus
 $\nabla^{\f}d\alpha= \sum_{j,\ell=1}^{3}\tfrac{\partial^2 \alpha}{\partial y_j\partial y_{\ell}}dy_j\otimes dy_{\ell} 
                      + \sum_{j=1}^3\tfrac{\partial \alpha}{\partial y_j}\nabla^{\f}dy_j$. 
 The last summand is $O(R^{-2}r^{-4})$, since the $\tfrac{\partial \alpha}{\partial y_j}$ are $O(r^{-4})$ 
 and the $|\nabla^{\f}dy_j|_{\f}$ are $O(R^{-2})$; 
 we thus focus on the hessian $\sum_{j,\ell=1}^{3}\tfrac{\partial^2 \alpha}{\partial y_j\partial y_{\ell}}dy_j\otimes dy_{\ell}$, 
 and all we need to prove is $\tfrac{\partial^2 \alpha}{\partial y_j\partial y_{\ell}}=O(r^{-4})$ (actually, $O(R^{-2})$) for all $j,\ell$. 
 Now in terms of the $y_j$ variables, 
  \begin{equation*}
   \alpha= 1+\frac{a}{\kappa+\big((y_1^2+y_2^2+y_3^2)^{1/2}\cosh(4my_1)+y_1\sinh(4my_1)\big)^2}, 
  \end{equation*}
 and using that $e^{4m|y_1|}=O(Rr^{-2})$, 
 proving that $\tfrac{\partial^2 \alpha}{\partial y_j\partial y_{\ell}}=O((R\cosh(4my_1)+y_1\sinh(4my_1)\big)^{-2})=O(r{-4})$ 
 for all $j,\ell$ amounts to an easy exercise. 
 This settles the cases of $\nabla^{\f}(dy_2-dy_2^{\flat})$ and $\nabla^{\f}(dy_3-dy_3^{\flat})$. 

 Since our treatment of $dy_j-dy_j^{\flat}$ is a little less conventional, 
 we shall see now how goes that of $\nabla^{\flat}(dy_j-dy_j^{\flat})$. 
 According to formulas \eqref{eqn_dy1flat1} and \eqref{eqn_dy1flat2} and the previous estimates on the derivatives of $r^2$, 
 it is enough to see that $dy_1^{\flat}=O(1)$ and $dR^{\flat}=O(1)$, 
 which are known for the previous step, 
 giving in particular $d\sinh[4m(y_1-y_1^{\flat})]=\cosh[4m(y_1-y_1^{\flat})]d(y_1-y_1^{\flat})$, 
 which is $O(R^{-1})$ (actually $O(R^{-2})$) for $\f$ since $\cosh[4m(y_1-y_1^{\flat})]\sim 1$ and $|d(y_1-y_1^{\flat})|_{\f}=O(R^{-2})$. 

 The treatment of $\eta^{\flat}$ is similar. 

 We prove moreover that $\big|(\nabla^{\f})^2(\f-\f^{\flat})\big|_{\f}=O(R^{-1})$ with the same techniques. 
 \cqfd 

~

 \section{Asymptotics of ALE hyperkähler metrics}   \label{part_ALE}
  We prove in this part an explicit version of Theorem \ref{thm_ALEintro};  
  we indeed compute explicitly the first non-vanishing terms of the hyperkähler data of the ALE gravitational instantons seen as deformation of Kleinian singularities.   
  This gives in particular the asymptotics stated in the previous part, Lemma \ref{lem_omega1Y}, 
  which are crucial in our construction of ALF metrics, as mentioned already. 

  \subsection{Kronheimer's ALE instantons}  \label{section_thmALE}
   \subsubsection{Basic facts and notations}  
   We introduce a few notions about the ALE gravitational instantons constructed by Kronheimer in \cite{kro1} 
   -- and which is exhaustive in the sense that any ALE gravitational instanton is isomorphic to one of Kronheimer's list --, 
   so as to state properly the main result of this part, i.e. Theorem \ref{thm_ALE} of next paragraph, 
   dealing with precise asymptotics of those asymptotically euclidean spaces.

   ~
   
   \noindent 
   \textit{Finite subgroup of $\SU(2)$, and McKay correspondence.} 
    The classification of the finite subgroup of $\SU(2)$ is well-known:
    up to conjugation, in addition to the binary dihedral groups $\mathcal{D}_k$ used in Part 1, 
    one has the cyclic groups of order $k\geq 2$, 
    generated by $\Big(\begin{smallmatrix} e^{2i\pi/k} & 0 \\ 0 & e^{-2i\pi/k} \end{smallmatrix}\Big)$, 
    on the one hand, 
    and the binary tetrahedral, octahedral and icosahedral groups of respective orders 24, 48 and 120, 
    which admit more complicated generators 
    -- all we need to notice for further purpose is that they respectively contain $\mathcal{D}_2$, 
    $\mathcal{D}_3$ and $\mathcal{D}_5$ (among others) as subgroups.  
    When no specification is needed, 
    we shall adopt the notation $\Gamma$ for any fixed group among these finite subgroups of $\SU(2)$. 
  
   ~
   
   \noindent 
   \textit{ALE instantons modelled on $\R^4/\Gamma$.}  
    Kronheimer's construction now consists in producing asymptotically euclidean hyperkähler metrics 
    on \textit{smooth deformations} of the Kleinian singularity $\C^2/\Gamma$, 
    which are diffeomorphic to the \textit{minimal resolution of} $\C^2/\Gamma$. 
    More precisely, the hyperkähler manifolds Kronheimer produces are parametrised as follow: 
    since $\Gamma$ is a finite subgroup of $\SU(2)$, 
    McKay's correspondence \cite{mck} associates a simple Lie algebra, $\mathfrak{g}_{\Gamma}$ say, to this group; 
    for instance, the Lie algebra associated to $\mathcal{D}_k$ is $\mathfrak{so}(2k+4)$ , 
    (also referred to as $D_{k+2}$ -- we prefer the $\mathfrak{so}$ notation which is less confusing when working with binary dihedral groups!). 
    Pick a (real) Cartan subalgebra $\mathfrak{h}$ of $\mathfrak{g}_{\Gamma}$. 
    Then: 
    
    \textit{For any $\zeta\in \mathfrak{h}\otimes\R3$ outside a codimension 3 set $D$, 
    there exists an ALE gravitational instanton $\big(X_{\zeta},g_{\zeta},I_1^{\zeta}, I_2^{\zeta},I_3^{\zeta}\big)$ 
    modelled on $\R^4/\Gamma$ at infinity in the sense that 
    there exists a diffeomorphisms $\Phi_{\zeta}$ between infinities of $X_[\zeta]$ and $\R^4/\Gamma$ such that: 
    ${\Phi_{\zeta}}_*g_{\zeta}-\e=O(r^{-4})$,  ${\Phi_{\zeta}}_*I_j^{\zeta}-I_j=O(r^{-4})$, $j=1,2,3$.}  
    
    The $O$ are here understood in the asymptotically euclidean setting, 
    i.e. $\vareps=O(r^{-a})$ means: for all $\ell\geq0$, $\big|(\nabla^{\e})^{\ell}\vareps\big|=O(r^{-a-\ell})$; 
    since we remain in this setting until the end of this part, 
    we shall keep this convention throughout the following sections \ref{section_hzeta} and \ref{section_ord3term}.

 \subsubsection{Asymptotics of ALE instantons: statement of the theorem}   \label{subsection_thmALE}
    Up to a judicious choice of the ALE diffeomorphism $\Phi_{\zeta}$, 
    which actually is part of Kronheimer's construction, 
    one can be more accurate about the $O(r^{-4})$-error term evoked above. 
    This is the purpose of the main result of this part:
     \begin{thm}  \label{thm_ALE}
      Given $\zeta\in\mathfrak{h}\otimes\R^3-D$, 
      one can choose the diffeomorphism $\Phi_{\zeta}$ between infinities of $X_{\zeta}$ and $\R^4/\Gamma$ such that 
      ${\Phi_{\zeta}}_*g_{\zeta}-\e=h_{\zeta}+O(r^{-6})$, 
      ${\Phi_{\zeta}}_*I_1^{\zeta}-I_1=\iota_1^{\zeta}+O(r^{-6})$ 
      and if $\omega_1^{\zeta}:=g_{\zeta}(I_1^{\zeta}\cdot,\cdot)$, 
      then ${\Phi_{\zeta}}_*\omega_1^{\zeta}-\omega_1^{\e}=\varpi_1^{\zeta}+O(r^{-6})$, 
      where: 
       \begin{equation}  \label{eqn_hzeta}
        \begin{aligned}
         h_{\zeta} = & -\|\Gamma\| \sum_{(j,k,\ell)\in\mathfrak{A}_3}|\zeta_j|^2\frac{(rdr)^2+\alpha_j^2-\alpha_k^2 - \alpha_{\ell}^2}{r^6} 
                       -\|\Gamma\|\langle\zeta_1,\zeta_2\rangle\frac{\alpha_1\cdot\alpha_2 -rdr\cdot\alpha_3}{r^6}                          \\
                     & -\|\Gamma\|\langle\zeta_1,\zeta_3\rangle\frac{\alpha_1\cdot\alpha_3 +rdr\cdot\alpha_2}{r^6}
                       -\|\Gamma\|\langle\zeta_2,\zeta_3\rangle\frac{\alpha_2\cdot\alpha_3-rdr\cdot\alpha_1}{r^6}, 
        \end{aligned}
       \end{equation}
      where $\iota_1^{\zeta}$ is given via the coupling:
       \begin{equation}   \label{eqn_iota1zeta}
        \begin{aligned}
         \e(\iota_1^{\zeta}\cdot,\cdot)=&  \|\Gamma\|(|\zeta_3|^2-|\zeta_2|^2) \frac{\alpha_2\cdot\alpha_3}{r^6} 
                                          - \|\Gamma\|(|\zeta_3|^2+|\zeta_2|^2)\frac{ rdr\cdot\alpha_1}{r^6} \\
                                        &  - \|\Gamma\|\langle\zeta_2,\zeta_3\rangle
                                           \frac{(rdr)^2+\alpha_3^2-\alpha_1^2 - \alpha_{2}^2}{r^6}, 
        \end{aligned}
       \end{equation}
      and 
       \begin{equation}  \label{eqn_varpi1zeta}
        \varpi_1^{\zeta} = -\|\Gamma\||\zeta_1|^2 \theta_1 -\|\Gamma\|\langle\zeta_1,\zeta_2\rangle \theta_2
                           -\|\Gamma\|\langle\zeta_1,\zeta_3\rangle \theta_3, 
       \end{equation}
      with $\|\Gamma\|=\cc|\Gamma|$ for a universal constant $\cc>0$. 
      Moreover, ${\Phi_{\zeta}}_{*}\vol^{g_{\zeta}}=\Omega_{\e}$, 
      and if $\Gamma$ is binary dihedral, tetrahedral, octahedral or icosahedral, 
      the error term can be taken of size $O(r^{-8})$. 
     \end{thm}

 Recall the notations $\alpha_j=I_1rdr$, $j=1,2,3$, and $\theta_a= \frac{rdr\wedge\alpha_a-\alpha_b\wedge\alpha_c}{r^6}$, 
 $(a,b,c)\in \big\{(1,2,3),(2,3,1).(3,1,2)\big\}$. 
 The scalar product on $\mathfrak{h}$ used in this statement is the one induced by the Killing form. 

The rest of this part is devoted to the proof of this result. 
In next section we specify the meaning of the space of parameters $\mathfrak{h}-D$; 
in particular we see how $\mathfrak{h}$ is identified to the degree 2 homology of our Kronheimer's instantons, 
which is helpful in computing the constant $\cc$ of the statement, 
as well as the coefficients appearing in formulas \eqref{eqn_hzeta}-\eqref{eqn_varpi1zeta}. 
We also fix the choice of the diffeomorphisms $\Phi_{\zeta}$, 
and check their properties on volume forms (Lemma \ref{lem_volform}). 
The explicit determination of $h_{\zeta}$, $\iota_1^{\zeta}$ and $\varpi_1^{\zeta}$ is the purpose of section \ref{section_hzeta}, 
and section \ref{section_ord3term}.

 \subsection{Precisions on Kronheimer's construction}   \label{section_kroconst}
  \subsubsection{The degree 2 homology/cohomology}      \label{subsec_deg2hom}
  \noindent\textit{The "forbidden set" $D$.}  
  We keep the notation $\Gamma$ for one of the subgroups of $\SU(2)$ mentioned in the previous section. 
  We saw that Kronheimer's ALE instantons asymptotic to $\R^4/\Gamma$ are parametrised by a triple $\zeta=(\zeta_1,\zeta_2,\zeta_3)\in \mathfrak{h}\otimes\R^3-D$.  
  with $\mathfrak{h}$ is a real Cartan subalgebra of the Lie algebra associated to $\Gamma$ by McKay correspondence; 
  for instance, if $\Gamma=\mathcal{D}_k$, $k\geq2$, then $\mathfrak{h}$ can be taken as the Cartan subalgebra of $\mathfrak{so}(2k+4)$ 
  constituted by matrices of shape $\diag(\lambda_1,\dots,\lambda_{k+2}, -\lambda_1,\dots, -\lambda_{k+2})$. 
  We shall first be more specific about the "forbidden set" $D$; 
  according to \cite[Cor. 2.10]{kro1}, it is the union of codimension 3 subspaces $D_{\theta}\otimes\R^3$ over a positive root system of $\mathfrak{h}$, 
  with $D_[\theta]$ the kernels of the concerned roots; 
  as such, it thus has codimension 3 in $\mathfrak{h}$. 
  
~
 
  \noindent
  \textit{Topology of $X_{\zeta}$.} 
  Recall the notation $\big(X_{\zeta}, g_{\zeta}, I_1^{\zeta}, I_2^{\zeta}, I_3^{\zeta}\big)$ for the hyperkähler manifold of admissible parameter $\zeta$ 
  -- this is actually also defined as a hyperkähler orbifold if $\zeta\in D$. 
  Those spaces are diffeomorphic to the minimal resolution of $\C^2/\Gamma$ (for $I_1$, say) \cite[Cor. 3.12]{kro1}; 
  as such they are simply connected and, again when $\Gamma=\mathcal{D}_k$, their rank 2 topology is given by the diagram:
    
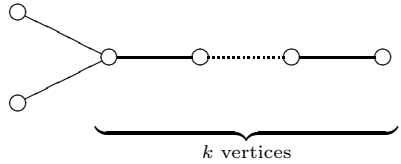
\begin{figure}[!h]
     \centering
      $\xy
       \POS(0,6)  *\cir<3pt>{}="a",
       \POS(12,0) *\cir<3pt>{}="b",
       \POS(24,0) *\cir<3pt>{}="c",
       \POS(36,0) *\cir<3pt>{}="d",
       \POS(48,0) *\cir<3pt>{}="e",
       \POS(0,-6) *\cir<3pt>{}="i",
       \POS(0,-15)*{}="j",
       \POS "a" \ar@{-}^<<<<{} "b",
       \POS "b" \ar@{-}^<<<<{} "c",
       \POS "i" \ar@{-}^<<{} "b",
       \POS "c" \ar@{.}^<<{} "d",
       \POS "d" \ar@{-}^<<{} "e" 
       \txt{\\ \\ \\ \\ \\ \hspace{25pt}  \footnotesize $\underbrace{\, \qquad \qquad \qquad \qquad \qquad  \quad \,}_{k \text{ vertices}}$ 
              \\ \\ } 
       \endxy$ 
    \end{figure}
 \\
(which is nothing but the Dynkin diagram associated to $\mathfrak{so}(2k+4)$), 
where each vertex represents the class of a sphere of $-2$ self-intersection, 
and where two vertices are linked by an edge if and only if the corresponding spheres intersect, 
in which case they intersect normally at one point.  

Furthermore, there is an identification between $H^2(X_{\zeta}, \R)$ and $\mathfrak{h}$ such that:
 \begin{itemize}
  \vspace{-3pt}
  \item the cohomology class of the Kähler form $\omega_j^{\zeta}:= g_{\zeta}\big(I_j^{\zeta}\cdot, \cdot\big)$ is $\zeta_j$, $j=1,2,3$; 
  \item $H_2(X_{\zeta}, \Z)$ is identified with the root lattice of $\mathfrak{h}$; 
        more precisely, given simple roots of $\mathfrak{h}$ and the corresponding basis of $H_2(X_{\zeta}, \Z)$, 
        the intersection matrix of this basis is exactly the opposite of the Cartan matrix of the simple roots, see \cite[p.678]{kro1}; 
        in the case $\Gamma=\mathcal{D}_k$, $k\geq2$, this matrix is thus:
    \begin{equation*}
     ^{k+2} 
     \left[ \begin{aligned} \\ \\ \\ \\ \\ \\ \\  \end{aligned} \right.
     \underbracket[0.75pt]{
     \begin{pmatrix}
         2      &  0       &  -1    &  0       &  \cdots  & 0       \\
         0      &  2       &  -1    &  0       &          & \vdots  \\
         -1     &  -1      &  2     &  -1      &  \ddots  & \vdots  \\
         0      &   0      &  -1    &  \ddots  &  \ddots  & 0       \\
        \vdots  &          & \ddots &  \ddots  &  \ddots  & -1      \\
         0      &  \cdots  & \cdots &     0    &    -1    & 2       \\
     \end{pmatrix} 
     }_{k+2}
     .
  \end{equation*}   
 \end{itemize}
From the latter, we deduce the following lemma, 
identifying cup-product on $H^2$ and the scalar product on $\mathfrak{h}$ induced by the Killing form, up to signs:
 \begin{lem}  \label{lem_cupprod}
  Consider $\alpha, \beta\in H^2(X_{\zeta}, \R)$ with compact support. 
  Then $\alpha\cup\beta=\int_{X_{\zeta}}\alpha\wedge\beta=-\langle\alpha,\beta\rangle$, 
  where the latter is computed with seeing $\alpha$ and $\beta$ in $\mathfrak{h}$ via the above identification.  
 \end{lem}
 \prf. We do it for $\Gamma=\mathcal{D}_k$, $k\geq2$. 
  By Poincaré duality, the computation of $\alpha, \beta\in H^2(X_{\zeta}, \R)$ amounts to that of intersection numbers for a basis of $H_2(X_{\zeta}, \Z)$. 
  But through the above identification between $H_2(X_{\zeta}, \Z)$ and the root lattice of $\mathfrak{h}$ above, 
  the matrix of intersection numbers on the one hand and that of scalar products of the corresponding basis (or dually, of the simple roots) are the same up to signs. 
  \cqfd

~

 \noindent 
 \textit{Period matrix.}  
  For $\zeta\in \mathcal{h}-D$, consider as above a basis $\sigma_j$, $j=1,\dots,r$ say, of $H_2(X_{\zeta}, \Z)$; 
  from the previous paragraph, the \textit{period matrix}
   \begin{equation*}
    P(\zeta)=(P_{j\ell}(\zeta))_{\stackrel{1\leq j\leq 3}{1\leq \ell\leq r}}:= 
      \bigg( \int_{\Sigma_{\ell}}\omega_j^{\zeta} \bigg)_{\stackrel{1\leq j\leq 3}{1\leq \ell\leq r}}
   \end{equation*}
  can be computed thanks to the identities $[\omega_j^{\zeta}]=\zeta_j$. 
  One easily sees that these $P(\zeta)=P(\xi)$ if and only if $\zeta=\xi$. 
  With this formalism Kronheimer's classification \cite[Thm. 1.3]{kro2} can be stated as : 
   \textit{two ALE gravitational instantons are isomorphic as hyperkähler manifolds if and only if they have the same period matrix}. 
  From this we deduce (see also \cite[p.8, (4)]{br}): 
  \begin{lem} \label{lem_zetaAzeta}
   Let $\zeta\in\mathfrak{h}-D$, and let $A\in SO(3)$ act on $\zeta$ and the complex structures $I_j^{\zeta}$ as in section \ref{section_prgrm}. 
   Then there exists a tri-holomorphic isometry between $\big(X_{\zeta},g_{\zeta},(AI^{\zeta})_1,(AI^{\zeta})_2,(AI^{\zeta})_3\big)$ 
   and $\big(X_{A\zeta},g_{A\zeta},I^{A\zeta}_1,I^{A\zeta}_2,I^{A\zeta}_3\big)$. 
  \end{lem}
 \prf. Just check that in both cases, the period matrix is $AP(\zeta)$, and apply Kronheimer's classification theorem. 
  \cqfd

 \subsubsection{Analytic expansions.}   \label{subsec_expnsns}
  \noindent
  \textit{Choice of the map at infinity.} 
  Consider a parameter $\zeta = (\zeta_1, \zeta_2, \zeta_3) \in \mathfrak{h}\otimes\R^3$
and set 
 \begin{equation*}
  \zeta'=(0, \zeta_2,\zeta_3), \qquad \zeta''=(0, 0,\zeta_3)
 \end{equation*}
-- we will keep these notations below. 
As described in \cite[p.677]{kro1}, there exist proper continuous maps: 
 \begin{equation*}
  \begin{aligned}
   & \lambda^{\zeta}_1 : \big(X_{\zeta},g_{\zeta},I_1^{\zeta}, I_2^{\zeta}, I_3^{\zeta}\big)           
       &\longrightarrow & \quad \big(X_{\zeta'},g_{\zeta'},I_1^{\zeta'}, I_2^{\zeta'}, I_3^{\zeta'}\big),           \\ 
   & \lambda^{\zeta'}_2 : \big(X_{\zeta'},g_{\zeta'},I_1^{\zeta'}, I_2^{\zeta'}, I_3^{\zeta'}\big)      
       &\longrightarrow & \quad \big(X_{\zeta''},g_{\zeta''},I_1^{\zeta''}, I_2^{\zeta''}, I_3^{\zeta''}\big),      \\
   & \lambda^{\zeta''}_3 : \big(X_{\zeta''},g_{\zeta''},I_1^{\zeta''}, I_2^{\zeta''}, I_3^{\zeta''}\big) 
       &\longrightarrow & \quad (\R^4/\Gamma, \e, I_1, I_2, I_3).  
  \end{aligned}
  ~
 \end{equation*} 
which are diffeomorphisms (at least) on $(\lambda^{\zeta''}_3\circ \lambda^{\zeta'}_2\circ \lambda^{\zeta}_1)^{-1}(\{0\})$, 
$(\lambda^{\zeta''}_3\circ \lambda^{\zeta'}_2)^{-1}(\{0\})$, and $(\lambda^{\zeta''}_3)^{-1}(\{0\})$ respectively. 
As soon as $\zeta''\notin D$ (resp. $\zeta', \zeta \notin D$), 
$\lambda^{\zeta''}_3$ (resp. $\lambda^{\zeta'}_2, \lambda^{\zeta}_1$) is a resolution of singularities for the third (resp. the second, the first) 
pair of complex structures; 
in particular, if $\zeta' \notin D$ (resp. if $\zeta \notin D$), 
then $\lambda^{\zeta'}_2$ (resp. $\lambda^{\zeta}_1$) is smooth, and holomorphic for the appropriate pair of complex structures. 

To get a "coordinate map" on $X_{\zeta}$ (or rather, to view objects on $\R^4/\Gamma$), one sets: 
 \begin{equation*} 
  F_{\zeta}:(\lambda^{\zeta''}_3\circ \lambda^{\zeta'}_2\circ \lambda^{\zeta}_1)^{-1}:\, (\R^4\backslash\{0\})/\Gamma \longrightarrow X_{\zeta}
 \end{equation*}
(beware this is not exactly the same order of composition as Kronheimer's "coordinate map", 
but this is not a problem by symmetry). 
  
~

 \noindent
 \textit{"Homogeneity" and consequences.} 
We shall see that the $F_{\zeta}$ are going be the $\Phi_{\zeta}$ of Theorem \ref{thm_ALE}. 
For now, according to Proposition 3.14 in \cite{kro1} and its proof, we have for any $\zeta$ the converging expansion
     \begin{equation*}
      {F_{\zeta}}^*g_{\zeta} = \e + \sum_{j=2}^{\infty} h_{\zeta}^{(j)},
     \end{equation*}
  with $h_{\zeta}^{(j)}$ a homogeneous polynomial of degree $j$ in $\zeta$ with coefficients homogeneous symmetric 2-tensors on $\R^4/\Gamma$ 
  -- more precisely, if $\kappa_s$ is the dilation $x\mapsto sx$ for any positive $s$, $\kappa_s^*h_{\zeta}^{(j)}= s^{-2(j-1)}h_{\zeta}^{(j)}$. 
  We will thus be concerned with determining explicitly the term $h_{\zeta}^{(2)}$, 
  and moreover show that when $\Gamma$ is binary dihedral then $h_{\zeta}^{(3)}=0$. 
  For now, observe that Kronheimer's arguments, consisting in analyticity and homogeneity properties of his construction, 
  can also be used to give the existence of an analogous expansion of other tensors such as the complex structures, 
  and therefore the Kähler forms, or the volume forms as well. 
  We can write for example
   \begin{equation}  \label{eqn_cpxstrexpsn}
    {F_{\zeta}}^*I_1^{\zeta} = I_1 + \sum_{j=1}^{\infty} \iota_{1,j}^{\zeta},
   \end{equation}
  where $\iota_{1,j}^{\zeta}$ is a homogeneous polynomial of degree $j$ in $\zeta$ with coefficients (1,1)-tensors, 
  satisfying $\kappa_s^*\iota_{1,j}^{\zeta}= s^{-2j}\iota_{1,j}^{\zeta}$
  (and again, the lower-order term $\iota_{1,1}^{\zeta}$ vanishes, but we will find this fact again below).   
  
 ~

  \subsubsection{Minimal resolutions, invariance of the holomorphic symplectic structure.} 
  We know that as soon as $\zeta\notin D$, $\lambda_1^{\zeta}:\big(X_{\zeta},I_1^{\zeta}\big)\to \big(X_{\zeta'},I_1^{\zeta'}\big)$ 
  is a minimal resolution, 
  and a similar statement holds for $\lambda_2^{\zeta'}:\big(X_{\zeta'},I_2^{\zeta'}\big)\to \big(X_{\zeta''},I_2^{\zeta''}\big)$ 
  and $\lambda_3^{\zeta''}:\big(X_{\zeta''},I_3^{\zeta''}\big)\to (\R^4/\Gamma, I_3)$ whenever $\zeta'\notin D$ or $\zeta''\notin D$, respectively 
  (\cite[p.675]{kro1}). 
  
  As seen already, those maps can happen to be smooth -- for instance $\lambda_1^{\zeta}$ is, when $\zeta, \zeta'\notin D$;  
  we are then only left with their holomorphicity property. 
  This can be used nevertheless with their asymptotic preserving of the hyperkähler structure, 
  to see that they do preserve the appropriate holomorphic symplectic structure:

 \begin{lem}  \label{lem_sympforms}
  Fix $\zeta\in\mathfrak{h}\otimes\R^3$, and assume that $\zeta''\notin D$. 
  Then the map $\lambda_3^{\zeta''}$ verifies:
  $(\lambda_3^{\zeta''})^*(\omega_1^{\e}+i\omega_2^{\e})= \omega_1^{\zeta''}+i\omega_2^{\zeta''}$. 
  
  Similarly, if $\zeta',\zeta''\notin D$, then 
  $(\lambda_2^{\zeta'})^*(\omega_3^{\zeta''}+i\omega_1^{\zeta''})=\omega_3^{\zeta'}+i\omega_1^{\zeta'}$;  
  if $\zeta,\zeta'\notin D$, then 
  $(\lambda_1^{\zeta})^*(\omega_2^{\zeta'}+i\omega_3^{\zeta'})=\omega_2^{\zeta}+i\omega_3^{\zeta}$. 
 \end{lem}
 \prf. 
  The assertion on $\lambda_3^{\zeta''}$ is actually classical, 
  and can be settled in the following elementary way. 
  Call $\theta$ the 2-form $(\lambda_3^{\zeta''})_*(\omega_1^{\zeta''}+i\omega_2^{\zeta''})$, 
  well-defined on $(\R^4\backslash\{0\})/\Gamma$, pulled-back to $\R^4\backslash\{0\}$. 
  Since $\lambda_3^{\zeta''}$ is holomorphic for the pair $\big(I_3^{\zeta''}, I_3\big)$ and $\omega_1^{\zeta''}+i\omega_2^{\zeta''}$ is a holomorphic (2,0)-form for $I_3^{\zeta''}$, 
  $\theta$ is a holomorphic (2,0)-form for $I_3$, 
  and can thus be written as $f(\omega_1^{\e}+i\omega_2^{\e})$, where 
  $f$ is thus holomorphic for $I_3$ on $\R^4\backslash\{0\}$. 
  By Hartogs' lemma it can be extended to the whole $\R^4$; 
  however, since $(\lambda_3^{\zeta''})_*\omega_j^{\zeta''}= {F_{\zeta''}}^*\omega_j^{\zeta''}\sim \omega_j^{\e}$ 
  near infinity on $\R^4/\Gamma$, $j=1,2$, 
  which can be seen as a consequence of the power series expansions analogous to \eqref{eqn_cpxstrexpsn} for Kähler forms, 
  we get that $f$ tends to 1 at infinity. 
  It is therefore constant, equal to 1, 
  which exactly means that $(\lambda_3^{\zeta''})^*(\omega_1^{\e}+i\omega_2^{\e})= \omega_1^{\zeta''}+i\omega_2^{\zeta''}$. 

  We deal with the assertion on $\lambda_2^{\zeta'}$ in a somehow similar way. 
  Since $\zeta',\zeta''\notin D$, $\lambda_2^{\zeta'}$ is a global diffeomorphism between the smooth $X_{\zeta'}$ and $X_{\zeta''}$,  
  holomorphic for the pair $\big(I_2^{\zeta'}, I_2^{\zeta''}\big)$; 
  since $\omega_3^{\zeta'}+i\omega_1^{\zeta'}$ trivialises $K_{(X_{\zeta'},I_2^{\zeta'})}$ and is a (2,0)-holomorphic form for $I_2^{\zeta'}$, 
  $(\lambda_2^{\zeta'})^*(\omega_3^{\zeta''}+i\omega_1^{\zeta''})$ can be written as $f(\omega_3^{\zeta'}+i\omega_1^{\zeta'})$ 
  with $f$ a holomorphic function on $\big(X_{\zeta'}, I_2^{\zeta'}\big)$. 
  Again $f$ tends to 1 near the infinity of $X_{\zeta'}$, since there $(\lambda_2^{\zeta'})^*\omega_j^{\zeta''}\sim\omega_j^{\zeta'}$, $j=1,3$. 
  Moreover $\omega_3^{\zeta''}+i\omega_1^{\zeta''}$ never vanishes on $X_{\zeta''}$, and so neither does $f$ on $X_{\zeta''}$. 
  We collect those observations by saying that $\log(|f|^2)$ is a $g_{\zeta'}$-harmonic function on $X_{\zeta'}$ tending to zero at infinity, 
  and thus identically vanishing. 
  Since $f$ is holomorphic, it is not hard seeing that it is therefore constant, thus $f\equiv1$, or in other words: 
  $(\lambda_2^{\zeta'})^*(\omega_3^{\zeta''}+i\omega_1^{\zeta''})=\omega_3^{\zeta'}+i\omega_1^{\zeta'}$. 
  
  The assertion on $\lambda_1^{\zeta}$ is done in the exact same way. 
\cqfd

~

  An easy but fundamental consequence of the construction of $F_{\zeta}$ via the $\lambda_j^{\zeta}$ and the previous lemma is the invariance of the volume form, 
which we state for $\zeta$ corresponding to smooth $X_{\zeta}$ so as to avoid useless technicalities:
 \begin{lem}  \label{lem_volform}
  The volume form ${F_{\zeta}}^* \vol^{g_{\zeta}}$ does not depend on $\zeta\in \mathfrak{h}\otimes\R^3-D$, 
  and is equal to the standard $\Omega_{\e}$. 
 \end{lem}
\prf. 
Notice first that once we know that ${F_{\zeta}}^* \vol^{g_{\zeta}}$ does not depend on $\zeta$, 
the equality ${F_{\zeta}}^* \vol^{g_{\zeta}}=\Omega_{\e}$ is a direct consequence of the 
expansion of ${F_{\zeta}}^* \vol^{g_{\zeta}}$ as a power series of $\zeta$, the constant term of which is $\Omega_{\e}$. 
To prove that ${F_{\zeta}}^* \vol^{g_{\zeta}}$ is independent of $\zeta$, 
we proceed within three steps, considering first $\zeta''$, and then $\zeta'$ and $\zeta$. 
Even if $\zeta\notin D$, $\zeta'$ or $\zeta''$ might lie in $D$; 
we can however assume this is not the case without loss of generality, 
since ${F_{\zeta}}^* \vol^{g_{\zeta}}$ can be written as a power series of $\zeta$.  
Now from the hyperkähler data $\big(X_{\zeta''},g_{\zeta''},I_1^{\zeta''}, I_2^{\zeta''}, I_3^{\zeta''}\big) $, 
we know that $\vol^{g_{\zeta''}}= \tfrac{1}{2}(\omega_1^{\zeta''})^2$. 
Since ${F_{\zeta}}^*(\omega_1^{\zeta''})= \omega_1^{\e}$ (the standard Kähler form on $\C^2$), 
we get that ${F_{\zeta''}}^* \vol^{g_{\zeta''}}=\Omega_{\e}$. 

Consider now $X_{\zeta'}$; 
we know that $\omega_3^{\zeta'}$ is "preserved" by $\lambda^{\zeta'}_2$, and therefore:
 \begin{align*}
  {F_{\zeta'}}^* \vol^{g_{\zeta'}} = \frac{1}{2} {F_{\zeta'}}^*(\omega_3^{\zeta'})^2 
                                   = \frac{1}{2} {F_{\zeta''}}^*(\lambda^{\zeta'}_2)_*(\omega_3^{\zeta'})^2 
                                   = \frac{1}{2} {F_{\zeta''}}^*(\omega_3^{\zeta''})^2 
                                   = {F_{\zeta''}}^*\vol^{g_{\zeta''}}, 
 \end{align*}
the last equality coming from the fact that $\omega_3^{\zeta''}$ is one of the Kähler forms of the hyperkähler structure 
$(g_{\zeta''},I_1^{\zeta''}, I_2^{\zeta''}, I_3^{\zeta''})$. 

To conclude, we notice that $\omega_3^{\zeta}$ is preserved by $\lambda_1^{\zeta}$ i.e. 
$\omega_3^{\zeta}=(\lambda_1^{\zeta})^*\omega_3^{\zeta'}$, and thus
 \begin{align*}
  {F_{\zeta}}^* \vol^{g_{\zeta}} = \frac{1}{2} {F_{\zeta}}^*(\omega_3^{\zeta})^2 
                                 = \frac{1}{2} {F_{\zeta'}}^*(\lambda^{\zeta}_1)_*(\omega_3^{\zeta})^2 
                                 = \frac{1}{2} {F_{\zeta'}}^*(\omega_3^{\zeta'})^2 
                                 = {F_{\zeta'}}^*\vol^{g_{\zeta'}}. 
 \end{align*}
Here we could also have used the forms $\omega_2^{\zeta}$ and $\omega_2^{\zeta'}$; 
to make a long story short, the reason for the volume form invariance is that 
at each step of the composition of the $\lambda_i^{\zeta}$, at least one Kähler form is preserved. 
\cqfd

~

 \subsection{Explicit determination of $h_{\zeta}$}   \label{section_hzeta}
  \subsubsection{Verifying a gauge}
   
We shall now work more precisely on the first possibly non-vanishing term of the expansion of ${F_{t\zeta}}^*g_{t\zeta}$; 
this allows us to \textit{redefine} $h^{(2)}_{\zeta}$ as follows: 
 \begin{df}
  Fix $\zeta\in \mathfrak{h}\otimes\R^3$, and set on $\R^4\backslash \{0\}$: 
   \begin{equation}  \label{df_hzeta}
     h_{\zeta}:= \frac{1}{2}\frac{d^2}{dt^2}\bigg|_{t=0} {F_{t\zeta}}^*g_{t\zeta}, 
   \end{equation}
  which is then $O(r^{-4})$, 
  with $\ell$th derivatives (for $\nabla^{\e}$) $O(r^{-4-{\ell}})$, near both 0 and infinity, 
  and verifies:
   \begin{equation*}
    {F_{\zeta}}^*g_{\zeta} = \e +  h_{\zeta} + \vareps_{\zeta},
   \end{equation*}
  with $(\nabla^{\e})^{\ell}\vareps_{\zeta}=O(r^{-6-{\ell}})$. 
  More precisely, $h_{\zeta}$ is a homogeneous polynomial of degree 2 in $\zeta$, 
  with coefficients symmetric 2-tensors homogeneous of degree 2 in the sense that $\kappa_s^*h_{\zeta}= s^{-2}h_{\zeta}$, 
  where $\kappa_s$ is the dilation $x\mapsto sx$ of $\R^4\backslash\{0\}$ for any $s>0$;  
  as for $\vareps_{\zeta}$, it is a sum of terms of degree at least 3 in $\zeta$.
 \end{df}

 As indicated by the title of this section, 
 given an admissible $\zeta$, 
 we want to analyse $h_{\zeta}$, which is the first (a priori, possibly) non-vanishing term in the expansion of $g_{\zeta}$ 
 (from now on, for the sake of simplicity, we forget about the $F_{\zeta}$ 
 -- we will be more accurate about this abuse of notation whenever needed). 
 There already exists a rather powerful theory of deformations of Kähler-Einstein metrics; 
 see in particular \cite[ch.12]{bes} for an overview on that subject. 
 Nonetheless, because of the diffeomorphisms action in general, 
 much of the theory is configured so as to work once a gauge is fixed, 
 precisely killing the ambiguity coming from the diffeomorphisms. 

 The following proposition asserts that the $h_{\zeta}$ are indeed in some gauge, 
 making us able of further considerations -- just as is done in paragraph \ref{subsec_Bgge}. 
 Let us specify though that in determining explicitly $h_{\zeta}$, 
 we will be more concerned with other specific properties of that tensor, 
 namely with its inductive decomposition into hermitian and skew-hermitian parts with respect to $I_1$, $I_2$ and $I_3$. 
 As we shall see though, the gauge and the decomposition are rather intricate with one another; 
 seeing the verification of the gauge as a guiding thread, we state: 
  \begin{prop}  \label{prop_h_gge}
   Fix $\zeta\in \mathfrak{h}\otimes\R^3$. 
   Then the lower order term $h_{\zeta}$ of the deformation $g_{\zeta}$ of $\e$ on $\R^4\backslash\{0\}$ 
   is in Bianchi gauge with respect to $\e$, and more precisely:
    \begin{equation*}
     \tr^{\e}(h_{\zeta}) = 0 \qquad \text{ and } \qquad \delta^{\e} h_{\zeta} = 0.
    \end{equation*}
   Moreover, the $I_1$-skew-hermitian part of $h_{\zeta}$ is $h_{\zeta'}$, 
   and the $I_2$-skew-hermitian part of $h_{\zeta'}$ is $h_{\zeta''}$, 
   and $h_{\zeta''}$ is $I_3$-hermitian, while the $I_1$-hermitian part of $h_{\zeta}$, 
   the $I_2$-hermitian part of $h_{\zeta'}$ and $h_{\zeta''}$ give rise to closed forms, that is:
    \begin{equation*}
     d\big(h_{\zeta}(I_1\cdot,\cdot)-h_{\zeta}(\cdot,I_1\cdot)\big)
        = d\big(h_{\zeta'}(I_2\cdot,\cdot)-h_{\zeta'}(\cdot,I_2\cdot)\big)
        = d\big(h_{\zeta''}(I_3\cdot,\cdot)\big) =0 \text{ on }\R^4\backslash\{0\}.
    \end{equation*}
 
  \end{prop}
 \begin{rmk}
  We took the liberty of possibly having $\zeta$ in $D$ since these statements are made on $\R^4\backslash\{0\}$. 
  More precisely, even if $X_{\zeta}$ is not smooth, 
  its orbifold singularities lie above $0\in \R^4$ via $F_{\zeta}$, 
  and $h_{\zeta}$ is smooth on the regular part of $X_{\zeta}$, i.e. $(F_{\zeta})^*h_{\zeta}$ is smooth on $\R^4\backslash\{0\}$. 
 \end{rmk}
\prf. 
Let us deal first with the assertion on $\tr^{\e}(h_{\zeta})$. 
At any point of $(\R^4 \backslash \{0\})/\Gamma$, for any $t$:
 \begin{equation*}
  \vol^{g_{t\zeta}} = \det^{\e} (g_{t\zeta}) \Omega_{\e} 
                    = \det^{\e} \big(\e + t^2h_{\zeta} + O(t^{3})\big)\Omega_{\e}
                    = \big(1+ t^2 \tr^{\e}(h_{\zeta}) + O(t^{3})\big)\Omega_{\e}.
 \end{equation*}
But we saw in Lemma \ref{lem_volform} that for all $t$, $\vol^{g_{t\zeta}}=\Omega_{\e}$; 
consequently, $\tr^{\e}(h_{\zeta})=0$. 

We now deal with the divergence assertion. 
As for the previous lemma, 
we proceed inductively on the shape of $\zeta$; 
the hermitian/skew-hermitian decomposition as well as the closedness property will come out along the different steps of the induction.  
For this we assume that $\zeta'=(0,\zeta_2,\zeta_3)$ and $\zeta''=(0,0,\zeta_3)$ are as well out of the "forbidden set" $D$. 
Again, since $h_{\zeta}$ can be written as a sum of quadratic polynomials of $\zeta$ times symmetric 2-forms independent of $\zeta$, 
this assumption does not actually lead to a loss of generality.  

~

\noindent
\textit{Step 1: $\delta^{\e} h_{\zeta''}=0$.} 
We hence start with $\zeta''=(0,0,\zeta_3)$. 
Since $I_3$ is parallel for $\e$, we have that $d^{*_{\e}}[h_{\zeta''}(\cdot,I_3\cdot)]=(\delta^{\e}h_{\zeta''})(I_3\cdot)$; 
indeed, given any local $\e$-orthonormal frame $(e_j)_{j=1,\dots,4}$: 
 \begin{equation}   \label{eqn_dvgces}
  d^{*_{\e}}[h_{\zeta''}(\cdot,I_3\cdot)] = -\sum_{j=1}^4 e_j\lrcorner \big[\nabla_{e_j}^{\e} \big(h_{\zeta''}(\cdot,I_3\cdot)\big)\big] 
   \quad\text{and}\quad 
  \delta^{\e}h_{\zeta''}= -\sum_{j=1}^4 (\nabla_{e_j}^{\e} h_{\zeta''})(e_j, \cdot), 
 \end{equation}
see for instance \cite[1.2.11]{biq} for the first equality, and \cite[1.2.13]{biq} for the second one.  
Moreover $h_{\zeta''}$ is clearly $I_3$-hermitian, since the $g_{t\zeta''}$ are, 
which is straightforward from the holomorphicity of the $\lambda_3^{t\zeta''}$ for the pairs $\big(I_3^{t\zeta''}, I_3\big)$; 
$h_{\zeta''}(\cdot,I_3\cdot)$ is therefore a (1,1)-form for $I_3$.  
It is furthermore closed, since the $g_{t\zeta''}(\cdot,  I_3\cdot)$ are. 
We can now use the Kähler identity with the structure $(\e, I_3)$  to write:
 \begin{equation*}
  d^{*_{\e}}\big(h_{\zeta''}(\cdot,I_3\cdot)\big) 
              = [\Lambda_{\omega_3^{\e}},d^c_{I_3}]\big(h_{\zeta''}(\cdot,I_3\cdot)\big).
 \end{equation*}
But $\Lambda_{\omega_3^{\e}}\big(h_{\zeta''}(\cdot,I_3\cdot)\big)=-\frac{1}{2}\tr^{\e}(h_{\zeta''})=0$, 
and since $h_{\zeta''}(\cdot,I_3\cdot)$ is $I_3$-hermitian and closed, 
$d^c_{I_3}\big(h_{\zeta''}(\cdot,I_3\cdot)\big)=d\big(h_{\zeta''}(\cdot,I_3\cdot)\big)=0$, hence the result. 

~

\noindent
\textit{Step 2: $\delta^{\e} h_{\zeta''}=0$.} 
We go on our induction and analyse $h_{\zeta'}$, 
where we recall the notation $\zeta'=(0,\zeta_2, \zeta_3)$. 
We proceed through the following lines:
 \begin{enumerate}[(i)]
  \vspace{-3pt}
  \item we come back momentarily to $h_{\zeta''}$ and prove it is $I_2$-skew-hermitian;
  \item we prove that the $I_2$-skew-hermitian part of $h_{\zeta'}$ is $h_{\zeta''}$, 
        which is known to be divergence-free for $\e$; 
  \item we conclude by proving that the $I_2$-hermitian part of $h_{\zeta'}$ is $\e$-divergence-free as well. 
 \end{enumerate}
We tackle Point (i). 
Recall that the map $\lambda_2^{\zeta'}:X_{\zeta'}\to X_{\zeta''}$ is holomorphic 
for the pair $\big(I_2^{\zeta'},I_2^{\zeta''}\big)$; 
since we forget about $F_{\zeta'}$ and $F_{\zeta''}$, 
this amounts to writing $I_2^{\zeta'}=I_2^{\zeta''}$. 
Recall that in the same way as for the metric, 
the complex structures admit an analytic expansion, 
which can be written as a power series of $\zeta$ with coefficients homogeneous (1,1)-tensors on $(\R^4\backslash\{0\})/\Gamma$.   
\textit{We assume momentarily that the first order variation vanishes}, 
and we thus write $I_2^{\zeta''}=I_2+\iota_2^{\zeta''}+\epsilon_2^{\zeta''}$, 
where $\iota_2^{\zeta''}=\frac{1}{2}\frac{d^2}{dt^2}\big|_{t=0} I_2^{t\zeta''}$, 
is $O(r^{-4})$ (with according decay on derivatives), 
and $(\nabla^{\e})^{\ell}\epsilon_2^{\zeta''}=O(r^{-6-\ell})$ for all $\ell\geq0$. 

Now $\iota_2^{\zeta''}$ splits into an $\e$-symmetric part and an $\e$-anti-symmetric part. 
But according to \cite[12.96]{bes}, to the anti-symmetric part, $(\iota_2^{\zeta''})^a$ say, 
corresponds an $I_2$-holomorphic $(2,0)$-form $\theta$ via 
the coupling $\e\big(\cdot,(\iota_2^{\zeta'})^a\cdot\big)=\theta$; 
this we get by considering the second order variation of the Kähler-Einstein deformation $\big(g_{t\zeta''},I_{2}^{t\zeta''}\big)$, 
satisfying the gauge $\tr^{\e}(h_{\zeta''})=\delta^{\e}h_{\zeta''}=0$, 
and observing that all the statements are \textit{ local}. 
We can lift $\theta$ on $\R^4\backslash\{0\}$, and then write $\theta=f dw_1\wedge dw_2$, 
where $w_1$ and $w_2$ are the standard $I_2$-holomorphic coordinates $x_1+ix_3$ and $x_4+ix_2$, 
and $f$ is thus $I_2$-holomorphic with decay $r^{-4}$ at infinity. 
By Hartogs' lemma we can extend $f$ through $0$; 
we thus have an entire function on $(\R^4, I_2)$, decaying at infinity: 
the only possibility is $f\equiv 0$, and therefore $(\iota_2^{\zeta''})^a=0$, or: 
$\iota_2^{\zeta''}$ is $\e$-symmetric. 

Here we would like to follow \cite[12.96]{bes} again, 
to see for example that $\iota_2^{\zeta''}$ then corresponds to the $I_2$-skew-hermitian part of $h_{\zeta''}$, 
via the coupling $\omega^{\e}_2(\cdot,\iota_2^{\zeta''}\cdot)$ -- 
this latter (2,0)-tensor being clearly $I_2$-skew-hermitian, 
because $\omega_2^{\e}$ is $I_2$-hermitian, 
and since for all $t$, $-1= (I_2^{t\zeta''})^2= I_2^2 +t^2(I_2\iota_2^{\zeta''} + \iota_2^{\zeta''}I_2) +O(t^3) $, 
thus $I_2\iota_2^{\zeta''}=-\iota_2^{\zeta''}I_2$. 
Since in our situation,  $\omega_2^{\zeta''}$ does not vary, 
we could also expect from  \cite[12.95]{bes} that the $I_2$-hermitian part of $h_{\zeta''}$ vanishes. 
Nonetheless some of the quoted arguments are of global nature, and one should check they can be adapted to our framework. 
This can be bypassed however by a rather simple computation, which we quote here: 
for any $t$, 
 \begin{equation*}
  g_{t\zeta''} = \left\{\begin{aligned} 
                           &\omega_2^{t\zeta''}(\cdot,I_2^{t\zeta''}\cdot) 
                           =\omega_2^{\e}(\cdot,I_2\cdot)+ t^2\omega_2^{\e}(\cdot,\iota_2^{\zeta''}\cdot) +O(t^3)
                            \quad\text{since}\quad\omega_2^{t\zeta''}=\omega_2^{\e} \\
                           &\e + t^2h_{\zeta''} +O(t^3),
                        \end{aligned} \right.
 \end{equation*}
and thus $h_{\zeta''}=\omega_2(\cdot,\iota_2^{\zeta''}\cdot) $ which is $I_2$-skew-hermitian, as announced.  

~

We now claim that the $I_2$-skew-hermitian part of $h_{\zeta'}$ is nothing but $h_{\zeta''}$, 
which is Point (ii) of the current step. 
Indeed, since for all $t$, $I_2^{t\zeta'}=I_2^{t\zeta''}$ (consider $\lambda_2^{t\zeta'}$), 
 \begin{align*}
  0=& g_{t\zeta'}(I_2^{t\zeta'}\cdot,I_2^{t\zeta'}\cdot) - g_{t\zeta'}    
   =  g_{t\zeta'}(I_2^{t\zeta''}\cdot,I_2^{t\zeta''}\cdot) - g_{t\zeta'}       \\
   =& \e(I_2^{t\zeta''}\cdot,I_2^{t\zeta''}\cdot) + t^2h_{\zeta'}(I_2^{t\zeta''}\cdot,I_2^{t\zeta''}\cdot)
     - \e - t^2h_{\zeta'} +O(t^{3})                                                                         \\
   =& \underbrace{\e(I_2\cdot, I_2\cdot)}_{=\e} + t^2 \e(I_2\cdot, \iota_2^{\zeta''}\cdot) + t^2 \e(\iota_2^{\zeta''}\cdot, I_2\cdot)  
      + t^2 h_{\zeta'}(I_2\cdot,I_2\cdot) -\e- t^2h_{\zeta'} +O(t^{3}),
 \end{align*}
and thus $h_{\zeta'}-h_{\zeta'}(I_2\cdot,I_2\cdot)=\e(I_2\cdot, \iota_2^{\zeta''}\cdot) + \e(\iota_2^{\zeta''}\cdot, I_2\cdot) $. 
We know that $\e(I_2\cdot, \iota_2^{\zeta'}\cdot)=\omega_2(\cdot,\iota_2^{\zeta''}\cdot)= h_{\zeta''}$. 
To conclude, use that $\e$ and $\e(\cdot,\iota_2^{\zeta''})$ are both symmetric, 
that $I_2\iota_2^{\zeta''}=-\iota_2^{\zeta''}I_2$, and that $\e$ is $I_2$-hermitian 
to see that for all $X$, $Y$, 
 \begin{equation*}
  \e(\iota_2^{\zeta''}X, I_2Y)= \e(I_2Y, \iota_2^{\zeta'}X)=-\e(Y, I_2\iota_2^{\zeta''}X)=\e(Y, \iota_2^{\zeta''}I_2X)
                              = \e(I_2X, \iota_2^{\zeta''}Y),
 \end{equation*}
i.e. $\e(\iota_2^{\zeta''}\cdot, I_2\cdot)=\e(I_2\cdot, \iota_2^{\zeta''}\cdot)=h_{\zeta''}$. 
We have proved that
 \begin{equation*}
  \frac{1}{2}\big(h_{\zeta'}-h_{\zeta'}(I_2\cdot,I_2\cdot)\big)= h_{\zeta''},
 \end{equation*}
as claimed.  
Since $\delta^{\e}h_{\zeta''}=0$, to see that $\delta^{\e}h_{\zeta'}=0$, 
we are only left with checking this identity on the $I_2$-hermitian part of $h_{\zeta'}$, 
which is Point (iii) of the current induction step.   

~

For this, let us call $\varphi$ this tensor twisted by $I_2$, 
namely $\varphi=\tfrac{1}{2}\big(h_{\zeta'}(I_2\cdot,\cdot)-h_{\zeta'}(\cdot,I_2\cdot)\big)$. 
As above, we want to see that $d^{*_{\e}}\varphi=0$. 
This is clearly an $I_2$-hermitian 2-form, that is an $I_2$-(1,1)-form. 
It is moreover trace-free with respect to $\e$, since $h_{\zeta'}$ is. 
If we check it is closed then we are done, using the Kähler identity $[\Lambda_{\omega_2},d^c_{I_2}]$. 
For this, we use an expansion of $\omega_2^{\zeta'}$: 
for all $t$, 
 \begin{align*}
  \omega_2^{t\zeta'}  
     =& \frac{1}{2}\big(g_{t\zeta'}(I_2^{t\zeta'}\cdot,\cdot)
                      -g_{t\zeta'}(\cdot,I_2^{t\zeta'}\cdot)\big)
     =  \frac{1}{2}\big(g_{t\zeta'}(I_2^{t\zeta''}\cdot,\cdot)
                      -g_{t\zeta'}(\cdot,I_2^{t\zeta''}\cdot)\big) \\
     =& \frac{1}{2}\big(\e(I_2\cdot,\cdot) + t^2\e(\iota_2^{\zeta''}\cdot,\cdot)+ t^2h_{\zeta'}(I_2\cdot,\cdot) \\
      &\quad       \,   
                       -\e(\cdot,I_2\cdot)- t^2\e(\cdot,\iota_2^{\zeta''}\cdot)
                       -t^2h_{\zeta'}(\cdot,I_2\cdot) \big) + O(t^{3})            \\
     =& \omega_2^{\e}+ t^2\varphi + O(t^{3}), \text{ since } \e(I_2\cdot,\cdot)=-\e(\cdot,I_2\cdot)=\omega_2^{\e} 
                                        \text{ and }\e(\cdot,\iota_2^{\zeta''}\cdot)= \e(\iota_2^{\zeta''}\cdot,\cdot);
 \end{align*}
this expansion can be differentiated term by term, so that $ t^2d\varphi + O(t^{3})=0$, 
hence $d\varphi=0$, as wanted. 

~

\noindent
\textit{Step 3: $\delta^{\e} h_{\zeta}=0$.} 
We now analyse $h_{\zeta}$. 
All the techniques to pass from $h_{\zeta''}$ to $h_{\zeta'}$ can actually be used again, 
and bring us to the desired conclusion:
 \begin{enumerate}
  \vspace{-3pt}
  \item
                   we first observe that $I_1^{\zeta}=I_1^{\zeta'}$, 
                   and we define $\iota_1^{\zeta'}=\tfrac{d^2}{dt^2}\big|_{t=0}I_1^{t\zeta'}$ 
                   which we assume again to be the possibly lower-order non-vanishing variation of $I_1^{t\zeta'}$; 
                   then $(\iota_1^{\zeta'})^a=0$, since otherwise we would have a non-trivial entire function on $\C^2$ 
                   going to 0 at infinity;
  \item
                   since $\omega_1^{\zeta'}= \omega_1^{\zeta''}=\omega_1^{\e}$, 
                   we get that $h_{\zeta'}$ is $I_1$-skew-hermitian, given by $\iota_1^{\zeta'}$ via 
                   the identity $h_{\zeta'}=\omega_1^{\e}(\cdot,\iota_1^{\zeta'}\cdot)$, 
                   and that the $I_1$-skew-hermitian component of $h_{\zeta}$ coincides with $h_{\zeta'}$, 
                   the $\delta^{\e}$ of which vanishes;  
                   we are thus left with the $I_1$-hermitian component of $h_{\zeta}$; 
  \item            this component is $\e$-trace-free ($h_{\zeta}$ is), 
                   and gives rises to an $I_1$-hermitian 2-form $\psi$, which is closed since the $\omega_1^{t\zeta}$ are; 
                   the Kähler identity $[\Lambda_{\omega_1^{\e}},d^c_{I_1}]=d^{*_{\e}}$ then leads us to 
                   $d^{*_{\e}}\psi=0$, which is equivalent to: 
                    \begin{equation*}
                     \delta^{\e}\big(I_1\text{-hermitian component of }h_{\zeta}\big)=0.
                    \end{equation*}
 \end{enumerate}
To finish this proof, we justify our assumption of the vanishing of the first order variation of the complex structures. 
For instance, let us not assume that $\iota:=\tfrac{d}{dt}\big|_{t=0}I_2^{t\zeta''}$ is a priori vanishing. 
Then it is defined on $\R^4\backslash\{0\}$, and is $O(r^{-2})$. 
Now as above, since $0=\tfrac{d}{dt}\big|_{t=0}g_{t\zeta''}$ has vanishing trace and divergence for $\e$, 
the $\e$-anti-symmetric part of $\iota$ has to vanish since it gives rise to a holomorphic function on $(\R^4, I_2)$ decaying at infinity. 
And we see as above that $\omega_2^{\e}(\cdot,\iota\cdot)=\tfrac{d}{dt}\big|_{t=0}g_{t\zeta''}=0$, and thus $\iota=0$. 
Similarly, the arguments for $\iota_1^{\zeta'}$ apply to $\tfrac{d}{dt}\big|_{t=0}I_1^{t\zeta'}$ with $0=\tfrac{d}{dt}\big|_{t=0}g_{t\zeta'}=0$ instead of $h_{\zeta'}$, 
so that $\tfrac{d}{dt}\big|_{t=0}I_1^{t\zeta'}=0$. 
\cqfd

 \begin{rmk}
  By contrast with what is usually done, we used properties already known of $h_{\zeta'}$ and $h_{\zeta''}$,  
  conjugated to properties of mappings between $X_{\zeta}$, $X_{\zeta'}$ and $X_{\zeta''}$ to show that indeed, 
  our first order deformations were in gauge, 
  which is also retroactively used in some places, e.g. in killing tensors like $(\iota^{\zeta''}_2)^a$. 
 \end{rmk}

 \subsubsection{Lower order variation of the Kähler forms: general shape}

 As seen when proving that the gauge was verified, 
given $\zeta\in \mathfrak{h}\otimes\R^3$, 
$h_{\zeta''}$ is $I_3$-hermitian, the $I_2$-skew-hermitian part of $h_{\zeta'}$ is $h_{\zeta''}$, 
and the $I_1$-skew-hermitian part of $h_{\zeta'}$ is $h_{\zeta''}$. 
In order to determine $h_{\zeta}$ completely, 
we are thus left with working on the respective $I_3$, $I_2$ and $I_1$-hermitian components of 
$h_{\zeta''}$, $h_{\zeta'}$ and $h_{\zeta}$, or equivalently on 
the respectively $I_3$, $I_2$ and $I_1$-(1,1) forms 
 \begin{align*}
  \varpi_3^{\zeta''}:= h_{\zeta''}&(I_3\cdot,\cdot), \qquad  \qquad
  \varpi_2^{\zeta'} := \frac{1}{2}\big(h_{\zeta'}(I_2\cdot,\cdot)- h_{\zeta'}(\cdot,I_2\cdot)\big)  \\ 
  \text{and}&\quad \varpi_1^{\zeta}  := \frac{1}{2}\big(h_{\zeta}(I_1\cdot,\cdot)-h_{\zeta}(\cdot,I_1\cdot) \big). 
 \end{align*}

We interpret these forms as the first (possibly) non-vanishing variation term of $\omega_3^{\zeta''}$, 
$\omega_2^{\zeta'}$ and $\omega_1^{\zeta}$; as such and as seen above, these are closed forms. 
More precisely, they follow a general common pattern:
 \begin{prop}  \label{prop_varpishape}
  There exist real numbers $a_{1j}(\zeta)$, $a_{2j}(\zeta')$, $a_{3j}(\zeta'')$, $j=1,2,3$, such that:
   \begin{align*}
    \varpi_3^{\zeta''}  = a_{31}(\zeta'')\theta_1 + a_{32}(\zeta'')&\theta_2 + a_{33}(\zeta'')\theta_3,  \quad 
    \varpi_2^{\zeta'}   = a_{21}(\zeta')\theta_1  + a_{22}(\zeta')\theta_2   + a_{23}(\zeta')\theta_3 , \\
     \text{and} \quad 
    \varpi_1^{\zeta}    &= a_{11}(\zeta)\theta_1  + a_{12}(\zeta)\theta_2   + a_{13}(\zeta)\theta_3.
   \end{align*}
 where we recall the notations:
   \begin{align*}
    \theta_1 = \frac{rdr\wedge \alpha_1-\alpha_2\wedge \alpha_3}{r^6}, \quad \theta_2 = \frac{rdr\wedge \alpha_2-\alpha_3\wedge \alpha_1}{r^6}, \quad
    \theta_3 = \frac{rdr\wedge \alpha_3-\alpha_1\wedge \alpha_2}{r^6}.
   \end{align*}
 \end{prop}
\prf. 
We do it for $\varpi_1^{\zeta}$, as it will be clear that the arguments would apply similarly to $\varpi_2^{\zeta'}$ and $\varpi_3^{\zeta''}$; 
we work on $\R^4\backslash\{0\}$. 
As $\varpi_1^{\zeta}$ is of type (1,1) for $I_1$, 
it is at any point a linear combination of $rdr\wedge \alpha_1$, 
$\alpha_2\wedge \alpha_3$, $rdr\wedge \alpha_2-\alpha_3\wedge \alpha_1$ and $rdr\wedge \alpha_3-\alpha_1\wedge \alpha_2$. 

The symmetric tensor $\frac{1}{2}\big(h_{\zeta}+h_{\zeta}(I_1\cdot,I_1\cdot)\big)$ corresponding to $\varpi_1^{\zeta}$ is moreover trace-free for $\e$, 
which translates into $\varpi_1^{\zeta}\wedge \omega_1^{\e}= 0$. 
Since $\omega_1= \frac{rdr\wedge \alpha_1+\alpha_2\wedge \alpha_3}{r^2}$, 
we have $(rdr\wedge \alpha_2-\alpha_3\wedge \alpha_1)\wedge\omega_1^{\e}= (rdr\wedge \alpha_3-\alpha_1\wedge \alpha_2)\wedge\omega_1^{\e}= 0$, 
whereas $rdr\wedge \alpha_1\wedge\omega_1^{\e}=\alpha_2\wedge \alpha_3\wedge\omega_1^{\e}$. 
As a consequence, the pointwise coefficient of $\alpha_2\wedge \alpha_3$ is the opposite of that of $\alpha_2\wedge\alpha_3$. 
To sum up, since the $\theta_j$ are $O(r^{-4})$ with corresponding decay (or growth, near 0) of their derivatives, 
which are precisely the orders of $\varpi_1^{\zeta}$, we know that:
 \begin{equation*}
  \varpi_1^{\zeta} = f\theta_1 + g\theta_2 + h\theta_3,
 \end{equation*}
for three bounded functions $f$, $g$, $h$, with euclidean $\ell$-th order derivatives of order $O(r^{-\ell})$, near 0 and infinity. 
We can be more precise here: 
from the properties of analytic expansions in play discussed in paragraph \ref{subsec_expnsns}, 
we have that $\kappa_s^*\varpi_1^{\zeta} = s^{-2}\varpi_1^{\zeta}$, where $\kappa_s$ is the dilation of factor $s>0$ on $\R^4$. 
But we exactly have $\kappa_s^*\theta_j = s^{-2}\theta_j$, $j=1,2,3$; 
therefore, $f$, $g$, $h$ are functions on the sphere $\mathbb{S}^3$. 
Notice that from this point, we also know that $\varpi_1^{\zeta}$ is anti-self-dual (for $\e$), 
since the $\theta_j$'s are. 
Therefore $\varpi_1^{\zeta}$ is $\e$-harmonic on $\R^4\backslash\{0\}$, 
which is the same as: $(\nabla^{\e})^*(\nabla^{\e})\varpi_1^{\zeta}= 0$. 
On the other hand, the $\theta_j$ are harmonic as well: 
they are anti-self-dual, and closed, since  
 \begin{equation} \label{eqn_theta_j}
  \theta_j=\frac{1}{4}dd_{I_j}^c \Big(\frac{1}{r^2}\Big), \qquad j=1,2,3.
 \end{equation}
Putting those facts together and setting $e_j=I_j \tfrac{1}{r}x_i\tfrac{\partial}{\partial x_i}$, $j=1,2,3$ 
-- forget about formulas \eqref{eqn_ei} -- 
so that $rdr(e_j)=0$ and $\alpha_k(e_j)=r\delta_{jk}$, $j,k=1,2,3$, 
we get that: 
 \begin{equation*}
  \Delta_{\e}(f\theta_1) = \frac{1}{r^2}(\Delta_{\mathbb{S}^3}f)\theta_1 - 2\sum_{k=1}^3 (e_k\cdot f) \nabla_{e_k}^{\e}\theta_1
 \end{equation*}
The $\nabla_{e_k}^{\e}\theta_1$ are easy to compute: 
since $e_k\cdot r=0$, $\nabla_{e_k}^{\e}\theta_1=\tfrac{1}{r^6}\nabla_{e_k}^{\e}(rdr\wedge \alpha_1-\alpha_2\wedge \alpha_3)$. 
Moreover since the $I_j$ are parallel, we just have to compute $\nabla_{e_k}^{\e}rdr$; 
since $\nabla^{\e}(rdr)= \e$, $\nabla_{e_k}^{\e}(rdr)= \e(e_k,\cdot)= \tfrac{1}{r}\alpha_k$. 
Therefore $\nabla_{e_1}^{\e}\theta_1=0$, $\nabla_{e_2}^{\e}\theta_1=\tfrac{2}{r}\theta_3$ and $\nabla_{e_3}^{\e}\theta_1=-\tfrac{2}{r}\theta_2$. 
Thus $\Delta(f\theta_1)= \frac{1}{r^2}(\Delta_{\mathbb{S}^3}f)\theta_1-\tfrac{2}{r}(e_2\cdot f)\theta_3+\tfrac{2}{r}(e_3\cdot f)\theta_2$. 
A circular permutation on the indices gives as well 
$\Delta(g\theta_2)= \frac{1}{r^2}(\Delta_{\mathbb{S}^3}g)\theta_2-\tfrac{2}{r}(e_3\cdot g)\theta_1+\tfrac{2}{r}(e_1\cdot g)\theta_3$ and 
$\Delta(h\theta_3)= \frac{1}{r^2}(\Delta_{\mathbb{S}^3}h)\theta_3-\tfrac{2}{r}(e_1\cdot h)\theta_2+\tfrac{2}{r}(e_2\cdot h)\theta_1$. 
Since the $\theta_j$ are linearly independent, $\Delta \varpi_1^{\zeta}=0$ translates into:
 \begin{equation} \label{eqn_laplfgh}
  \Delta_{\mathbb{S}^3}f - 4(e_3\cdot g - e_2\cdot h) = 
  \Delta_{\mathbb{S}^3}g - 4(e_1\cdot h - e_3\cdot f) = 
  \Delta_{\mathbb{S}^3}h - 4(e_2\cdot f - e_1\cdot g) = 0. 
 \end{equation}
On the other hand, $d\varpi_1^{\zeta}=0$ is equivalent to $e_1\cdot f+e_2\cdot g+e_3\cdot h= e_2\cdot f -e_1\cdot g=e_3\cdot g -e_2\cdot h=e_1\cdot h -e_3\cdot f=0$; 
the latter three equalities, 
plugged into equations \eqref{eqn_laplfgh}, exactly give $\Delta_{\mathbb{S}^3}f=\Delta_{\mathbb{S}^3}g=\Delta_{\mathbb{S}^3}h=0$, hence: 
$f$, $g$ and $h$ are constant. 
\cqfd

 \begin{rmk}
  We have not used the $\Gamma$-invariance of the tensors here; 
  nonetheless, since the $\theta_j$ are $\SU(2)$-invariant, which comes from the identities $\theta_j=dd^c_{I_j}(\tfrac{1}{r^2})$, 
  this does not give us any further information. 
 \end{rmk}

 \subsubsection{Lower order variation of the Kähler forms: determination of the coefficients}
We know form the formal expansion of $g_{\zeta}$ (or those of $g_{\zeta'}$ and $g_{\zeta''}$) 
that the $a_{jk}$ coefficients of Proposition \ref{prop_varpishape} are quadratic homogeneous polynomials in their arguments. 
Their explicit form is given as follows. 
 \begin{prop}   \label{prop_varpicoefs}
  With the same notations as in Proposition \ref{prop_varpishape}, 
    \begin{align*}
     a_{31}(\zeta'') &=  0,                    &  a_{32}(\zeta'')=& \,0,                                        
      & a_{32}(\zeta'')=&-\|\Gamma\||\zeta_3|^2, \\
     a_{21}(\zeta')  &=  0,                    &  a_{22}(\zeta' )=& \,-\|\Gamma\||\zeta_2|^2,                    
      & a_{32}(\zeta'')=&-\|\Gamma\|\langle\zeta_2,\zeta_3\rangle, \\
     a_{11}(\zeta)   &= -\|\Gamma\||\zeta_1|^2, &  a_{12}(\zeta  )=& \,-\|\Gamma\|\langle\zeta_1,\zeta_2\rangle,          
      & a_{32}(\zeta'')=&-\|\Gamma\|\langle\zeta_1,\zeta_3\rangle.
    \end{align*}
  where $\|\Gamma\|:=\tfrac{|\Gamma|}{2\vl(\mathbb{B}^4)}=\tfrac{|\Gamma|}{\pi^2}$. 
 \end{prop}
 \prf. 
 We shall first prove the assertion on the $a_{3j}(\zeta'')$, 
 and then apply the same techniques to determine the $a_{2j}(\zeta')$ -- 
 the $a_{1j}(\zeta)$ being dealt with in a similar way. 
 One more time we can assume that $\zeta$ is chosen in $\mathfrak{h}-D$, 
 and so that $\zeta', \zeta''\notin D$. 
 
 ~

 \noindent\textit{The coefficient $a_{33}(\zeta'')$}. 
 To begin with, 
 set $a=a_{31}(\zeta'')$, $b=a_{32}(\zeta'')$ and $c=a_{33}(\zeta'')$. 
 We consider on $X_{\zeta''}$ (which is smooth by our assumption $\zeta''\notin D$) 
 a closed form $\lambda$ with compact support representing $\zeta_3$ by Poincaré duality; 
 this is possible since minimal resolutions of $\C^2/\Gamma$ have compactly supported cohomology 
 \cite[Thm 8.4.3]{joy}, 
 and $X_{\zeta''}$ is diffeomorphic to such a resolution 
 (this is actually a minimal resolution of $(\C^2/\Gamma,I_3)$, 
 but we will not use this fact). 
 Next, consider a smooth cut-off function $\chi$, 
 vanishing on $(-\infty,1]$, equal to 1 on $[2, +\infty)$. 
 From the equality $\omega_3^{\e}=\tfrac{1}{2}dd^c_{I_3}(r^2)$, and from formulas \eqref{eqn_theta_j},  
 we have that 
  \begin{equation*} 
   \vareps:=\omega_3^{\zeta''}-\lambda-d\bigg[\frac{1}{4}I_3d\big(\chi(r)r^2\big) + 
                                       \frac{1}{4}(aI_1+bI_2+cI_3)d\big(\chi(r)r^{-2}\big)\bigg]
  \end{equation*}
is well-defined on $X_{\zeta'}$, 
has cohomology class 0, and is $O(r^{-6})$ at infinity, with appropriate decay on its derivatives; 
here we write $r$ instead of $(\lambda_3^{\zeta''})^*r$. 
As we need it further, we shall also see now that $\vareps$ admits a primitive decaying at infinity.  
 
From \cite[Thm 8.4.1]{joy}, 
$\vareps$ can indeed be written as $h+d\beta+d^{*_{g_{\zeta''}}}\gamma$, 
where $h$ is in $C^{\infty}_{3}(X_{\zeta},\Lambda^2)$ and is $g_{\zeta''}$-harmonic, 
and $\beta$ and $\gamma$ are in $C^{\infty}_{2}(X_{\zeta},\Lambda^2)$; 
we used here classical notations for weighted spaces: 
for example, $\beta=O(r^{-2})$, $\nabla^{\e}\beta=O(r^{-3})$ and so on. 
The harmonic form $h$ is actually decaying fast enough so that we can say it is closed and co-closed: 
write  (all the operations and tensors are computed with respect to $g_{\zeta''}$) for all $r$
 \begin{equation*}
  0=\int_{\mathbb{B}_{X_{\zeta''}}(r)} (h, \Delta h) \vol 
   =\int_{\mathbb{B}_{X_{\zeta''}}(r)} \big(|dh|^2+|d^{*}h|^2\big) \vol+ 
    \int_{\mathbb{S}_{X_{\zeta''}}(r)} (h\odot dh+ h\odot d^{*}h) \vol,
 \end{equation*}
where $\mathbb{B}_{X_{\zeta''}}(r)=(\lambda_3^{\zeta''})^{-1}\big(\mathbb{B}^4(r)/\Gamma\big)$, 
and $\mathbb{S}_{X_{\zeta''}}(r)$ is its boundary. 
From what precedes, the boundary integral is easily seen to be $O(r^{3-3-4})=O(r^{-4})$, and thus $dh=d^*h=0$. 
Hence $0=d\vareps=dd^{*}\gamma$; 
an integration by parts similar to the previous, but with boundary term of size $O(r^{-2})$, leads us to $d^*\gamma=0$, 
and thus $\vareps=h+d\beta$. 
According to \cite[thm 8.4.1]{joy} again, $\mathcal{H}^2(X_{\zeta''}) \to H^2(X_{\zeta''})$, 
$h\mapsto[h]$ is an isomorphism; 
now here $[h]=[\vareps-d\beta]=0$. 
Therefore $h=0$, and $\vareps=d\beta$, with $\beta=O(r^{-2})$. 

Write $\mathbb{B}(r)$ for $\mathbb{B}_{X_{\zeta''}}(r)$ to simplify notations; 
we shall now compute the integrals $\int_{\mathbb{B}(r)}(\omega_3^{\zeta''})^2$ in two different ways.  
First, recall that $(\omega_3^{\zeta''})^2=2\vol^{\e}$, 
and thus\footnote{ 
Indeed, $\int_{\mathbb{B}(r)}(\omega_3^{\zeta''})^2$ is the limit as $s$ goes to 0 of 
$\int_{\mathbb{B}(r)-\mathbb{B}(s)}(\omega_3^{\zeta''})^2$,
since as an $s$-tubular neighbourhood of $E:=(\lambda_3^{\zeta''})^{-1}\big(\{0\}\big)$ which is of real dimension 2, 
$\mathbb{B}(s)$ has its volume tending to 0 when $s$ goes to 0. 
Now we can also see $\int_{\mathbb{B}(r)-\mathbb{B}(s)}(\omega_3^{\zeta''})^2$ on $\R^4/\Gamma$ via $\lambda_3^{\zeta''}$ which is diffeomorphic away from $E$, 
and since $(\lambda_3^{\zeta''})_*(\omega_3^{\zeta''})^2=2\Omega_{\e}$, 
$\int_{\mathbb{B}(r)-\mathbb{B}(s)}(\omega_3^{\zeta''})^2$ is twice the euclidean volume of the annulus of radii $s$ and $r$ in $\R^4/\Gamma$, 
hence the result when $s\to 0$.}
$\int_{\mathbb{B}(r)}(\omega_3^{\zeta''})^2= \frac{2r^4}{|\Gamma|}\vl\big(\mathbb{B}^4\big)$. 
On the other hand, since $\omega_3^{\zeta''}=\lambda+d\varphi+\vareps$, 
with $\varphi=\tfrac{1}{4}I_3d\big(\chi(r)r^2\big) + \tfrac{1}{4}(aI_1+bI_2+cI_3)d\big(\chi(r)r^{-2}\big)$, 
we have:
 \begin{equation}    \label{eqn_intgrlsum1}
  \begin{aligned}
   \int_{\mathbb{B}(r)}(\omega_3^{\zeta''})^2 = &
     \int_{\mathbb{B}(r)} \lambda^2 +2\int_{\mathbb{B}^4(r)/\Gamma}\lambda\wedge d\varphi
      + 2\int_{\mathbb{B}(r)}\lambda\wedge\vareps                                           \\
     &\qquad + \int_{\mathbb{B}(r)} (d\varphi)^2
             + 2\int_{\mathbb{B}(r)}\vareps\wedge d\varphi + \int_{\mathbb{B}^4(r)/\Gamma}\vareps^2. 
  \end{aligned}
 \end{equation}
Let us analyse those summands separately. 

For $r$ large enough, 
$\int_{\mathbb{B}(r)}\lambda^2 = \int_{X_{\zeta''}}\lambda^2=\lambda\cup\lambda= -|\zeta_3|^2$, 
by Lemma \ref{lem_cupprod}, and the fact that $[\lambda]=[\omega_3^{\zeta''}]=(\zeta'')_3= \zeta_3$. 

The integral $\int_{\mathbb{B}(r)}\lambda\wedge d\varphi$ equals 
$\int_{\mathbb{S}(r)} \lambda\wedge\varphi$ by Stokes'
theorem, where $\mathbb{S}(r)$ stands for $\mathbb{S}_{X_{\zeta''}}(r)$,  
and this vanishes for $r$ large enough; 
similarly, $\int_{\mathbb{B}(r)}\lambda\wedge\vareps=\int_{\mathbb{B}(r)}\lambda\wedge d\beta= 
\int_{\mathbb{S}(r)}\lambda\wedge \beta=0$ for $r$ large enough. 

We now come to $\int_{\mathbb{B}(r)} (d\varphi)^2$. 
By Stokes, this is equal to $ \int_{\mathbb{S}(r)} d\varphi\wedge \varphi$, 
which we view back on $\R^4/\Gamma$ via $\lambda_3^{\zeta''}$. 
For $r\geq 2$, the integrand is
 \begin{align*}
  (\omega_3^{\e}+a\theta_1+b\theta_2+c\theta_3)\wedge 
     \Big[\frac{1}{2}\alpha_3-&\frac{1}{2}\big(a\frac{\alpha_1}{r^4}+b\frac{\alpha_2}{r^4}+c\frac{\alpha_3}{r^4}\big)\Big]  \\
    & = \frac{1}{2r^2}\Big(1-\frac{c}{r^4}\Big)\alpha_1\wedge\alpha_2\wedge\alpha_3+O(r^{-7}),   
 \end{align*}
since $\omega_3^{\e}\wedge\alpha_3=\tfrac{1}{r^2}\alpha_1\wedge\alpha_2\wedge\alpha_3$, 
$\omega_3^{\e}\wedge\alpha_1=\tfrac{1}{r^2}rdr\wedge\alpha_3\wedge\alpha_1=0$ and 
$\omega_3^{\e}\wedge\alpha_2=\tfrac{1}{r^2}rdr\wedge\alpha_3\wedge\alpha_2=0$ on $\mathbb{S}^3(r)/\Gamma$; 
$\theta_3\wedge\alpha_3= -\tfrac{\alpha_1\wedge\alpha_2\wedge\alpha_3}{r^6}$, 
$\theta_1\wedge\alpha_3=\tfrac{rdr\wedge\alpha_1\wedge\alpha_3}{r^6}=0$, 
$\theta_2\wedge\alpha_3=\tfrac{rdr\wedge\alpha_2\wedge\alpha_3}{r^6}=0$ on $\mathbb{S}^3(r)/\Gamma$; 
and $\theta_j\wedge\alpha_k=O(r^{-7})$ for $j,k=1,2,3$. 
Observe that $\int_{\mathbb{S}^3(r)/\Gamma}\alpha_1\wedge\alpha_2\wedge\alpha_3=4r^2\vl\big(\mathbb{B}^4(r)\big)/|\Gamma| $ 
-- for instance, compute $\int_{\mathbb{B}^4(r)}(\omega_3^{\e})^2=\frac{1}{2}\int_{\mathbb{B}^4(r)}\omega_3^{\e}\wedge d\alpha_3$ 
by Stokes; 
we thus end up with:
 \begin{equation*}
  \int_{\mathbb{B}(r)} (d\varphi)^2 
              = 2\big(1-\frac{c}{r^4}\big)\frac{\vl\big(\mathbb{B}^4(r)\big)}{|\Gamma|}+O(r^{-4})
              = \frac{2(r^4-c)}{|\Gamma|}\vl\big(\mathbb{B}^4\big)+O(r^{-4}).
 \end{equation*} 

We conclude by the last two summands of \eqref{eqn_intgrlsum1}. 
On the one hand, $\int_{\mathbb{B}(r)}\vareps\wedge d\varphi= \int_{\mathbb{S}(r)}\vareps\wedge \varphi= 
O(r^{3-6+1})=O(r^{-2})$, since $\vareps=O(r^{-6})$ and $\varphi=O(r)$. 
On the other hand, 
$\int_{\mathbb{B}(r)}\vareps^2=\int_{\mathbb{B}(r)}\vareps\wedge d\beta
  =\int_{\mathbb{S}(r)}\vareps\wedge \beta $: this is $O(r^{3-6-2})=O(r^{-5})$ 
  (and this is actually the only place where we need an estimate on the decay of a primitive of $\vareps$). 
  
Collecting the different estimates and letting $r$ go to $\infty$, we have:
 \begin{equation*}
  \frac{2r^4}{|\Gamma|}\vl\big(\mathbb{B}^4\big)= 
       -|\zeta_3|^2 + \frac{2\big(r^4-c\big)}{|\Gamma|}\vl\big(\mathbb{B}^4\big),
 \end{equation*}
that is: $c= -\tfrac{|\Gamma|}{2\vl(\mathbb{B}^4)}|\zeta_3|^2$. 

~

\noindent\textit{The coefficients $a_{31}(\zeta'')$ and $a_{31}(\zeta'')$}. 
We use the same techniques to see that $a_{31}(\zeta'')$ and $a_{32}(\zeta'')$, 
used above as $a$ and $b$ in $\varphi$, vanish. 
Recall for example that $\omega_1^{\zeta''}=\omega_1^{\e}$; 
therefore, $\mu:=\omega_1^{\zeta''}- d\big[\tfrac{1}{4}I_1d\big(\chi(r)r^2\big)\big]$ 
has compact support in $X_{\zeta''}$, and cohomology class $[\omega_1^{\zeta''}]= (\zeta'')_1=0$. 

We can compute $\int_{\mathbb{B}(r)}\omega_3^{\zeta''}\wedge\omega_1^{\zeta''}$ in two different ways. 
First, this is 0 for all $r$, since $\omega_1^{\zeta''}$ and $\omega_3^{\zeta''}$ are Kähler forms for the same hyperkähler 
metric and anti-commuting complex structures, and thus $\omega_3^{\zeta''}\wedge\omega_1^{\zeta''}\equiv0$. 
Secondly, using the sums $\omega_3^{\zeta''}=\lambda+d\varphi+\vareps$ 
and $\omega_1^{\zeta''}=\mu+d\psi$, with $\psi=\tfrac{1}{4}I_1d\big(\chi(r)r^2\big)$, we write for all $r$: 
 \begin{equation}    \label{eqn_intgrlsum2}
  \begin{aligned}
   \int_{\mathbb{B}(r)}\omega_3^{\zeta''}\wedge\omega_1^{\zeta''} = &
     \int_{\mathbb{B}(r)} \lambda\wedge\mu 
      + \int_{\mathbb{B}(r)}d\varphi\wedge \mu
      + \int_{\mathbb{B}(r)}\vareps \wedge \mu                                         \\
     &\qquad + \int_{\mathbb{B}(r)} \lambda\wedge d\psi
      + \int_{\mathbb{B}(r)}d\varphi\wedge d\psi + \int_{\mathbb{B}(r)}\vareps\wedge d\psi. 
  \end{aligned}
 \end{equation}
We briefly examine each summand. 
The first integral, $\int_{\mathbb{B}(r)} \lambda\wedge\mu$, 
is equal to $\int_{X_{\zeta''}} \lambda\wedge\mu$ for $r$ large enough, and this is 0, since $[\mu]=0$. 
The second integral can be rewritten as $\int_{\mathbb{S}(r)}d\varphi\wedge \lambda$, 
which vanishes for $r$ large enough; the same is true for the third (since $\vareps=d\beta$) and fourth integrals. 
As for the sixth integral, it can be written as $ \int_{\mathbb{S}(r)}\vareps\wedge \psi$, 
and this is $O(r^{3-6+1})=O(r^{-2})$. 

We are left with the fifth summand of \eqref{eqn_intgrlsum2}, 
which we rewrite as $\int_{\mathbb{S}(r)} d\varphi\wedge\psi$, and view back on $\R^4/\Gamma$. 
For $r\geq2$, the integrand is 
 \begin{equation*}
  (\omega_3^{\e}+a\theta_1+b\theta_2+c\theta_3)\wedge \frac{1}{2}\alpha_1= -a\frac{\alpha_1\wedge\alpha_2\wedge\alpha_3}{2r^6}
 \end{equation*}
on $\mathbb{S}^3(r)$, since there $\omega_3^{\e}\wedge\alpha_1= \theta_2\wedge\alpha_1=\theta_3\wedge\alpha_1= 0$. 
As a consequence, 
$\int_{\mathbb{B}(r)}d\varphi\wedge d\psi= -a\vl(\mathbb{B}^4)/|\Gamma|$. 

This tells us that the right-hand-side of \eqref{eqn_intgrlsum2} is $-a\vl(\mathbb{B}^4)/|\Gamma|$ 
when we let $r$ go to $\infty$; as the left-hand-side is always 0, $a_{31}(\zeta'')=a=0$. 
One proves that $a_{32}(\zeta'')=0$ in the same way. 

~

\noindent\textit{The $a_{2j}(\zeta')$ coefficients}. 
Let us come now to the $a_{2j}(\zeta')$ coefficients. 
Since $I_2^{\zeta'}=I_2^{\zeta''}$ and 
since we have the equality of $I_2^{\zeta'}$-holomorphic (2,0)-forms 
$\omega_3^{\zeta''}+i\omega_1^{\zeta''}=\omega_3^{\zeta'}+i\omega_1^{\zeta'}$, 
we have the writing:
 \begin{equation*}
  \omega_1^{\zeta'} = \mu+d\psi  \qquad \text{and} \qquad
  \omega_3^{\zeta'} = \lambda+ d\varphi+\vareps; 
 \end{equation*}
taking $\nu$ a compactly supported closed 2-form representing $\zeta_2=[\omega_2^{\zeta'}]$ (since $(\zeta')_2=\zeta_2$), 
we can write, for the same reasons as those invoked when proving the analogous formula for $\omega_3^{\zeta'}$: 
 \begin{equation*}
  \omega_2^{\zeta'} = \nu + d\xi + \eta
 \end{equation*}
with $\xi=\tfrac{1}{4}\big[I_2d\big(\chi(r)r^2\big)+(aI_1+bI_2+cI_3)d\big(\chi(r)r^{-2}\big)\big]$ 
where this time, $a=a_{21}(\zeta')$, $b=a_{22}(\zeta')$, $c=a_{23}(\zeta')$, 
and with $\eta\in C^{\infty}_6(X_{\zeta'},\Lambda^2)$ of the form $d\gamma$ 
for some $\gamma\in C^{\infty}_2(X_{\zeta'},\Lambda^2)$. 

Exactly as in what precedes, we get from computing respectively 
$\int_{\mathbb{B}(r)}(\omega_2^{\zeta'})^2$ and 
$\int_{\mathbb{B}(r)} \omega_2^{\zeta'}\wedge\omega_1^{\zeta'}$, 
where now $\mathbb{B}(r)=(\lambda_3^{\zeta''}\circ\lambda_2^{\zeta'})^{-1}\big(\mathbb{B}^4(r)/\Gamma\big)$,  
that $a_{22}(\zeta')=b=-\tfrac{|\Gamma|}{2\vl(\mathbb{B}^4)}|\zeta_2|^2$ and $a_{21}(\zeta')=a=0$. 
The slightly new computation is that of $\int_{\mathbb{B}^4(r)/\Gamma} \omega_2^{\zeta'}\wedge\omega_3^{\zeta'}$, 
though it goes through similar lines. 
Indeed, $\omega_2^{\zeta'}\wedge\omega_3^{\zeta'}=0$ identically, whereas: 
 \begin{equation}    \label{eqn_intgrlsum3}
  \begin{aligned}
   \int_{\mathbb{B}(r)}\omega_2^{\zeta'}\wedge\omega_3^{\zeta'} = &
     \int_{\mathbb{B}(r)} \nu\wedge \lambda
      + \int_{\mathbb{B}(r)} \nu\wedge d\varphi
      + \int_{\mathbb{B}(r)}\nu \wedge \vareps                                         \\
     &+ \int_{\mathbb{B}(r)a}  d\xi\wedge\lambda
      + \int_{\mathbb{B}(r)}d\xi\wedge d\varphi 
      + \int_{\mathbb{B}(r)} d\xi\wedge\vareps                                \\
     &+ \int_{\mathbb{B}(r)} \eta\wedge\lambda 
      + \int_{\mathbb{B}(r)} \eta\wedge d\varphi 
      + \int_{\mathbb{B}(r)} \eta\wedge \vareps.
  \end{aligned}
 \end{equation}
The first integral of the right-hand-side equals 
$\int_{X_{\zeta'}} \nu\wedge \lambda= -\langle[\nu],[\lambda]\rangle=-\langle\zeta_2,\zeta_3\rangle$ for $r$ large enough. 
For $r$ large enough again, the second, the third, the fourth and the seventh integrals vanish (use Stokes' theorem). 
Further uses of Stokes' theorem and the estimates $\varphi,\xi=O(r)$, 
$\vareps,\eta=O(r^{-6})$ and $\beta,\gamma=O(r^{-2})$ (with $\vareps=d\beta$, $\eta=d\gamma$) 
give us that $\int_{\mathbb{B}(r)} d\xi\wedge\vareps$ and $\int_{\mathbb{B}(r)} \eta\wedge d\varphi$ 
are $O(r^{-2})$, and that $\int_{\mathbb{B}(r)} \eta\wedge \vareps=O(r^{-5})$. 

We hence are left with computing the contribution of the fifth summand, 
namely $\int_{\mathbb{B}(r)}d\xi\wedge d\varphi $, 
in the right-hand-side of \eqref{eqn_intgrlsum3}. 
This can be rewritten as $\int_{\mathbb{S}(r)}\xi\wedge d\varphi$; seen on $\R^4/\Gamma$, 
the integrand is:
 \begin{equation*}
  \frac{1}{2}\Big[\alpha_2-\frac{1}{r^4}\big( a_{22}(\zeta')\alpha_2+c\alpha_3\big)\Big]
    \wedge \big(\omega_3^{\e}+ a_{33}(\zeta')\theta_3\big). 
 \end{equation*}
All computations done, this can be rewritten on $\mathbb{S}^3(r)/\Gamma$ as: 
$-c\tfrac{\alpha_1\wedge \alpha_2\wedge \alpha_3}{2r^6}$. 
Thus $\int_{\mathbb{B}(r)} d\xi\wedge d\varphi= -\tfrac{2c\vl(\mathbb{B}^4)}{|\Gamma|}+O(r^{-4})$ 
(recall that $\int_{\mathbb{S}^3(r)/\Gamma}\alpha_1\wedge \alpha_2\wedge \alpha_3=\tfrac{4r^6\vl(\mathbb{B}^4)}{|\Gamma|}$).

Collecting the estimates of the last two paragraphs and letting $r$ go to $\infty$, 
we get: $0=-\langle\zeta_2,\zeta_3\rangle -\tfrac{2c\vl(\mathbb{B}^4)}{|\Gamma|}$, 
that is: $c=-\tfrac{|\Gamma|}{2\vl(\mathbb{B}^4)}\langle\zeta_2,\zeta_3\rangle$. 

~

\noindent\textit{The $a_{1j}(\zeta)$ coefficients}. 
After noticing that $\omega_2^{\zeta}+i\omega_3^{\zeta}= \omega_2^{\zeta'}+i\omega_3^{\zeta'}$, 
and that $\omega_1^{\zeta}$ can be written as a sum 
 \begin{equation*}
  \mu + d\psi + \varsigma,
 \end{equation*}
where $\mu$ still has compact support but this time has class $[\omega_1^{\zeta}]=\zeta_1$, 
$\psi$ is the 1-form $\tfrac{1}{4}\big[I_1d\big(\chi(r)r^2\big)+(a_{11}(\zeta)I_1+a_{12}(\zeta)I_2+a_{13}(\zeta)I_3)d\big(\chi(r)r^{-2}\big)\big]$, 
and $\varsigma\in C^{\infty}_6(X_{\zeta},\Lambda^2)$ can be written $d\sigma$ with $\sigma\in C^{\infty}_2(X_{\zeta},\Lambda^2)$, 
we get that 
 \begin{align*}
  a_{11}(\zeta)= -\frac{|\Gamma|}{2\vl(\mathbb{B}^4)}&|\zeta_1|^2, 
   \qquad               a_{12}(\zeta)= -\frac{|\Gamma|}{2\vl(\mathbb{B}^4)}\langle\zeta_1,\zeta_2\rangle \\
    \text{and }\quad    a_{13}(\zeta)=& -\frac{|\Gamma|}{2\vl(\mathbb{B}^4)}\langle\zeta_1,\zeta_3\rangle,
 \end{align*}
from respective  computations of $\int_{\mathbb{B}(r}(\omega_1^{\zeta'})^2$, 
$\int_{\mathbb{B}(r)}\omega_1^{\zeta'}\wedge\omega_2^{\zeta'}$ and 
$\int_{\mathbb{B}(r)}\omega_1^{\zeta'}\wedge\omega_3^{\zeta'}$, 
where here $\mathbb{B}(r)$ is seen in $X_{\zeta}$, following the same lines as above. 
\cqfd

 \subsubsection{Conclusion: proof of Theorem \ref{thm_ALE} (general $\Gamma$)}

Let us sum up the situation. 
If we take $\Phi_{\zeta}=F_{\zeta}^{-1}:X_{\zeta}\backslash F_{\zeta}^{-1} (\{0\})\to (\R^4\backslash\{0\})/\Gamma$ 
and keep the notations introduced in this section, 
we have ${\Phi_{\zeta}}_{*}g_{\zeta}=\e+h_{\zeta}+O(r^{-6})$, 
${\Phi_{\zeta}}_{*}I_1^{\zeta}=I_1+\iota_1^{\zeta}+O(r^{-6})$ 
and ${\Phi_{\zeta}}_{*}\omega_1^{\zeta}=\omega_1^{\e}+\varpi_1^{\zeta}+O(r^{-6})$. 
The $I_1$-hermitian component of $h_{\zeta}$ is $\varpi_1^{\zeta}(\cdot,I_1\cdot)$, which we know, 
and its $I_1$-skew-hermitian component is $h_{\zeta'}$. 
Now the $I_2$-hermitian component of $h_{\zeta'}$ is $\varpi_2^{\zeta'}(\cdot,I_2\cdot)$, which we also know, 
and its $I_2$-skew-hermitian component is $h_{\zeta''}$. 
Finally, $h_{\zeta''}$ is $I_3$-hermitian, equal to $\varpi_3^{\zeta''}(\cdot,I_3\cdot)$, which we know as well. 
In a nutshell, we are able to write down explicitly $h_{\zeta}$: 
 \begin{align*}
  h_{\zeta}  = & \varpi_3^{\zeta''}(\cdot,I_3\cdot) + \varpi_2^{\zeta'}(\cdot,I_2\cdot) + \varpi_1^{\zeta}(\cdot,I_1\cdot) \\
             = & -\|\Gamma\| \bigg(\sum_{j=1}^3|\zeta_j|^2\theta_j(\cdot, I_j\cdot)
                         + \sum_{1\leq j<k\leq3}\langle\zeta_j,\zeta_j\rangle \theta_j(\cdot, I_j\cdot)\bigg), 
 \end{align*} 
which gives exactly formula \eqref{eqn_hzeta}, with $\cc=\frac{1}{2\vl(\mathbb{B}^4)}=\pi^{-2}$. 

From this and the formula for $\varpi_1^{\zeta}$ proved in \ref{prop_varpicoefs} 
-- which gives formula \eqref{eqn_varpi1zeta} of Theorem \ref{thm_ALE} --, 
we deduce the expected formula for $\iota_1^{\zeta}$.  
We know indeed that $\iota_1^{\zeta}= \iota_1^{\zeta'}$, 
and that $h_{\zeta'}= \omega_1^{\e}(\cdot,\iota_1^{\zeta'}\cdot)$ and $\iota_1^{\zeta'}$ is $\e$-symmetric, hence
 \begin{align*}
  \e(\iota_1^{\zeta}\cdot,\cdot) =& \e(\cdot,\iota_1^{\zeta'}\cdot) = -\omega_1^{\e}(I_1\cdot, \iota_1^{\zeta'}\cdot)= -h_{\zeta'}(I_1\cdot,\cdot) \\
                                 =& \|\Gamma\|\big(|\zeta_3|^2\theta_3(I_1\cdot, I_3\cdot)+ |\zeta_2|^2\theta_2(I_1\cdot, I_2\cdot) 
                                          + \langle\zeta_2,\zeta_3\rangle \theta_3(I_1\cdot,I_2\cdot ) \big)                                  \\
                                 =& \|\Gamma\|\big(|\zeta_3|^2\theta_3(\cdot, I_2\cdot)- |\zeta_2|^2\theta_2(\cdot, I_3\cdot) 
                                          - \langle\zeta_2,\zeta_3\rangle \theta_3(\cdot,I_3\cdot ) \big)                                     \\
                                 =& \|\Gamma\||\zeta_3|^2\frac{\alpha_2\cdot\alpha_3-rdr\cdot\alpha_1}{r^6} 
                                    - \|\Gamma\||\zeta_2|^2\frac{\alpha_2\cdot\alpha_3+ rdr\cdot\alpha_1}{r^6} \\
                                  &  - \|\Gamma\|\langle\zeta_2,\zeta_3\rangle\frac{(rdr)^2+\alpha_3^2-\alpha_1^2 - \alpha_{2}^2}{r^6},       
 \end{align*}
of which formula \eqref{eqn_iota1zeta} is just a rewriting.

 \subsection{Vanishing of the third order terms when $\Gamma$ is not cyclic} \label{section_ord3term}

We shall see in this section that in the expansion $g_{\zeta} = \e + h_{\zeta} + \sum_{j=3}^{\infty} h_{\zeta}^{(j)}$, 
\textit{when $\Gamma$ is one of the $\mathcal{D}_k$, $k\geq2$, or contains one of these as is the case when $\Gamma$ 
is binary tetrahedral, octahedral or icosahedral}, then the third order term $h_{\zeta}^{(3)}$ vanishes, 
and that this holds as well for complex structures and Kähler forms. 
Keeping working the diffeomorphisms $F_{\zeta}$ of the previous part even if we omit them to simplify notations, 
we claim:
 \begin{prop}   \label{prop_ord3term}
  Suppose $\Gamma$ contains $\mathcal{D}_k$, $k\geq 2$, as a subgroup. 
  Then $g_{\zeta} = \e + h_{\zeta} + O(r^{-8})$, $I_1^{\zeta}= I_1 + \iota_1^{\zeta} +O(r^{-8})$, 
  $\omega_1^{\zeta}= \omega_1 + \varpi_1^{\zeta} + O(r^{-8})$, 
  where by $O(r^{-8})$ we mean tensors whose $\ell$th-order derivatives (for $\nabla^{\e}$) are $O(r^{-8-\ell})$. 
 \end{prop}
 \prf. 
  We shall first see that, for a general $\Gamma$, 
  the crucial considerations made in section \ref{section_hzeta} on the second order term $h_{\zeta}$ of the expansion of $g_{\zeta}$ 
  still hold for $h_{\zeta}^{(3)}$; 
  first recall that $h_{\zeta}^{(3)}$ is a homogeneous polynomial of $\zeta$ of order 3, with values in tensors $O(r^{-6})$, 
  with according decay on the derivatives, \textit{independent of $\zeta$}. 
  We start by the claim that: 
   \begin{equation*}
     \tr^{\e}(h_{\zeta}^{(3)}) = 0 \qquad \text{ and } \qquad \delta^{\e} h_{\zeta}^{(3)} = 0.
   \end{equation*}
  Indeed for the trace assertion, once $\zeta\in \mathfrak{h}\otimes\R^3-D$ is fixed, 
  one has for all $t$ that
   \begin{align*}
    \Omega_{\e} &= \vol^{g_{t\zeta}}=\det^{\e}\big(\e+ t^2 h_{\zeta} + t^3h_{\zeta}^{(3)} + O(t^4) \big)\Omega_{\e} \\
                &= \big(1 + t^2 \tr^{\e}(h_{\zeta}) + t^3\tr^{\e}(h_{\zeta}^{(3)}) + O(t^4)\big)\Omega_{\e}, 
   \end{align*}
 since the higher order contributions of $t^2h_{\zeta}$ 
 are included in the $O(t^4)$, hence the claim.  
 
 We thus see that $h_{\zeta}^{(3)}$ shares this property with $h_{\zeta}$ 
 because the non-linear contributions of the $h_{t\zeta}$, which are of order at least 4 in $t$, 
 do not interfere with the linear contribution of $h_{t\zeta}^{(3)}$. 
 We thus generalize this observation to prove that $h_{\zeta}^{(3)}$ shares other properties with $h_{\zeta}$, and to start with, 
 that $\delta^{\e}h_{\zeta}^{(3)}=0$, as promised. 
 Again we proceed within three steps, considering first $\zeta''=(0,0,\zeta_3)$, 
 and then $\zeta'=(0,\zeta_2,\zeta_3)$ and $\zeta=(\zeta_1,\zeta_2, \zeta_3)$. 

 The case of $h_{\zeta''}^{(3)}$ is immediate, and merely amounts to the fact that it is an $I_3$-hermitian tensor (the $g_{t\zeta''}$ are) 
 with vanishing trace for $\e$, 
 used with the Kähler identity $[\Lambda_{\omega_3},d^c_{I_3}]=d^{*_{\e}}$ applied to $h_{\zeta''}^{(3)}(I_3\cdot,\cdot)$. 
 
 For the case of $h_{\zeta'}^{(3)}$, remember the following: 
 we first saw that the second order variation of $I_2^{\zeta'}=I_2^{\zeta''}$ was $\e$-symmetric; 
 this still holds for the third order term, since otherwise we would end up with a non-trivial $I_2$-entire function on $\C^2$, 
 decaying (like $r^{-6}$) at infinity, which is absurd. 
 Then we identified the $I_2$-skew-hermitian part of $h_{\zeta'}$ with $h_{\zeta''}$; 
 again, this holds for $h_{\zeta'}^{(3)}$ with $h_{\zeta''}^{(3)}$ (and the latter is indeed $I_2$-skew-hermitian). 
 This amounts to looking at the third order term in $t$ of: 
  \begin{itemize}
   \vspace{-3pt}
   \item the expansion of $g_{t\zeta''}=\omega_2(\cdot,I_2^{t\zeta''}\cdot)$ to see that $h_{\zeta''}^{(3)}$ is indeed $I_2$-skew-hermitian 
         (recall $\omega_2^{\zeta}=\omega_2$); 
   \item the expansion of $g_{t\zeta'}(I_2^{t\zeta'}\cdot,I_2^{t\zeta'}\cdot)-g_{t\zeta'}$ 
         to see that $\tfrac{1}{2}\big(h_{\zeta'}^{(3)}+h_{\zeta'}^{(3)}(I_2\cdot,I_2\cdot)\big)=h_{\zeta''}^{(3)}$. 
  \end{itemize}
 We concluded by using the usual Kähler identity (for $I_2$) on the $\e$-trace-free $I_2$-(1,1) form 
 $\tfrac{1}{2}\big(h_{\zeta'}(I_2\cdot,\cdot)-h_{\zeta'}(\cdot,I_2\cdot)\big)$, 
 after seeing it was closed; 
 we can do the same on its analogue $\tfrac{1}{2}\big(h_{\zeta'}^{(3)}(I_2\cdot,\cdot)-h_{\zeta'}^{(3)}(\cdot,I_2\cdot)\big)$, 
 which is also an $\e$-trace-free $I_2$-(1,1) form, 
 and is closed as seen when looking at the third order in $t$ of the expansion of 
 $\omega_2^{t\zeta'}=\tfrac{1}{2}\big(g_{t\zeta'}(I_2^{t\zeta'}\cdot,\cdot)-g_{t\zeta'}(I_2^{t\zeta'}\cdot,I_2\cdot)\big)$. 

 One deals with $h_{\zeta}$ in analogous way. 
 In particular, we get in passing that the third order variation of $I_1^{\zeta}=I_1^{\zeta'}$, 
 $\jmath_1^{\zeta}$ say, is $\e$-symmetric and anti-commutes to $I_1$, 
 that the $I_1$-skew-hermitian part of $h_{\zeta}^{(3)}$ is $h_{\zeta'}^{(3)}$, 
 related to $\jmath_1^{\zeta}$ by $h_{\zeta'}^{(3)}=\omega_1(\cdot,\jmath_1^{\zeta}\cdot)$, 
 and that its $I_1$-hermitian part gives rise to an $\e$-trace-free closed $I_1$-(1,1) form. 

 ~

 Running backward this description, 
 we will thus be done if we show that the third order variations of the Kähler forms vanish when $\Gamma$ contains a binary dihedral group. 
 In general though, we know these are $O(r^{-6})$ near 0 and infinity with corresponding decay on their derivatives, 
 that they are of type (1,1) for one of the $I_j$ and trace-free; 
 they are thus $*_{\e}$-anti-self-dual, and therefore can be written as $f\theta_1+g\theta_2+h\theta_3$, 
 where this time, $r^2f$, $r^2g$ and $r^2h$ depend only of the spherical coordinate of their argument. 
 Our form are moreover closed, hence in particular harmonic; 
 using again that the Laplace-Beltrami operator and the rough laplacian coincide on $(\R^4,\e)$, 
 and that the $\theta_j$ are harmonic, 
 we have this time:
  \begin{equation*}
   \Delta_{\e}(f\theta_1)= (\Delta_{\e}f)\theta_1- 2\sum_{k=0}^3 (e_k\cdot f)\nabla_{e_k}^{\e}\theta_1, 
  \end{equation*}
 with $e_0=\tfrac{x_j}{r}\tfrac{\partial}{\partial x_j}$. 
 We set $\tilde{f}=r^2f$; this is a function on $\mathbb{S}^3$, and 
 $e_0\cdot f= e_0\cdot(r^{-2}\tilde{f})= e_0\cdot(r^{-2})\tilde{f}= -2r^{-3}\tilde{f}= -2r^{-1}f$. 
 Since on functions, $\Delta_{\e}=-\tfrac{1}{r^3}\partial_r\big(r^3\partial_r\cdot\big)+\tfrac{1}{r^2}\Delta_{\mathbb{S}^3}$, 
 one has: 
  \begin{equation*}
   \Delta_{\e}f= \Delta_{\e}(r^{-2}\tilde{f}) = -\frac{1}{r^3}\partial_r\big(r^3\partial_r(r^{-2})\big)\tilde{f}+ \frac{1}{r^4}\Delta_{\mathbb{S}^3}\tilde{f}
               = \frac{1}{r^4}\Delta_{\mathbb{S}^3}\tilde{f}, 
  \end{equation*}
 since $\partial_r\big(r^3\partial_r(r^{-2})\big)=0$ ($r^{-2}$ is the Green function on $\R^4$). 

 Moreover, 
 $\nabla_{e_0}^{\e}\theta_1= \partial_r(r^{-6})(rdr\wedge\alpha_1-\alpha_2\wedge\alpha_3)+r^{-6}\nabla_{e_0}^{\e}(rdr\wedge\alpha_1-\alpha_2\wedge\alpha_3)
  =-\tfrac{6}{r}\theta_1+\tfrac{2}{r}\theta_1= -\tfrac{4}{r}\theta_1 $. 
 We recall that $\nabla_{e_1}^{\e}\theta_1=0$, $\nabla_{e_2}^{\e}\theta_1= \tfrac{2}{r}\theta_3$ and $\nabla_{e_3}^{\e}\theta_1= -\tfrac{2}{r}\theta_2$, 
 therefore: 
  \begin{equation*}
   \Delta_{\e}(f\theta_1)=\frac{1}{r^4}\big(\Delta_{\mathbb{S}^3}\tilde{f} -16\tilde{f}\big)\theta_1 
                          - \frac{2}{r^3}\big((e_2\cdot\tilde{f})\theta_3 - (e_3\cdot\tilde{f})\theta_2\big).  
  \end{equation*}
 Writing the analogous equations on $\tilde{g}=r^2g$, $\tilde{h}=r^2h$, 
 the equation $\Delta_{\e}(f\theta_1+g\theta_2+h\theta_3)=0$ is equivalent to the system: 
  \begin{equation}  \label{eqn_laplfgh2}
   \begin{aligned} 
    \Delta_{\mathbb{S}^3}\tilde{f} -16\tilde{f} - 4(e_3\cdot \tilde{g}) + 4(e_2\cdot \tilde{h}) &= 0,  \\
    \Delta_{\mathbb{S}^3}\tilde{g} -16\tilde{g} - 4(e_1\cdot \tilde{h}) + 4(e_3\cdot \tilde{f}) &= 0,  \\
    \Delta_{\mathbb{S}^3}\tilde{h} -16\tilde{h} - 4(e_1\cdot \tilde{f}) + 4(e_3\cdot \tilde{g}) &= 0.
   \end{aligned}
  \end{equation}
 Now the closure assertion on $f\theta_1+g\theta_2+h\theta_3$ is equivalent to 
 $(e_1\cdot f)+(e_2\cdot g)+(e_3\cdot h)= (e_0\cdot f)-(e_3\cdot g)+(e_2\cdot h)
                                        = (e_0\cdot g)-(e_1\cdot h)+(e_3\cdot f)
                                        = (e_0\cdot h)-(e_2\cdot f)+(e_1\cdot g)=0$. 
 
 Since $e_0\cdot u= -\tfrac{2}{r^3}\tilde{u}$ and $e_k\cdot u= \tfrac{1}{r^2}e_k\cdot\tilde{u}$ for $u=f,g,h$ and $k=1,2,3$, 
 we deduce from the latter equalities and the system \eqref{eqn_laplfgh2} the equations:
  \begin{equation*}  
   \Delta_{\mathbb{S}^3}\tilde{f} -8\tilde{f}= \Delta_{\mathbb{S}^3}\tilde{g} -8\tilde{g} =\Delta_{\mathbb{S}^3}\tilde{h} -8\tilde{h}=0. 
  \end{equation*}
 Setting $\hat{f}=r^2\tilde{f}$ and likewise for $\hat{g}$ and $\hat{h}$, 
 we get that $\hat{f}$, $\hat{g}$, $\hat{h}$ are harmonic (on the whole $r^4$) and homogeneous of degree 2. 
 This is not hard seeing that they are thus linear combinations of the $x_1^2-x_j^2$, $j=2,3,4$, and the $x_jx_k$, $1\leq j< k\leq 4$. 

 The $\theta_j$ are $\Gamma$-invariant; $f$, $g$ and $h$, and consequently $\hat{f}$, $\hat{g}$ and $\hat{h}$, 
 must thus be as well. 
 But \textit{if $\Gamma$ contains a binary dihedral group as a subgroup, 
 then there is no non-trivial linear combination of the above polynomials which is $\Gamma$-invariant}. 
 We use first the $\tau$-invariance; 
 if indeed $\mathcal{D}_k<\Gamma$ for some $k\geq2$ and $u=\sum_{j=1}^3 a_j (x_1^2-x_j^2) + \sum_{1\leq j< \ell\leq 4} a_{j\ell}x_jx_{\ell}$ is $\Gamma$-invariant, 
 then $2u=u+\tau^*u=a_2(x_1^2-x_2^2+x_3^2-x_4^2)+ a_3(x_1^2-x_3^2+x_3^2-x_1^2)+ a_4(x_1^2-x_4^2+x_3^2-x_2^2)
                    +a_{12}(x_1x_2+x_3x_4)+a_{13}(x_1x_3-x_3x_1)+a_{14}(x_1x_4-x_3x_2)
                    +a_{23}(x_2x_3-x_4x_1)+a_{24}(x_2x_4-x_4x_2)+a_{34}(x_3x_4+x_1x_2)$, 
 that is: $u$ has shape $a(x_1^2-x_2^2+x_3^2-x_4^2)+2b(x_1x_2+x_3x_4)+2c(x_1x_4-x_3x_2)$, 
 i.e. $a\mathfrak{Re}(z_1^2+z_2^2)+b\mathfrak{Im}(z_1^2+z_2^2)+c\mathfrak{Im}(\overline{z_1}z_2)$, 
 $a,b,c\in\R$, in complex notations. 
 We now use the $\zeta_k$-action and write: 
  \begin{align*}
   ku =& \sum_{\ell=0}^k\zeta_k^*u= \sum_{\ell=0}^k a(\zeta_k^{\ell})^*\mathfrak{Re}(z_1^2+z_2^2)
         + b(\zeta_k^{\ell})^*\mathfrak{Im}(z_1^2+z_2^2)+c(\zeta_k^{\ell})^*\mathfrak{Im}(\overline{z_1}z_2) \\
      =& \mathfrak{Re}\bigg(\sum_{j=0}^k a(e^{2i\ell\pi/k}z_1^2+e^{-2i\ell\pi/k}z_2^2)\bigg) \\ 
       & +  \mathfrak{Im}\bigg(\sum_{j=0}^k b(e^{2i\ell\pi/k}z_1^2 + e^{-2i\ell\pi/k}z_2^2) + ce^{-2i\ell\pi/k}\overline{z_1}z_2 \bigg) \\
      =& 0,        
  \end{align*}
 since $e^{2i\pi/k}\neq1$ ($k\geq2$). 
 In particular, the third order variation term of $\omega_1^{\zeta}$ vanishes. 
 In other words, $\omega_1^{\zeta}=\omega_1+\varpi_1^{\zeta}+O(r^{-8})$. 
 
 Moreover, since $\jmath_1^{\zeta}$ is determined by $h_{\zeta'}^{(3)}$ which is also 0, 
 this third order variation of complex structure vanishes as well, or: $I_1^{\zeta}=I_1+\iota_1^{\zeta}+O(r^{-8})$. 
 \cqfd
 
~

 This completes the proof of Theorem \ref{thm_ALE}. 
 Notice however that in view of the previous two sections, 
 we could also have given similar statements on the second and third complex structures and Kähler forms of $X_{\zeta}$. 
 We chose to focus and the first ones since this is what is needed in our construction of Part 1, 
 see in particular Lemma \ref{lem_omega1Y}, which is just a specialization  of Theorem \ref{thm_ALE}: 
 take $\zeta=\xi$ verifying conditions \eqref{eqn_xi}, and $\Phi_Y=\Phi_{\zeta}$. 
 
 Nonetheless, the asymptotics of the second and third complex Kähler forms are available via Proposition \ref{prop_varpicoefs}, 
 from which the asymptotics of the corresponding complex structures easily follow, 
 since  the asymptotics of the metric are known.

 \subsection{Comments on Lemma \ref{lem_beth}}  \label{section_cplxstr}
  \subsubsection{The condition \eqref{eqn_xi}}
   The first comment we want to make about Lemma \ref{lem_beth} concerns the reason why we state it under the condition \eqref{eqn_xi}, 
   which we can recall as $(|\zeta_2|^2-|\zeta_3|^2)+2i\langle\zeta_2,\zeta_3\rangle=0$ 
   (if one takes $\zeta$ instead of $\xi$ as parameter). 
   
   One could instead try to generalize the proof we give in section \ref{section_gluing} with help of the asymptotic given by Theorem \ref{thm_ALE},  
   with $\zeta$ a generic element of $\mathfrak{h}\otimes\R^3-D$. 
   This is formally possible, by leads to include terms such $\tfrac{1}{r^2z_1}$, 
   $\tfrac{1}{r^2\overline{z_1}}$, $\tfrac{1}{r^2z_2}$, $\tfrac{1}{r^2\overline{z_2}}$ in $\vareps_1$ and $\vareps_2$ of that proof, 
   which is obviously not compatible with the requirement that  $\beth$ is a diffeomorphism of $\R^4$. 

   In others words, $(|\zeta_2|^2-|\zeta_3|^2)+2i\langle\zeta_2,\zeta_3\rangle$ appears as an obstruction for $I_1^{\zeta}$ 
   to be approximated to higher orders by $I_1$, even with some liberty on the diffeomorphism between infinities of 
   $X_{\zeta}$ and $\R^4/\Gamma$, 
   which reveals some link between the parametrisation of the $X_{\zeta}$ and the general problem of the approximation of their complex structures.

  \subsubsection{Links with the parametrisation}
  Conversely we interpret of Lemma \ref{lem_beth} as follows: when $\Gamma=\mathcal{D}_k$ 
  -- this would be true also in the tetrahedral, octahedral and icosahedral cases --
  and $(|\zeta_2|^2-|\zeta_3|^2)+2i\langle\zeta_2,\zeta_3\rangle=0$, 
  then the complex structure $I_1^{\zeta}$ can be viewed as approximating the standard complex structure $I_1$ 
  with precision twice that of the general case, i.e. with an error $O(r^{-8})$ instead of $O(r^{-4})$, 
  up to an adjustment of the ALE diffeomorphism given in Kronheimer's contruction. 
  Now $(|\zeta_2|^2-|\zeta_3|^2)+2i\langle\zeta_2,\zeta_3\rangle=\langle\zeta_2+i\zeta_3,\zeta_2+i\zeta_3\rangle$, 
  and this precisely the coefficient $a_k$ in the equation of $X_{\zeta}$ seen as a submanifold of $\C^3$, 
  which is
   \begin{equation}  \label{eqn_defC2/Dk}
    u^2+v^2w+w^{k+1}= a_0 + a_1w +\cdots + a_kw^k +b v 
   \end{equation}
  ($a_j$ being given by symmetric functions of the $(k+2)$ first diagonal values of $\zeta_2+i\zeta_3\in \mathfrak{h}_{\C}$ or 
  $\mathfrak{h}_{\C}/(\text{Weyl group})$ of degree $(k+2-j)$, and $b$ by their Pfaffian). 
  
  Denote by $X_{\mathcal{D}_k}$ the orbifold defined in $\C^3$ by the equation 
  $u^2+v^2w+w^{k+1}=0$, i.e. equation \eqref{eqn_defC2/Dk} with $a_0=\cdots=a_k=b=0$. 
  This is identified to $\C^2/\mathcal{D}_k$ via the map  
  $(z_1,z_2)\mapsto(u,v,w)
     :=\big(\tfrac{1}{2}(z_1^{2k+1}z_2-z_2^{2k+1}z_1), \tfrac{i}{2}(z_1^{2k}+z_2^{2k}), z_1^2z_2^2\big)$. 
  This suggests that $(u,v,w)$ in equation \eqref{eqn_defC2/Dk} should somehow have respective degrees $2k+2$, $2k$ and $4$ 
  in the $z_1,z_2$ variables, and this equation remains homogeneous if we give formal degree 2 to $\zeta$.  
  When $a_k=0$, the right-hand-side member of \eqref{eqn_defC2/Dk} is therefore formally conferred "pure" degree at most $ 4k-4$,  
  instead of $4k$. 
  
  We suspect that this corresponds to the improvement by four orders in the approximation of $I_1^{\zeta}$ by $I_1$ 
  in the sense of Lemma \ref{lem_beth}. 
  It would thus be of interest to draw a rigorous picture out of these informal considerations, 
  establishing a more direct link of this kind between the parameter $\zeta$ and the associated complex structures, 
  without passing by the analysis of $g_{\zeta}$, 
  which we unfortunately have not been able to so far.

\begin{small}

 \renewcommand{\refname}{\large{References}}

\end{small}

\small \textsc{Max-Planck-Institut für Mathematik Bonn}

\url{auvray@mpim-bonn.mpg.de }


\begin{thebibliography}{99}
         \bibitem[Auv]{auv1} H. Auvray, \textit{From ALF to ALE gravitational instantons. I}, arXiv:1210.1654 [math.DG], preprint 2012. 
         \bibitem[BKN]{bkn} S. Bando, A. Kasue, H. Nakajima, \textit{On a construction of coordinates at infinity on manifolds with fast curvature decay and maximal volume growth}, 
                 Invent. Math. 97 (1989), no. 2, 313-349. 
         \bibitem[Bes]{bes} A.L. Besse, \textit{Einstein manifolds}, Ergebnisse der Mathematik und ihrer Grenzgebiete (3) [Results in Mathematics and Related Areas (3)], 10. 
                            Springer-Verlag, Berlin, 1987.
         \bibitem[Biq]{biq} O. Biquard, \textit{Differential Geometry and Global Analysis}, lecture notes available at \url{http://www.math.ens.fr/~biquard/dgga2007.pdf}. 
         \bibitem[BM]{bm} O. Biquard, V. Minerbe, \textit{A Kummer construction for gravitational instantons}, to appear in Comm. Math. Phys.
         \bibitem[BR]{br} O. Biquard, Y. Rollin, \textit{Smoothing singular extremal Kähler surfaces and minimal Lagrangians}, arXiv:1211.6957 [math.DG]. 
         \bibitem[BK]{bk} P. Buser, H. Karcher, \textit{Gromov's almost flat manifolds}, Astérisque, 81. Société Mathématique de France, Paris, 1981.
         \bibitem[CH]{ch-hi} S. Cherkis, N. Hitchin, \textit{Gravitational instantons of type $D_k$}, 
                 Comm. Math. Phys. 260 (2005), no. 2, 299-317.
         \bibitem[CK]{ch-ka} S. Cherkis, A. Kapustin, \textit{Singular monopoles and gravitational instantons of type $D_k$}, 
                 Comm. Math. Phys. 203 (1999), 713-728.
         \bibitem[GH]{gh} G.W. Gibbons, S.W. Hawking, \textit{Gravitational multi-instantons}, Phys. Lett. B78, 430-432 (1976).
         \bibitem[GLP]{glp} M. Gromov, J. Lafontaine, P. Pansu, 
                            \textit{Metric structures for Riemannian and non-Riemannian spaces}, Transl. from the french by Sean Michael Bates. 
                            With appendices by M. Katz, P. Pansu, and S. Semmes. Edited by J. Lafontaine and P. Pansu. (English). 
                            Birkhäuser, Basel, Progress in Mathematics (1999).
         \bibitem[GV]{gv} M. Gursky, J.A. Viaclovsky, \textit{Critical metrics on connected sums of Einstein four-manifolds}, arXiv:1303.0827 [math.DG], preprint 2013.
         \bibitem[Haw]{haw} S.W. Hawking, \textit{Gravitational instantons}, Phys. Lett. A 60 (1977), no. 2, 81-83. 
         \bibitem[Hei]{hein} H. J. Hein, \textit{On gravitational instantons}, PhD Thesis, Princeton University. 2010.
         \bibitem[Joy]{joy} D. Joyce, {\it Compact Manifolds with Special Holonomy}, Oxford Mathematical Monographs, Oxford University Press, 2000.
         \bibitem[Kro1]{kro1} P.B. Kronheimer, \textit{The construction of ALE spaces as hyper-Kähler quotients}, 
                 J. Differential Geom. 29 (1989), no. 3, 665-683
         \bibitem[Kro2]{kro2} \sameauthor, \textit{A Torelli-type theorem for gravitational instantons}, J. Differential Geom. 29 (1989), no. 3, 685-697.
         \bibitem[LeB]{leb} C. LeBrun, \textit{Complete Ricci-flat Kähler metrics on $\C^n$ need not be flat}, Proc. Symp. Pure Math. 52.2 (1991) 297-304. 
         \bibitem[McK]{mck} J. McKay, \textit{Graphs, singularities, and finite groups}. 
                             The Santa Cruz Conference on Finite Groups (Univ. California, Santa Cruz, Calif., 1979), pp. 183-186, Proc. Sympos. Pure Math., 37, Amer. Math. Soc., Providence, R.I., 1980.
         \bibitem[Min1]{min1} V. Minerbe, \textit{Sur l'effondrement à l'infini des variétés asymptotiquement plates}, 
                 PhD thesis, available at \url{http://www.math.jussieu.fr/~minerbe/theseminerbe.pdf}. 
         \bibitem[Min2]{min2} \sameauthor, \textit{On some asymptotically flat manifolds with non maximal volume growth}, arXiv:0709.1084 [math.DG]
         \bibitem[Min3]{min3} \sameauthor, \textit{Rigidity for Multi-Taub-NUT metrics}, J. Reine Angew. Math. 656 (2011), 47-58.
         \bibitem[{\c{S}}uv]{suv} I. \c{S}uvaina,  \textit{ALE Ricci-flat Kähler metrics and deformations of quotient surface singularities}, 
         \bibitem[TY1]{tian-yau1} G. Tian, S.T. Yau, \textit{Complete Kähler manifolds with zero Ricci curvature. I}.  J. Amer. Math. Soc.  3  (1990),  no. 3, 579-609.
         \bibitem[TY2]{tian-yau2} \sameauthor, \textit{Complete Kähler manifolds with zero Ricci curvature. II}.  Invent. Math. 106 (1991), no. 1, 27-60.
         \bibitem[Yau]{yau} Yau, S.T.: \textit{On the Ricci curvature of a compact Kähler manifold and the complex Monge-Ampère equation}, 
	         $I^*$, Comm. Pure Appl. Math. 31 (1978) 339-441.

   
 \end{thebibliography}
\end{document}